\documentclass[12pt,leqno]{article}
\usepackage{amssymb}
\usepackage[mathscr]{eucal}
\usepackage{amsmath,amssymb,latexsym,theorem,bbm} %
\usepackage{color}
\usepackage{appendix}

\newcommand{\SC}{\scriptstyle}

\newcommand{\CC}{\mathsf{C}}
\newcommand{\DD}{\mathsf{D}}

\newcommand{\NN}{\mathbb{N}}

\newcommand{\RR}{\mathbb{R}}

\newcommand{\ZZ}{\mathbb{Z}}

\newcommand{\bA}{{\boldsymbol{A}}}
\newcommand{\bcA}{{\boldsymbol{\cA}}}

\newcommand{\tbA}{\widetilde{\bA}}

\newcommand{\bB}{{\boldsymbol{B}}}

\newcommand{\bC}{{\boldsymbol{C}}}

\newcommand{\bD}{{\boldsymbol{D}}}

\newcommand{\be}{{\boldsymbol{e}}}

\newcommand{\bI}{{\boldsymbol{I}}}
\newcommand{\bcI}{{\boldsymbol{\cI}}}

\newcommand{\bm}{{\boldsymbol{m}}}
\newcommand{\bM}{{\boldsymbol{M}}}
\newcommand{\bP}{{\boldsymbol{P}}}

\newcommand{\bQ}{{\boldsymbol{Q}}}
\newcommand{\bR}{{\boldsymbol{R}}}
\newcommand{\bu}{{\boldsymbol{u}}}

\newcommand{\bv}{{\boldsymbol{v}}}
\newcommand{\bV}{{\boldsymbol{V}}}
\newcommand{\obV}{{\overline{\bV}}}
\newcommand{\tbV}{\widetilde{\bV}}

\newcommand{\bw}{{\boldsymbol{w}}}

\newcommand{\bx}{{\boldsymbol{x}}}
\newcommand{\bX}{{\boldsymbol{X}}}
\newcommand{\tbX}{\widetilde{\bX}}
\newcommand{\tX}{\widetilde{X}}
\newcommand{\by}{{\boldsymbol{y}}}
\newcommand{\bY}{{\boldsymbol{Y}}}

\newcommand{\bz}{{\boldsymbol{z}}}
\newcommand{\bZ}{{\boldsymbol{Z}}}
\newcommand{\tbZ}{\widetilde{\bZ}}
\newcommand{\bU}{{\boldsymbol{U}}}
\newcommand{\bgamma}{{\boldsymbol{\gamma}}}

\newcommand{\bxi}{{\boldsymbol{\xi}}}

\newcommand{\bzeta}{{\boldsymbol{\zeta}}}
\newcommand{\bvare}{{\boldsymbol{\vare}}}
\newcommand{\bPi}{{\boldsymbol{\Pi}}}

\newcommand{\bzero}{{\boldsymbol{0}}}

\newcommand{\cA}{{\mathcal A}}
\newcommand{\cB}{{\mathcal B}}

\newcommand{\cC}{{\mathcal C}}
\newcommand{\cD}{{\mathcal D}}
\newcommand{\bcD}{\boldsymbol{\cD}}

\newcommand{\cI}{\mathcal{I}}

\newcommand{\cF}{{\mathcal F}}

\newcommand{\cM}{{\mathcal M}}

\newcommand{\bcM}{\boldsymbol{\cM}}
\newcommand{\cN}{{\mathcal N}}
\newcommand{\bcN}{\boldsymbol{\cN}}
\newcommand{\cP}{{\mathcal P}}
\newcommand{\cQ}{{\mathcal Q}}
\newcommand{\bcP}{\boldsymbol{\cP}}
\newcommand{\bcQ}{\boldsymbol{\cQ}}
\newcommand{\cS}{{\mathcal S}}
\newcommand{\cU}{{\mathcal U}}

\newcommand{\bcU}{\boldsymbol{\cU}}

\newcommand{\cX}{{\mathcal X}}
\newcommand{\bcX}{\boldsymbol{\cX}}
\newcommand{\tbcX}{\widetilde{\bcX}}
\newcommand{\cY}{{\mathcal Y}}
\newcommand{\cZ}{{\mathcal Z}}
\newcommand{\cW}{{\mathcal W}}
\newcommand{\bcW}{\boldsymbol{\cW}}
\newcommand{\tbcW}{\widetilde{\bcW}}

\newcommand{\bcY}{\boldsymbol{\cY}}

\newcommand{\bcZ}{\boldsymbol{\cZ}}

\newcommand{\tcW}{\widetilde{\cW}}

\newcommand{\dd}{\mathrm{d}}

\newcommand{\LEFT}{{\mathrm{left}}}
\newcommand{\RIGHT}{\mathrm{right}}
\newcommand{\slu}{{\SC\mathrm{lu}}}

\newcommand{\INARtwo}{\textup{INAR(2)}}
\newcommand{\INARp}{\textup{INAR($p$)}}

\newcommand{\EE}{\operatorname{\mathbb{E}}}
\newcommand{\PP}{\operatorname{\mathbb{P}}}

\newcommand{\OO}{\operatorname{O}}
\newcommand{\var}{\operatorname{Var}}

\newcommand{\Tr}{\operatorname{Tr}}

\newcommand{\halpha}{\widehat{\alpha}}
\newcommand{\hbeta}{\widehat{\beta}}
\newcommand{\hgamma}{\widehat{\gamma}}
\newcommand{\hdelta}{\widehat{\delta}}
\newcommand{\hvarrho}{\widehat{\varrho}}

\newcommand{\hmxi}{\widehat{m}_\xi}
\newcommand{\hbmbxi}{\widehat{\bm}_\bxi}

\newcommand{\tOmega}{\widetilde{\Omega}}

\newcommand{\tV}{\widetilde{V}}

\newcommand{\vare}{\varepsilon}

\renewcommand{\mid}{\,|\,}
\newcommand{\bmid}{\,\big|\,}

\renewcommand{\leq}{\leqslant}
\renewcommand{\geq}{\geqslant}

\newcommand{\stoch}{\stackrel{\PP}{\longrightarrow}}
\newcommand{\distr}{\stackrel{\cD}{\longrightarrow}}
\newcommand{\distre}{\stackrel{\cD}{=}}

\newcommand{\lu}{\stackrel{\slu}{\longrightarrow}}
\newcommand{\as}{\stackrel{{\mathrm{a.s.}}}{\longrightarrow}}
\newcommand{\ase}{\stackrel{{\mathrm{a.s.}}}{=}}

\newcommand{\bbone}{\mathbbm{1}}
\newcommand{\ns}{{\lfloor ns\rfloor}}
\newcommand{\nt}{{\lfloor nt\rfloor}}
\newcommand{\nT}{{\lfloor nT\rfloor}}

\newcommand{\proofend}{\hfill\mbox{$\Box$}}

\numberwithin{equation}{section}

\theoremstyle{change} \theorembodyfont{\em}
\newtheorem{Thm}{Theorem.}[section]
\newtheorem{Lem}[Thm]{Lemma.}
\newtheorem{Pro}[Thm]{Proposition.}
\newtheorem{Cor}[Thm]{Corollary.}

\theorembodyfont{\rm}
\newtheorem{Rem}[Thm]{Remark.}

\begin{document}

\begin{center}
 {\bfseries\Large Statistical inference of 2-type critical \\[2mm]
                   Galton--Watson processes with immigration} \\[5mm]
 {\sc\large Krist\'of $\text{K\"ormendi}^{*,\diamond}$,
            \ Gyula $\text{Pap}^{**}$}
\end{center}

\vskip0.2cm

\noindent
 * MTA-SZTE Analysis and Stochastics Research Group,
   Bolyai Institute, University of Szeged,
   Aradi v\'ertan\'uk tere 1, H--6720 Szeged, Hungary.

\noindent
 ** Bolyai Institute, University of Szeged,
    Aradi v\'ertan\'uk tere 1, H-6720 Szeged, Hungary.

\noindent e-mail: kormendi@math.u-szeged.hu (K. K\"ormendi),
                  papgy@math.u-szeged.hu (G. Pap).

\noindent $\diamond$ Corresponding author.



\renewcommand{\thefootnote}{}
\vspace*{0.2cm}
\footnote{The research of G. Pap was realized in the frames of
 T\'AMOP 4.2.4.\ A/2-11-1-2012-0001 ,,National Excellence Program --
 Elaborating and operating an inland student and researcher personal support
 system''.
The project was subsidized by the European Union and co-financed by the
 European Social Fund.}

\vspace*{-10mm}

\begin{abstract}
In this paper the asymptotic behavior of the conditional least squares
 estimators of the offspring mean matrix for a 2-type critical positively
 regular Galton--Watson branching process with immigration is described.
We also study this question for a natural estimator of the spectral radius of
 the offspring mean matrix, which we call criticality parameter.
We discuss the subcritical case as well.
\end{abstract}

\section{Introduction}
\label{section_intro}

Branching processes have a number of applications in biology, finance,
 economics, queueing theory etc., see e.g.\ Haccou, Jagers and Vatutin
 \cite{HJV}.
Many aspects of applications in epidemiology, genetics and cell kinetics were
 presented at the 2009 Badajoz Workshop on Branching Processes, see
 \cite{Badajoz}.
 
The estimation theory for single-type Galton--Watson branching processes with
 immigration has a long history, see the survey paper of Winnicki \cite{Win1}.
The critical case is the most interesting and complicated one.
There are two multi-type critical Galton--Watson processes with immigration for which
 statistical inference is available: the unstable integer-valued autoregressive models
 of order 2 (which can be considered as a special 2-type Galton--Watson branching
 process with immigration), see Barczy et al.\ \cite{BarIspPap2} and the 2-type doubly
 symmetric process, see Isp\'any et al.\ \cite{IspKorPap}.
In the present paper the asymptotic behavior of the conditional least squares (CLS)
 estimator of the offspring means and criticality parameter for the general 2-type
 critical positively regular Galton--Watson process with immigration is described,
 see Theorem \ref{main}.
It turns out that in a degenerate case this estimator is not even weakly
 consistent.
We also study the asymptotic behavior of a natural estimator of the spectral
 radius of the offspring mean matrix, which we call criticality parameter.
We discuss the subcritical case as well, but the supercritical case still remains
 open.

Let us recall the results for a single-type Galton--Watson branching process
 \ $(X_k)_{k \in \ZZ_+}$ \ with immigration.
Assuming that the immigration mean \ $m_\vare$ \ is known, the CLS estimator
 of the offspring mean \ $m_\xi$ \ based on the observations \ $X_1, \dots, X_n$
 \ has the form
 \[
   \hmxi^{(n)} = \frac{\sum_{k=1}^n X_{k-1} (X_k - m_\vare)}{\sum_{k=1}^n X_{k-1}^2} 
 \]
 on the set \ $\sum_{k=1}^n X_{k-1}^2 > 0$, \ see Klimko ans Nelson \cite{KliNel}.
Suppose that \ $m_\vare > 0$, \ and the second moment of the branching and immigration
 distributions are finite.
 
If the process is subcritical, i.e., \ $m_\xi < 1$, \ then the probability of the
 existence of the estimator \ $\hmxi^{(n)}$ \ tends to 1 as \ $n \to \infty$, \ and
 the estimator \ $\hmxi^{(n)}$ \ is strongly consistent, i.e.,
 \ $\hmxi^{(n)} \as m_\xi$ \ as \ $n \to \infty$.
\ If, in addition, the third moments of the branching and immigration distributions
 are finite, then 
 \begin{equation}\label{m_xi_sub1}
   n^{1/2} (\hmxi^{(n)} - m_\xi)
   \distr
   \cN\biggl(0, \frac{V_\xi \EE(\tX^3) + V_\vare \EE(\tX^2)}
                     {\bigl[\EE(\tX^2)\bigr]^2}\biggr)
   \qquad \text{as \ $n \to \infty$,}
 \end{equation}
 where \ $V_\xi$ \ and \ $V_\vare$ \ denote the offspring and immigration variance,
 respectively, and the distribution of the random variable \ $\tX$ \ is the unique
 stationary distribution of the Markov chain \ $(X_k)_{k \in \ZZ_+}$.
\ Klimko and Nelson \cite{KliNel} contains a similar results for the CLS estimator
 of \ $(m_\xi, m_\vare)$, \ and \eqref{m_xi_sub1} can be derived by the method of
 that paper, see also Theorem \ref{main_sub}.
Note that \ $\EE(\tX^2)$ \ and \ $\EE(\tX^3)$ \ can be expressed by the first three
 moments of the branching and immigration distributions.

If the process is critical, i.e., \ $m_\xi = 1$, \ then the probability of the
 existence of the estimator \ $\hmxi^{(n)}$ \ tends to 1 as \ $n \to \infty$, \ and
 \begin{equation}\label{m_xi_1}
  n (\hmxi^{(n)} - 1)
  \distr
  \frac{\int_0^1 \cX_t \, \dd (\cX_t - m_\vare t)}{\int_0^1 \cX_t^2 \, \dd t}
  \qquad \text{as \ $n \to \infty$,}
 \end{equation}
 where the process \ $(\cX_t)_{t\in\RR_+}$ \ is the unique strong solution of the
 stochastic differential equation (SDE)
 \[
   \dd \cX_t = m_\vare \, \dd t + \sqrt{ V_\xi \cX_t^+ } \, \dd \cW_t ,
   \qquad t \in \RR_+ ,
 \]
 with initial value \ $\cX_0 = 0$, \ where \ $(\cW_t)_{t\in\RR_+}$ \ is a
 standard Wiener process, and \ $x^+$ \ denotes the positive part of \ $x \in \RR$.
\ Wei and Winnicki \cite{WW2} proved a similar results for the CLS estimator of the
 offspring mean when the immigration mean is unknown, and \eqref{m_xi_1} can be
 derived by the method of that paper.
Note that \ $\cX^{(n)} \distr \cX$ \ as \ $n \to \infty$ \ with
 \ $\cX^{(n)}_t := n^{-1} X_{\nt}$ \ for \ $t \in \RR_+$, \ $n \in \NN$, \ where
 \ $\lfloor x \rfloor$ \ denotes the (lower) integer part of \ $x \in \RR$, \ see Wei
 and Winnicki \cite{WW}.
If, in addition, \ $V_\xi = 0$, \ then
 \begin{equation}\label{m_xi_1_0}
  n^{3/2} (\hmxi^{(n)} - 1)
  \distr
  \cN\biggl(0, \frac{3 V_\vare}{m_\vare^2}\biggr)
  \qquad \text{as \ $n \to \infty$,}
 \end{equation}
 see Isp\'any et al.\ \cite{IspPapZui1}.

If the process is supercritical, i.e., \ $m_\xi > 1$, \ then the probability of the
 existence of the estimator \ $\hmxi^{(n)}$ \ tends to 1 as \ $n \to \infty$, \ the
 estimator \ $\hmxi^{(n)}$ \ is strongly consistent, i.e., \ $\hmxi^{(n)} \as m_\xi$
 \ as \ $n \to \infty$, \ and 
 \begin{equation}\label{m_xi_super1r}
  \biggl(\sum_{k=1}^n X_{k-1}\biggr)^{1/2} (\hmxi^{(n)} - m_\xi)
  \distr
  \cN\biggl(0, \frac{(m_\xi+1)^2}{m_\xi^2+m_\xi+1} \, V_\xi\biggr)
  \qquad \text{as \ $n \to \infty$.}
 \end{equation}
Wei and Winnicki \cite{WW2} showed the same asymptotic behavior for the CLS estimator
 of the offspring mean when the immigration mean is unknown, and \eqref{m_xi_super1r}
 can be derived by the method of that paper.

In Section \ref{section_preliminaries} we recall some preliminaries on
 2-type Galton--Watson models with immigration.
Section \ref{section_main_results} contains our main results.
Section \ref{section_deco_proof} contains a useful decomposition of the
 process.
Sections  \ref{section_proof_main}, \ref{section_proof_maint} and 
 \ref{section_proof_sub} contain the proofs. 
In Appendix \ref{section_moments} we present estimates for the moments of the
 processes involved.
Appendix \ref{section_estimators} is devoted to the CLS estimators.
Appendix \ref{app_C} and \ref{section_conv_step_processes} is for a version of
 the continuous mapping theorem and for convergence of random step processes,
 respectively.

\section{Preliminaries on 2-type Galton--Watson models
            with immigration}
\label{section_preliminaries}

Let \ $\ZZ_+$, \ $\NN$, \ $\RR$ \ and \ $\RR_+$ \ denote the set of
 non-negative integers, positive integers, real numbers and non-negative real
 numbers, respectively.
Every random variable will be defined on a fixed probability space
 \ $(\Omega,\cA,\PP)$.

For each \ $k, j \in \ZZ_+$ \ and \ $i, \ell \in \{ 1, 2 \}$, \ the number of
 individuals of type \ $i$ \ in the \ $k^\mathrm{th}$ \ generation will be
 denoted by \ $X_{k,i}$, \ the number of type \ $\ell$ \ offsprings produced by
 the \ $j^\mathrm{th}$ \ individual who is of type \ $i$ \ belonging to the
 \ $(k-1)^\mathrm{th}$ \ generation will be denoted by \ $\xi_{k,j,i,\ell}$, \ and
 the number of type \ $i$ \ immigrants in the \ $k^\mathrm{th}$ \ generation 
 will be denoted by \ $\vare_{k,i}$.
\ Then we have
 \begin{equation}\label{GWI(2)}
  \begin{bmatrix}
   X_{k,1} \\
   X_{k,2}
  \end{bmatrix}
  = \sum_{j=1}^{X_{k-1,1}}
     \begin{bmatrix}
      \xi_{k,j,1,1} \\
      \xi_{k,j,1,2}
     \end{bmatrix}
    +\sum_{j=1}^{X_{k-1,2}}
       \begin{bmatrix}
        \xi_{k,j,2,1} \\
        \xi_{k,j,2,2}
       \end{bmatrix}
    + \begin{bmatrix}
       \vare_{k,1} \\
       \vare_{k,2}
      \end{bmatrix} ,
  \qquad k \in \NN .
 \end{equation}
Here
 \ $\bigl\{ \bX_0, \, \bxi_{k,j,i}, \, \bvare_k
            : k, j \in \NN, \, i \in \{ 1, 2 \} \bigr\}$
 \ are supposed to be independent, where 
 \[
   \bX_k := \begin{bmatrix}
             X_{k,1} \\
             X_{k,2}
            \end{bmatrix} , \qquad
   \bxi_{k,j,i} := \begin{bmatrix}
                  \xi_{k,j,i,1} \\
                  \xi_{k,j,i,2}
                 \end{bmatrix} , \qquad
   \bvare_k := \begin{bmatrix}
                \vare_{k,1} \\
                \vare_{k,2}
               \end{bmatrix} .
 \]
Moreover, \ $\{ \bxi_{k,j,1} : k, j \in \NN \}$,
 \ $\{ \bxi_{k,j,2} : k, j \in \NN \}$ \ and \ $\{ \bvare_k : k \in \NN \}$
 \ are supposed to consist of identically distributed random vectors.

We suppose \ $\EE(\|\bxi_{1,1,1}\|^2) < \infty$,
 \ $\EE(\|\bxi_{1,1,2}\|^2) < \infty$ \ and \ $\EE(\|\bvare_1\|^2) < \infty$.
\ Introduce the notations
 \begin{gather*}
  \bm_{\bxi_i} := \EE\bigl(\bxi_{1,1,i}\bigr) \in \RR^2_+ , \qquad
  \bm_{\bxi} := \begin{bmatrix}
                \bm_{\bxi_1} & \bm_{\bxi_2}
               \end{bmatrix} \in \RR^{2 \times 2}_+ , \qquad
  \bm_{\bvare} := \EE\bigl(\bvare_1\bigr) \in \RR^2_+ , \\
  \bV_{\!\!\bxi_i} := \var\bigl(\bxi_{1,1,i}\bigr) \in \RR^{2 \times 2} ,
  \qquad
  \bV_{\!\!\bvare} := \var\bigl(\bvare_1\bigr) \in \RR^{2 \times 2}  , \qquad
  i \in \{1, 2\} .
 \end{gather*}
Note that many authors define the offspring mean matrix as \ $\bm^\top_\bxi$. 
\ For \ $k \in \ZZ_+$, \ let
 \ $\cF_k := \sigma\bigl( \bX_0, \bX_1 , \dots, \bX_k \bigr)$.
\ By \eqref{GWI(2)}, 
 \begin{equation}\label{mart}
  \EE(\bX_k \mid \cF_{k-1}) 
  = X_{k-1,1} \, \bm_{\bxi_1} + X_{k-1,2} \, \bm_{\bxi_2} +\bm_{\bvare}
  = \bm_{\bxi} \, \bX_{k-1} + \bm_{\bvare} .
 \end{equation}
Consequently,
 \[
   \EE(\bX_k) = \bm_{\bxi} \EE(\bX_{k-1}) + \bm_{\bvare} , \qquad k \in \NN ,
 \]
 which implies
 \begin{equation}\label{EXk}
  \EE(\bX_k) 
  = \bm_{\bxi}^k \, \EE(\bX_0)
    + \sum_{j=0}^{k-1} \bm_{\bxi}^j \, \bm_{\bvare} , \qquad
  k \in \NN .
 \end{equation}
Hence, the asymptotic behavior of the sequence \ $(\EE(\bX_k))_{ k \in \ZZ_+ }$
 \ depends on the asymptotic behavior of the powers \ $(\bm_{\bxi}^k)_{k \in \NN}$
 \ of the offspring mean matrix, which is related to the spectral radius
 \ $r(\bm_\bxi) =: \varrho \in \RR_+$ \ of \ $\bm_\bxi$ \ (see the 
 Frobenius--Perron theorem, e.g., Horn and Johnson
 \cite[Theorems 8.2.11 and 8.5.1]{HJ}).
A 2-type Galton--Watson process \ $(\bX_k)_{k \in \ZZ_+}$ \ with immigration
 is referred to respectively as \emph{subcritical}, \emph{critical} or
 \emph{supercritical} if \ $\varrho < 1$, \ $\varrho = 1$ \ or \ $\varrho > 1$
 \ (see, e.g., Athreya and Ney \cite[V.3]{AN} or Quine \cite{Q}).
We will write the offspring mean matrix of a 2-type Galton--Watson process
 with immigration in the form
 \begin{equation}\label{bA}
  \bm_\bxi :=\begin{bmatrix}
             \alpha & \beta \\
             \gamma  & \delta 
            \end{bmatrix} .
 \end{equation}
We will focus only on \emph{positively regular} 2-type Galton--Watson
 processes with immigration, i.e., when there is a positive integer
 \ $k \in \NN$ \ such that the entries of \ $\bm_\bxi^k$ \ are positive (see
 Kesten and Stigum \cite{KesSti1}), which is equivalent to
 \ $\beta, \gamma \in (0, \infty)$, \ $\alpha, \delta \in \RR_+$ \ with
 \ $\alpha + \delta > 0$.
\ Then the matrix \ $\bm_\bxi$ \ has eigenvalues
 \[
   \lambda_+
   := \frac{\alpha + \delta + \sqrt{(\alpha - \delta)^2 + 4 \beta \gamma}}{2} ,
   \qquad
   \lambda_-
   := \frac{\alpha + \delta - \sqrt{(\alpha - \delta)^2 + 4 \beta \gamma}}{2} ,
 \]
 satisfying \ $\lambda_+ > 0$ \ and \ $-\lambda_+ < \lambda_- < \lambda_+$,
 \ hence the spectral radius of \ $\bm_\bxi$ \ is
 \begin{equation}\label{varrho}
  \varrho = r(\bm_\bxi) = \lambda_+
  = \frac{\alpha + \delta + \sqrt{(\alpha - \delta)^2 + 4 \beta \gamma}}{2} .
 \end{equation}
By the Perron theorem (see, e.g., Horn and Johnson
 \cite[Theorems 8.2.11 and 8.5.1]{HJ}), 
 \[
   \lambda_+^{-k} \bm_{\bxi}^k \to \bu_\RIGHT \bu_\LEFT^\top \qquad
   \text{as \ $k \to \infty$,}
 \]
 where \ $\bu_\RIGHT$ \ is the unique right eigenvector of \ $\bm_\bxi$
 \ (called the right Perron vector of \ $\bm_\bxi$) corresponding to the
 eigenvalue \ $\lambda_+$ \ such that the sum of its coordinates is 1, and
 \ $\bu_\LEFT$ \ is the unique left eigenvector of \ $\bm_\bxi$ \ (called the
 left Perron vector of \ $\bm_\bxi$) corresponding to the eigenvalue
 \ $\lambda_+$ \ such that \ $\langle \bu_\RIGHT, \bu_\LEFT \rangle = 1$,
 \ hence we have
 \[
   \bu_\RIGHT
   = \frac{1}{\beta+\lambda_+-\alpha}
     \begin{bmatrix} \beta \\ \lambda_+-\alpha \end{bmatrix} , \qquad
   \bu_\LEFT
   = \frac{1}{\lambda_+-\lambda_-}
     \begin{bmatrix}
      \gamma + \lambda_+ - \delta \\ 
      \beta + \lambda_+ - \alpha
     \end{bmatrix} .
 \]
More exactly, using the so-called Putzer's spectral formula, see, e.g., Putzer
 \cite{Put}, the powers of \ $\bm_\bxi$ \ can be written in the form
 \begin{equation}\label{Putzer}
   \begin{aligned}
    \bm_\bxi^k
    &= \frac{\lambda_+^k}{\lambda_+-\lambda_-}
       \begin{bmatrix}
        \lambda_+ - \delta & \beta \\
        \gamma & \lambda_+ - \alpha
       \end{bmatrix}
     + \frac{\lambda_-^k}{\lambda_+-\lambda_-}
       \begin{bmatrix}
        \lambda_+ - \alpha & - \beta \\
        - \gamma & \lambda_+ - \delta
       \end{bmatrix} \\
  &= \lambda_+^k \bu_\RIGHT \bu_\LEFT^\top 
     + \lambda_-^k \bv_\RIGHT \bv_\LEFT^\top ,
   \qquad k \in \NN ,
  \end{aligned}
 \end{equation}
 where \ $\bv_\RIGHT$ \ and \ $\bv_\LEFT$ \ are appropriate right and left
 eigenvectors of \ $\bm_\bxi$, \ respectively, belonging to the eigenvalue
 \ $\lambda$, \ for instance,
 \[
   \bv_\RIGHT
   = \frac{1}{\lambda_+-\lambda_-}
     \begin{bmatrix}
      -\beta - \lambda_+ + \alpha \\
      \gamma + \lambda_+ - \delta 
     \end{bmatrix} ,
   \qquad
   \bv_\LEFT
   = \frac{1}{\beta+\lambda_+-\alpha}
     \begin{bmatrix} - \lambda_+ + \alpha \\ \beta \end{bmatrix} .
 \]
The process \ $(\bX_k)_{k \in \ZZ_+}$ \ is critical and positively regular if and
 only if \ $\alpha, \delta \in [0, 1)$ \ and \ $\beta, \gamma \in (0, \infty)$
 \ with \ $\alpha + \delta > 0$ \ and
 \ $\beta \gamma = (1 - \alpha) (1 - \delta)$, \ and then the matrix \ $\bm_\bxi$
 \ has eigenvalues \ $\lambda_+ = 1$ \ and
 \[
   \lambda_- = \alpha + \delta - 1 \in (-1, 1) =: \lambda .
 \]
Next we will recall a convergence result for critical and positively regular 
 2-type CBI processes.
A function \ $f : \RR_+ \to \RR^d$ \ is called \emph{c\`adl\`ag} if it is right
 continuous with left limits.
\ Let \ $\DD(\RR_+, \RR^d)$ \ and \ $\CC(\RR_+, \RR^d)$ \ denote the space of
 all \ $\RR^d$-valued c\`adl\`ag and continuous functions on \ $\RR_+$,
 \ respectively.
Let \ $\cD_\infty(\RR_+, \RR^d)$ \ denote the Borel $\sigma$-field in
 \ $\DD(\RR_+, \RR^d)$ \ for the metric characterized by Jacod and Shiryaev
 \cite[VI.1.15]{JacShi} (with this metric \ $\DD(\RR_+, \RR^d)$ \ is a complete
 and separable metric space).
For \ $\RR^d$-valued stochastic processes \ $(\bcY_t)_{t\in\RR_+}$ \ and
 \ $(\bcY^{(n)}_t)_{t\in\RR_+}$, \ $n \in \NN$, \ with c\`adl\`ag paths we write
 \ $\bcY^{(n)} \distr \bcY$ \ as \ $n \to \infty$ \ if the distribution of
 \ $\bcY^{(n)}$ \ on the space \ $(\DD(\RR_+, \RR^d), \cD_\infty(\RR_+, \RR^d))$
 \ converges weakly to the distribution of \ $\bcY$ \ on the space
 \ $(\DD(\RR_+, \RR^d), \cD_\infty(\RR_+, \RR^d))$ \ as \ $n \to \infty$.
\ Concerning the notation \ $\distr$ \ we note that if \ $\xi$ \ and \ $\xi_n$,
 \ $n \in \NN$, \ are random elements with values in a metric space
 \ $(E, \rho)$, \ then we also denote by \ $\xi_n \distr \xi$ \ the weak
 convergence of the distributions of \ $\xi_n$ \ on the space \ $(E, \cB(E))$
 \ towards the distribution of \ $\xi$ \ on the space \ $(E, \cB(E))$ \ as
 \ $n \to \infty$, \ where \ $\cB(E)$ \ denotes the Borel $\sigma$-algebra on
 \ $E$ \ induced by the given metric \ $\rho$.

For each \ $n \in \NN$, \ consider the random step process
 \[
   \bcX^{(n)}_t := n^{-1} \bX_\nt , \qquad t \in \RR_+ .
 \]
The following theorem is a special case of the main result in Isp\'any and Pap
 \cite[Theorem 3.1]{IspPap2}.

\begin{Thm}\label{conv}
Let \ $(\bX_k)_{k\in\ZZ_+}$ \ be a 2-type Galton--Watson process with immigration 
 such that \ $\alpha, \delta \in [0, 1)$ \ and \ $\beta, \gamma \in (0, \infty)$
 \ with \ $\alpha + \delta > 0$ \ and \ $\beta \gamma = (1 - \alpha) (1 - \delta)$
 \ (hence it is critical and positively regular), \ $\bX_0 = \bzero$,
 \ $\EE(\|\bxi_{1,1,1}\|^2) < \infty$, \ $\EE(\|\bxi_{1,1,2}\|^2) < \infty$ \ and
 \ $\EE(\|\bvare_1\|^2) < \infty$.
\ Then
 \begin{equation}\label{Conv_X}
  (\bcX_t^{(n)})_{t\in\RR_+}
  \distr (\bcX_t)_{t\in\RR_+} := (\cZ_t \bu_\RIGHT)_{t\in\RR_+} \qquad
  \text{as \ $n \to \infty$}
 \end{equation}
 in \ $\DD(\RR_+, \RR^d)$, \ where \ $(\cZ_t)_{t \in \RR_+}$ \ is the pathwise
 unique strong solution of the SDE
 \begin{equation}\label{SDE_Y}
  \dd \cZ_t
  = \langle \bu_\LEFT, \bm_\bvare \rangle \, \dd t
    + \sqrt{ \langle \obV_{\!\!\bxi} \bu_\LEFT, \bu_\LEFT \rangle \cZ_t^+ }
      \, \dd \cW_t ,
  \qquad t \in \RR_+ , \qquad \cZ_0 = 0 ,
 \end{equation}
 where \ $(\cW_t)_{t \in \RR_+}$ \ is a standard Brownian motion and
 \begin{equation}\label{obV}
  \obV_{\!\!\bxi}
  := \sum_{i=1}^2 \langle \be_i, \bu_\RIGHT \rangle \bV_{\!\!\bxi_i}
  = \frac{\beta \bV_{\!\!\bxi_1} + (1 - \alpha) \bV_{\!\!\bxi_2}}
          {\beta + 1 - \alpha}
 \end{equation}
 is a mixed offspring variance matrix.
\end{Thm}

In fact, in Isp\'any and Pap \cite[Theorem 3.1]{IspPap2}, the above result has been 
 prooved under the higher moment assumptions
 \ $\EE(\|\bxi_{1,1,1}\|^4) < \infty$, \ $\EE(\|\bxi_{1,1,2}\|^4) < \infty$ \ and
 \ $\EE(\|\bvare_1\|^4) < \infty$, \ which have been relaxed in Danka and Pap 
 \cite[Theorem 3.1]{DanPap}.

\begin{Rem}\label{REMARK_SDE}
The SDE \eqref{Y} has a unique strong solution \ $(\cZ_t^{(z)})_{t \in \RR_+}$
 \ for all initial values  \ $\cZ_0^{(z)} = z \in \RR$, \ and if \ $z \geq 0$,
 \ then \ $\cZ_t^{(z)}$ \ is nonnegative for all \ $t \in \RR_+$ \ with
 probability one, hence \ $\cZ_t^+$ \ may be replaced by \ $\cZ_t$ \ under the
 square root in \eqref{Y}, see, e.g.,
 Barczy et al.~\cite[Remark 3.3]{BarIspPap0}.
\end{Rem}

Clearly, \ $\obV_{\!\!\bxi}$ \ depends only on the branching distributions,
 i.e., on the distributions of \ $\bxi_{1,1,1}$ \ and \ $\bxi_{1,1,2}$.
\ Note that \ $\obV_{\!\!\bxi}= \var(\bY_1 \mid \bY_0 = \bu_\RIGHT)$, \ where
 \ $(\bY_k)_{k\in\ZZ_+}$ \ is a 2-type Galton--Watson process without immigration
 such that its branching distributions are the same as that of
 \ $(\bX_t)_{k\in\ZZ_+}$, \ since for each \ $i \in \{1, 2\}$,
 \ $\bV_{\!\!\bxi_i} = \var(\bY_1 \mid \bY_0 = \be_i)$.

For the sake of simplicity, we consider a zero start Galton--Watson process
 with immigration, that is, we suppose \ $\bX_0 = \bzero$.
\ The general case of nonzero initial value may be handled in a similar way,
 but we renounce to consider it.
In the sequel we always assume \ $\bm_\vare \ne \bzero$, \ otherwise
 \ $\bX_k = \bzero$ \ for all \ $k \in \NN$.

\section{Main results}
\label{section_main_results}

For each \ $n \in \NN$, \ any CLS estimator
 \[
   \hbmbxi^{(n)} = \begin{bmatrix}
                   \halpha_n & \hbeta_n \\
                   \hgamma_n & \hdelta_n
                  \end{bmatrix}
 \]
 of the offspring mean matrix \ $\bm_\bxi$ \ based on a sample
 \ $\bX_1, \ldots, \bX_n$ \ has the form
 \begin{equation}\label{CLSEabc}
  \hbmbxi^{(n)} = \bB_n \bA_n^{-1} 
 \end{equation}
 on the set
 \begin{equation}\label{H_n}
  \Omega_n := \left\{ \omega \in \Omega : \det(\bA_n(\omega)) > 0 \right\} ,
 \end{equation}
 where
 \begin{align}\label{A_n_B_n}
  \bA_n := \sum_{k=1}^n \bX_{k-1} \bX_{k-1}^\top , \qquad
  \bB_n := \sum_{k=1}^n (\bX_k - \bm_\bvare) \bX_{k-1}^\top ,
 \end{align}
 see Lemma \ref{CLSE1}.
The spectral radius \ $\varrho$ \ given in \eqref{varrho} can be called
 \emph{criticality parameter}, and its natural estimator is the spectral radius
 of \ $\hbmbxi^{(n)}$, \ namely,
 \begin{equation}\label{hvarrho}
  \hvarrho_n := r\bigl(\hbmbxi^{(n)}\bigr)
  = \frac{\halpha_n + \hdelta_n
           + \sqrt{(\halpha_n - \hdelta_n)^2 + 4 \hbeta_n \hgamma_n}}
          {2} ,
 \end{equation}
 on the set \ $\Omega_n \cap \tOmega_n$ \ with
 \begin{equation}\label{tH_n}
  \tOmega_n
  := \bigl\{ \omega \in \Omega_n
             : (\halpha_n(\omega) - \hdelta_n(\omega))^2
               + 4 \hbeta_n(\omega) \hgamma_n(\omega) \geq 0 \bigr\} .
 \end{equation}
 
First we consider the critical and positively regular case.
By Lemma \ref{ExUn}, \ $\PP(\Omega_n) \to 1$ \ and \ $\PP(\tOmega_n) \to 1$
 \ as \ $n \to \infty$ \ under appropriate assumptions.

\begin{Thm}\label{main}
Let \ $(\bX_k)_{k \in \ZZ_+}$ \ be a 2-type Galton--Watson process with immigration
 such that \ $\alpha, \delta \in [0, 1)$ \ and \ $\beta, \gamma \in (0, \infty)$
 \ with \ $\alpha + \delta > 0$ \ and \ $\beta \gamma = (1 - \alpha) (1 - \delta)$
 \ (hence it is critical and positively regular), \ $\bX_0 = \bzero$,
 \ $\EE(\|\bxi_{1,1,1}\|^8) < \infty$, \ $\EE(\|\bxi_{1,1,2}\|^8) < \infty$,
 \ $\EE(\|\bvare_1\|^8) < \infty$, \ and \ $\bm_\bvare \ne \bzero$.

If
 \ $\langle \obV_{\!\!\bxi} \bv_\LEFT, \bv_\LEFT \rangle
    + \langle \bV_{\!\!\bvare} \bv_\LEFT, \bv_\LEFT \rangle
    + \langle \bv_\LEFT, \bm_\bvare \rangle^2 > 0$,
 \ then the probability of the existence of the estimator \ $\hbmbxi^{(n)}$
 \ tends to 1 as \ $n \to \infty$.

If
 \ $\langle \obV_{\!\!\bxi} \bv_\LEFT, \bv_\LEFT \rangle > 0$,
 \ then the probability of the existence of the estimator \ $\hvarrho_n$
 \ tends to 1 as \ $n \to \infty$.

If \ $\langle \obV_{\!\!\bxi} \bv_\LEFT, \bv_\LEFT \rangle > 0$, \ then
 \begin{align}\label{m_xi}
  n^{1/2} (\hbmbxi^{(n)}- \bm_\bxi)
  &\distr
   \frac{(1 - \lambda^2)^{1/2}}
        {\langle \obV_{\!\!\bxi} \bv_\LEFT, \bv_\LEFT \rangle^{1/2}} \,
   \frac{\obV_{\!\!\bxi}^{1/2} \int_0^1 \cY_t \, \dd\tbcW_t}
        {\int_0^1 \cY_t \, \dd t}
   \bv_\LEFT^\top , \\
  \label{varrho_n}
  n (\hvarrho_n - 1)
  &\distr
   \frac{\int_0^1 \cY_t \, \dd(\cY_t - t \langle \bu_\LEFT, \bm_\bvare \rangle)}
        {\int_0^1 \cY_t^2 \, \dd t} ,
 \end{align}
 as \ $n \to \infty$, \ with
 \ $\cY_t := \langle \bu_\LEFT, \bcM_t + t \bm_\bvare \rangle$, \ $t \in \RR_+$,
 \ where \ $(\bcM_t)_{t\in\RR_+}$ \ is the unique strong solution of the SDE
 \begin{equation}\label{Y}
  \dd\bcM_t
  = (\langle \bu_\LEFT, \bcM_t + t \bm_\bvare \rangle^+)^{1/2} \, \obV_{\!\!\bxi}^{1/2}
    \, \dd\bcW_t ,
  \qquad t \in \RR_+ \qquad \bcM_0 = \bzero ,
 \end{equation}
 where \ $(\bcW_t)_{t \in \RR_+}$ \ and \ $(\tbcW_t)_{t \in \RR_+}$ \ are independent
 2-dimenional standard Wiener processes.

If \ $\langle \obV_{\!\!\bxi} \bv_\LEFT, \bv_\LEFT \rangle = 0$ \ and
 \ $\langle \bV_{\!\!\bvare} \bv_\LEFT, \bv_\LEFT \rangle
    + \langle \bv_\LEFT, \bm_\bvare \rangle^2 > 0$,
 \ then
 \begin{equation}\label{m_xi_m}
  \hbmbxi^{(n)} - \bm_\bxi
  \distr
  \frac{\bcI_3 + \bcI_4}{\cI_1 + \cI_2}
  \bv_\LEFT^\top \qquad \text{as \ $n \to \infty$,}
 \end{equation}
 where
 \begin{align*}
  \cI_1 &:= \frac{\langle \bv_\LEFT, \bm_\bvare \rangle^2}
                 {(1 - \lambda)^2}
            \left[ \int_0^1 \cY_t^2 \, \dd t
                   - \left(\int_0^1 \cY_t \, \dd t\right)^2 \right] , \\ 
  \cI_2 &:= \frac{\langle \bV_{\!\!\bvare} \bv_\LEFT, \bv_\LEFT \rangle}
                 {1 - \lambda^2}
            \int_0^1 \cY_t^2 \, \dd t , \\
  \bcI_3 &:= \frac{\langle \bv_\LEFT, \bm_\bvare \rangle}{1 - \lambda}
            \left( \int_0^1 \cY_t^2 \, \dd t \int_0^1 1 \, \dd \bcM_t
                   - \int_0^1 \cY_t \, \dd t
                     \int_0^1 \cY_t \, \dd \bcM_t \right) , \\  
  \bcI_4 &:= \frac{\langle \bV_{\!\!\bvare} \bv_\LEFT, \bv_\LEFT \rangle^{1/2}}
                       {(1 - \lambda^2)^{1/2}}
            \int_0^1 \cY_t^2 \, \dd t \,
            \obV_{\!\!\bxi}^{1/2} \!\!
            \int_0^1 \cY_t \, \dd\tbcW_t .
 \end{align*}
\end{Thm}

\begin{Rem}\label{REMARK00}
If
 \ $\langle \obV_{\!\!\bxi} \bv_\LEFT, \bv_\LEFT \rangle
    + \langle \bV_{\!\!\bvare} \bv_\LEFT, \bv_\LEFT \rangle
    + \langle \bv_\LEFT, \bm_\bvare \rangle^2 = 0$
 \ then, by Lemma \ref{main_VVt},
 \ $(1 - \alpha) X_{k,1} \ase \beta X_{k,2}$ \ for all \ $k \in \NN$, \ hence
 there is no unique CLS estimator for the offspring mean matrix \ $\bm_\bxi$.
\proofend
\end{Rem}

\begin{Rem}\label{REMARK1} 
By It\^o's formula, \eqref{Y} yields that \ $(\cY_t)_{t \in \RR_+}$ \ satisfies
 the SDE \eqref{SDE_Y} with initial value \ $\cY_0 = 0$.
\ Indeed, by It\^o's formula and the first 2-dimensional equation of the SDE
 \eqref{MNPSDE} we obtain
 \[
   \dd \cY_t
   = \langle \bu_\LEFT, \bm_\bvare \rangle \, \dd t
     + (\cY_t^+)^{1/2} \bu_\LEFT^\top \obV_{\!\!\bxi}^{1/2} \, \dd \bcW_t , \qquad
   t \in \RR_+ .
 \]
If
 \ $\langle \obV_{\!\!\bxi} \bu_\LEFT, \bu_\LEFT \rangle
    = \|\bu_\LEFT^\top \obV_{\!\!\bxi}^{1/2}\|^2 = 0$
 \ then \ $\bu_\LEFT^\top \obV_{\!\!\bxi}^{1/2} = \bzero$, \ hence
 \ $\dd \cY_t = \langle \bu_\LEFT, \bm_\bvare \rangle \dd t$, \ $t \in \RR_+$,
 \ implying that the process \ $(\cY_t)_{t\in\RR_+}$ \ satisfies the SDE
 \eqref{SDE_Y}.
If \ $\langle \obV_{\!\!\bxi} \bu_\LEFT, \bu_\LEFT \rangle \ne 0$ \ then the process
 \[
   \tcW_t := \frac{\langle \obV_{\!\!\bxi}^{1/2} \bu_\LEFT, \bcW_t \rangle}
                  {\langle \obV_{\!\!\bxi} \bu_\LEFT, \bu_\LEFT \rangle^{1/2}} ,
   \qquad t \in \RR_+ ,
 \]
 is a (one-dimensional) standard Wiener process, hence the process
 \ $(\cY_t)_{t\in\RR_+}$ \ satisfies the SDE \eqref{SDE_Y}.

Consequently, \ $(\cY_t)_{t\in\RR_+} \distre (\cZ_t)_{t\in\RR_+}$, \ where
 \ $(\cZ_t)_{t \in \RR_+}$ \ is the unique strong solution of the SDE
 \eqref{SDE_Y} with initial value \ $\cZ_0 = 0$, \ hence, by Theorem
 \ref{conv},
 \begin{equation}\label{convX}
  (\bcX^{(n)}_t)_{t\in\RR_+}
  \distr (\bcX_t)_{t\in\RR_+}
  \distre (\cY_t \bu_\RIGHT)_{t\in\RR_+}
  \qquad \text{as \ $n \to \infty$.}
 \end{equation}
 
If \ $\langle \obV_{\!\!\bxi} \, \bu_\LEFT, \bu_\LEFT \rangle = 0$ \ then the
 unique strong solution of \eqref{SDE_Y} is the deterministic function
 \ $\cY_t = \langle \bu_\LEFT, \bm_\bvare \rangle \, t$, \ $t \in \RR_+$,
 \ and then, in \eqref{m_xi}, we have eventually asymptotic normality, since
 \[
   \frac{\int_0^1 \cY_t \obV_{\!\!\bxi}^{1/2} \, \dd\tbcW_t}
        {\int_0^1 \cY_t \, \dd t}
    \distre
    \cN\left(0, \frac{4}{3} \obV_{\!\!\bxi} \right) ,
 \]
 and by \eqref{varrho_n}, \ $n (\hvarrho_n- 1) \distr 0$, \ i.e., the scaling
 \ $n$ \ is not proper.
\proofend
\end{Rem}

\begin{Rem}\label{REMARK0}
Note that \ $\langle \obV_{\!\!\bxi} \bu_\LEFT, \bu_\LEFT \rangle = 0$ \ is
 fulfilled if and only if 
 \ $\langle \bu_\LEFT, \bxi_{k,j,i} \rangle
    \ase \langle \bu_\LEFT, \EE(\bxi_{1,1,i}) \rangle$
 \ for each \ $k, j \in \NN$ \ and \ $i \in \{1, 2\}$, \ i.e., the offspring
 point \ $ \bxi_{k,j,i}$ \ lies almost surely on a line, orthogonal to
 \ $\bu_\LEFT$ \ and containing the point \ $\EE(\bxi_{1,1,i})$.
\ Indeed, \ $\langle \obV_{\!\!\bxi} \bu_\LEFT, \bu_\LEFT \rangle = 0$
 \ is fulfilled if and only if
 \ $\langle \bV_{\!\!\bxi_i} \bu_\LEFT, \bu_\LEFT \rangle
    = \EE[\langle \bu_\LEFT, \bxi_{k,j,i} - \EE(\bxi_{k,j,i}) \rangle^2]
    = 0$
 \ for each \ $k, j \in \NN$ \ and \ $i \in \{1, 2\}$.

In a similar way, \ $\langle \obV_{\!\!\bxi} \bv_\LEFT, \bv_\LEFT \rangle = 0$
 \ is fulfilled if and only if
 \ $\langle \bv_\LEFT, \bxi_{k,j,i} \rangle
    \ase \langle \bv_\LEFT, \EE(\bxi_{1,1,i}) \rangle$
 \ for each \ $k, j \in \NN$ \ and \ $i \in \{1, 2\}$, \ i.e., the offspring
 point \ $ \bxi_{k,j,i}$ \ lies almost surely on a line, orthogonal to
 \ $\bv_\LEFT$ \ and containing the point \ $\EE(\bxi_{1,1,i})$.
\proofend
\end{Rem}

\begin{Rem}\label{REMARK3}
We note that in the critical positively regular case the limit distribution
 for the CLS estimator of the offspring mean matrix \ $\bm_\bxi$ \ is
 concentrated on the 2-dimensional subspace
 \ $\RR^2 \bv_\LEFT^\top \subset \RR^{2\times2}$.
\ Surprisingly, the scaling factor of the CLS estimators of \ $\bm_\bxi$ \ is
 \ $\sqrt{n}$, \ which is the same as in the subcritical case.
The reason of this strange phenomenon can be understood from the joint
 asymptotic behavior of \ $\det(\bA_n)$ \ and \ $\bD_n \tbA_n$ \ given in
 Theorems \ref{main_Ad} and \ref{maint_Ad}.
One of the decisive tools in deriving the needed asymptotic behavior is a good
 bound for the moments of the involved processes, see Corollary
 \ref{EEX_EEU_EEV}.
\proofend
\end{Rem}

\begin{Rem}\label{REMARK4}
The shape of
 \ $\int_0^1 \cY_t \, \dd\bigl(\obV_{\!\!\bxi}^{1/2} \tbcW_t \bv_\LEFT^\top\bigr)
    \big/ \int_0^1 \cY_t \, \dd t$
 \ in \eqref{m_xi} is similar to the limit distribution of the Dickey--Fuller
 statistics for unit root test of AR(1) time series, see, e.g.,
 Hamilton \cite[17.4.2 and 17.4.7]{Ham} or Tanaka
 \cite[(7.14) and Theorem 9.5.1]{Tan}, but it contains two independent
 2-dimensional standard Wiener processes, which can be reduced to three
 1-dimensional independent standard Wiener processes.
This phenomenon is very similar to the appearance of two independent standard
 Wiener processes in limit theorems for CLS estimators of the variance of the
 offspring and immigration distributions for critical branching processes with
 immigration in Winnicki \cite[Theorems 3.5 and 3.8]{Win2}.
Finally, note that the limit distribution of the CLS estimator of the
 offspring mean matrix \ $\bm_\bxi$ \ is always symmetric, although non-normal
 in \eqref{m_xi}.
Indeed, since \ $(\bcW_t)_{t \in \RR_+}$ \ and \ $(\tbcW_t)_{t \in \RR_+}$ \ are
 independent, by the SDE \eqref{Y}, the processes \ $(\cY_t)_{t \in \RR_+}$
 \ and \ $(\tbcW_t)_{t \in \RR_+}$ \ are also independent, which yields that the
 limit distribution of the CLS estimator of the offspring mean matrix
 \ $\bm_\bxi$ \ in \eqref{m_xi} is symmetric.
\proofend
\end{Rem}

\begin{Rem}\label{REMARK5}
We note that an eighth order moment condition on the offspring and immigration
 distributions in Theorem \ref{main} is supposed (i.e., we suppose
 \ $\EE(\|\bxi_{1,1,1}\|^8) < \infty$, \ $\EE(\|\bxi_{1,1,2}\|^8) < \infty$ \ and
 \ $\EE(\|\bvare_1\|^8) < \infty$).
However, it is important to remark that this condition is a technical one, we
 suspect that Theorem \ref{main} remains true under lower order moment
 condition on the offspring and immigration distributions, but we renounce to
 consider it.
\proofend
\end{Rem}

In the subcritical case we have the folowing result.
By Lemma \ref{ExUn_sub}, \ $\PP(\Omega_n) \to 1$ \ and \ $\PP(\tOmega_n) \to 1$ \ as
 \ $n \to \infty$ \ under appropriate assumptions.

\begin{Thm}\label{main_sub}
Let \ $(\bX_k)_{k \in \ZZ_+}$ \ be a 2-type Galton--Watson process with immigration
 such that \ $\alpha, \delta \in [0, 1)$ \ and \ $\beta, \gamma \in (0, \infty)$ \ with
 \ $\alpha + \delta > 0$ \ and \ $\beta \gamma < (1 - \alpha) (1 - \delta)$ \ (hence
 it is subcritical and positively regular), \ $\bX_0 = \bzero$,
 \ $\EE(\|\bxi_{1,1,1}\|^2) < \infty$, \ $\EE(\|\bxi_{1,1,2}\|^2) < \infty$,
 \ $\EE(\|\bvare_1\|^2) < \infty$, \ $\bm_\bvare \ne \bzero$, \ and at least one of the
 matrices \ $\bV_{\!\!\bxi_1}$, \ $\bV_{\!\!\bxi_2}$, \ $\bV_{\!\!\bvare}$ \ is
 invertible.
Then the probability of the existence of the estimators \ $\hbmbxi^{(n)}$ \ and
 \ $\hvarrho_n$ \ tends to 1 as \ $n \to \infty$, \ and the estimators
 \ $\hbmbxi^{(n)}$ \ and \ $\hvarrho_n$ \ are strongly consistent, i.e.,
 \ $\hbmbxi^{(n)} \as \bm_\bxi$ \ and \ $\hvarrho_n \as \varrho$ \ as \ $n \to \infty$.

If, in addition, \ $\EE(\|\bxi_{1,1,1}\|^6) < \infty$,
 \ $\EE(\|\bxi_{1,1,2}\|^6) < \infty$ \ and \ $\EE(\|\bvare_1\|^6) < \infty$, \ then
 \begin{gather}\label{m_xi_sub}
  n^{1/2} (\hbmbxi^{(n)} - \bm_\bxi) \distr \bZ , \\
  n^{1/2} (\hvarrho_n - \varrho)
  \distr \Tr(\bR \bZ)
  \distre \cN\bigl(0, \Tr\bigl[\bR^{\otimes2} \EE(\bZ^{\otimes2})\bigr]\bigr) ,
  \label{varrho_n_sub}
 \end{gather}
 as \ $n \to \infty$, \ where \ $\bZ$ \ is a \ $2 \times 2$ \ random matrix having a
 normal distribution with zero mean and with
 \begin{align*}
  \EE(\bZ^{\otimes2})
  = \Biggl\{& \sum_{i=1}^2 
               \EE\bigl[(\bxi_{1,1,i} - \EE(\bxi_{1,1,i}))^{\otimes2}\bigr]
               \EE\biggl[\tX_i \Bigl(\tbX^\top\Bigr)^{\otimes2}\biggr] \\
            &+ \EE\bigl[(\bvare_1 - \EE(\bvare_1))^{\otimes2}\bigr]
               \EE\biggl[\Bigl(\tbX^\top\Bigr)^{\otimes2}\biggr] \Biggr\}
     \biggl(\Bigl[\EE\Bigl(\tbX \tbX^\top\Bigr)\Bigr]^{\otimes2}\biggr)^{-1} ,
 \end{align*}
 where the distribution of the 2-dimensional random vector \ $\tbX$ \ is the unique
 stationary distribution of the Markov chain \ $(\bX_k)_{k\in\ZZ_+}$, \ and
 \[
   \bR := (\nabla r(\bm_\bxi))^\top
        = \frac{1}{2} \bI_2
          + \frac{1}{2\sqrt{(\alpha-\delta)^2+4\beta\gamma}}
            \begin{bmatrix}
             \alpha - \delta & 2 \beta \\
             2 \gamma & \delta - \alpha
            \end{bmatrix} .
 \]
\end{Thm}

We suspect that the moment assumptions might be relaxed to
 \ $\EE(\|\bxi_{1,1,1}\|^3) < \infty$, \ $\EE(\|\bxi_{1,1,2}\|^3) < \infty$ \ and
 \ $\EE(\|\bvare_1\|^3) < \infty$ \ by the method of Danka and Pap 
 \cite[Theorem 3.1]{DanPap}.

\section{Decomposition of the process}
\label{section_deco_proof}

Applying \eqref{mart}, let us introduce the sequence
 \begin{equation}\label{Mk}
  \bM_k
  := \bX_k - \EE(\bX_k \mid \cF_{k-1})
   = \bX_k - \bm_\bxi \bX_{k-1} - \bm_\bvare ,
  \qquad k \in \NN ,
 \end{equation}
 of martingale differences with respect to the filtration
 \ $(\cF_k)_{k \in \ZZ_+}$.
\ By \eqref{Mk}, the process \ $(\bX_k)_{k \in \ZZ_+}$ \ satisfies the recursion
 \begin{equation}\label{regr}
  \bX_k = \bm_\bxi \bX_{k-1} + \bm_\bvare + \bM_k ,
  \qquad k \in \NN .
 \end{equation}
By \eqref{CLSEabc}, for each \ $n \in \NN$, \ we have
 \[
   \hbmbxi^{(n)} - \bm_\bxi = \bD_n \bA_n^{-1} 
 \]
 on the set \ $\Omega_n$ \ given in \eqref{H_n}, where \ $\bA_n$ \ is defined
 in \eqref{A_n_B_n}, and
 \[
   \bD_n := \sum_{k=1}^n \bM_k \bX_{k-1}^\top , \qquad
   n \in \NN .
 \]
In the critical case, by \eqref{convX} and Lemmas \ref{Conv2Funct} and
 \ref{Marci}, one can derive
 \[
   n^{-3} \bA_n
   \distr
   \int_0^1 \cY_t^2 \, \dd t \, \bu_\RIGHT \bu_\RIGHT^\top
   =: \bcA
   \qquad \text{as \ $n \to \infty$.}   
 \]
Indeed, let us apply Lemma \ref{Conv2Funct} and Lemma \ref{Marci} with
 \ $K : [0,1] \times \RR^2 \to \RR^{2\times2}$,
  \[
    K(s, \bx) := \bx \bx^\top , \qquad (s, \bx) \in [0, 1] \times \RR^2 ,
  \]
 and \ $\cU := \bcX$, \ $\cU^n := \bcX^{(n)}$, \ $n \in \NN$.
\ Then
 \begin{align*}
  \| K(s, \bx) - K(t, \by) \|
  &= \| \bx \bx^\top - \by \by^\top \|
   = \| \bx \bx^\top - \by \bx^\top + \by \bx^\top - \by \by \| \\
  &\leq \| \bx - \by \| \| \bx \| + \| \by \| \| \bx - \by \|
   \leq 2 R \big( | t - s | + \| \bx - \by \| \big)
 \end{align*}
 for all \ $s, t \in [0, 1]$ \ and \ $\bx, \by \in \RR^2$ \ with
 \ $\| \bx \| \leq R$ \ and  \ $\| \by \| \leq R$, \ where \ $R > 0$.
\ Further, for all \ $n \in \NN$,
 \begin{align*}
  &\Phi_n(\bcX^{(n)})
   = \left( \bcX_1^{(n)},
            \frac{1}{n} \sum_{k=1}^n \bcX_{k/n}^{(n)} (\bcX_{k/n}^{(n)})^\top \right)
   = \left( \frac{1}{n} \bX_n,
            \frac{1}{n^3} \sum_{k=1}^n \bX_k \bX_k^\top \right) , \\
  &\Phi(\bcX)       
   = \left( \bcX_1, \int_0^1 \bcX_u \bcX_u^\top \, \dd u \right) .
 \end{align*}
Using that, by \eqref{convX}, \ $(\bcX^n, \bcX^n) \distr (\bcX, \bcX)$ \ as
 \ $n \to \infty$ \ and that the process \ $(\bcX)_{t\in\RR_+}$ \ admits
 continuous paths with probability one, Lemma \ref{Conv2Funct}
 (with the choice \ $C := \CC(\RR_+, \RR)$) \ and Lemma \ref{Marci} yield that
 \[
   n^{-3} \bA_n
   = \frac{1}{n^3} \sum_{k=1}^n \bX_k \bX_k^\top
   \distr \int_0^1 \bcX_u \bcX_u^\top \, \dd u
   \distre \int_0^1 (\cY_u \bu_\RIGHT) (\cY_u \bu_\RIGHT)^\top \dd u
   = \bcA
 \]
 as \ $n \to \infty$.
\ However, since \ $\det(\bcA) = 0$, \ the continuous mapping theorem can not
 be used for determining the weak limit of the sequence
 \ $(n^3 \bA_n^{-1})_{n \in \NN}$.
\ We can write
 \begin{equation}\label{hbmbxideco}
  \hbmbxi^{(n)} - \bm_\bxi
  = \bD_n \bA_n^{-1}
  = \frac{1}{\det(\bA_n)} \bD_n \tbA_n , \qquad n \in \NN ,
 \end{equation}
 on the set \ $\Omega_n$, \ where \ $\tbA_n$ \ denotes the adjugate of
 \ $\bA_n$ \ (i.e., the matrix of cofactors) given by
 \[
   \tbA_n
   :=\sum_{k=1}^n
      \begin{bmatrix}
       X_{k-1,2}^2 & - X_{k-1,1} X_{k-1,2} \\
       - X_{k-1,1} X_{k-1,2} & X_{k-1,1}^2
      \end{bmatrix},
   \qquad n \in \NN .
 \]
In order to prove Theorem \ref{main} we will find the asymptotic behavior of
 the sequence \ $(\det(\bA_n), \bD_n \tbA_n)_{n \in \NN}$.
\ First we derive a useful decomposition for \ $X_k$, \ $k \in \NN$.
\ Let us introduce the sequence
 \[
   U_k := \langle \bu_\LEFT, \bX_k \rangle
        = \frac{(\gamma+1-\delta) X_{k,1} + (\beta+1-\alpha) X_{k,2}}
               {1-\lambda} ,
   \qquad k \in \ZZ_+ .
 \]
One can observe that \ $U_k \geq 0$ \ for all \ $k \in \ZZ_+$, \ and
 \begin{equation}\label{rec_U}
  U_k = U_{k-1} + \langle \bu_\LEFT, \bm_\bvare \rangle
        + \langle \bu_\LEFT, \bM_k \rangle ,
  \qquad k \in \NN ,
 \end{equation}
 since
 \ $\langle \bu_\LEFT, \bm_\bxi \bX_{k-1} \rangle = \bu_\LEFT^\top \bm_\bxi \bX_{k-1}
    = \bu_\LEFT^\top \bX_{k-1} = U_{k-1}$,
 \ because \ $\bu_\LEFT$ \ is a left eigenvector of the mean matrix \ $\bm_\bxi$
 \ belonging to the eigenvalue 1.
Hence \ $(U_k)_{k \in \ZZ_+}$ \ is a nonnegative unstable AR(1) process with
 positive drift \ $\langle \bu_\LEFT, \bm_\bvare \rangle$ \ and with
 heteroscedastic innovation \ $(\langle \bu_\LEFT, \bM_k \rangle)_{k \in \NN}$.
\ Note that the solution of the recursion \eqref{rec_U} is
 \begin{equation}\label{rec_U_sol}
  U_k = \sum_{j=1}^k \langle \bu_\LEFT, \bM_j + \bm_\bvare \rangle , \qquad
  k \in \NN ,
 \end{equation}
 and, by Remark \ref{REMARK1} and Lemma \ref{Conv2Funct},
 \begin{equation}\label{convU}
  (n^{-1} U_{\nt})_{t\in\RR_+}
  = (\langle \bu_\LEFT, \bcX^{(n)}_t \rangle)_{t\in\RR_+}
  \distr (\langle \bu_\LEFT, \bcX_t \rangle)_{t\in\RR_+}
  \distre (\cY_t)_{t\in\RR_+}
  \qquad \text{as \ $n \to \infty$,}
 \end{equation}
 where \ $(\cY_t)_{t\in\RR_+}$ \ is the pathwise unique strong solution of the
 SDE \eqref{SDE_Y}.
Moreover, let
 \[
   V_k := \langle \bv_\LEFT, \bX_k \rangle
        = \frac{-(1-\alpha) X_{k,1} + \beta X_{k,2}}
               {\beta+1-\alpha} ,
   \qquad k \in \ZZ_+ .
 \]
Note that we have
 \begin{equation}\label{rec_V}
  V_k = \lambda V_{k-1} + \langle \bv_\LEFT, \bm_\bvare \rangle
        + \langle \bv_\LEFT, \bM_k \rangle ,
  \qquad k \in \NN ,
 \end{equation}
 since
 \ $\langle \bv_\LEFT, \bm_\bxi \bX_{k-1} \rangle = \bv_\LEFT^\top \bm_\bxi \bX_{k-1}
    = \lambda \bv_\LEFT^\top \bX_{k-1} = \lambda V_{k-1}$,
 \ because \ $\bv_\LEFT$ \ is a left eigenvector of the mean matrix \ $\bm_\bxi$
 \ belonging to the eigenvalue \ $\lambda$.
\ Thus \ $(V_k)_{k \in \NN}$ \ is a stable AR(1) process with drift
 \ $\langle \bv_\LEFT, \bm_\bvare \rangle$ \ and with heteroscedastic innovation
 \ $(\langle \bv_\LEFT, \bM_k \rangle)_{k \in \NN}$.
\ Note that the solution of the recursion \eqref{rec_V} is
 \begin{equation}\label{rec_V_sol}
  V_k = \sum_{j=1}^k \lambda^{k-j} \langle \bv_\LEFT, \bM_j + \bm_\bvare \rangle ,
  \qquad k \in \NN .
 \end{equation}
By \eqref{GWI(2)} and \eqref{Mk}, we obtain the decomposition
 \begin{equation}\label{Mdeco}
  \bM_k = \sum_{j=1}^{X_{k-1,1}} \big( \bxi_{k,j,1} - \EE(\bxi_{k,j,1}) \big)
          + \sum_{j=1}^{X_{k-1,2}} \big( \bxi_{k,j,2} - \EE(\bxi_{k,j,2}) \big)
          + \big( \bvare_k - \EE(\bvare_k) \big) , \qquad k \in \NN .
 \end{equation}
Under the assumption \ $\langle \obV_{\!\!\bxi} \bv_\LEFT, \bv_\LEFT\rangle = 0$,
 \ by Remark \ref{REMARK0}, we have
 \ $\langle \bv_\LEFT, \bxi_{k,j,i} - \EE(\bxi_{k,j,i})\rangle \ase 0$ \ for all
 \ $k, j \in \NN$ \ and \ $i \in \{1, 2\}$, \ and hence
 \ $\langle \bv_\LEFT, \bM_k\rangle
    \ase \langle \bv_\LEFT, \bvare_k - \EE(\bvare_k)\rangle$,
 \ implying
 \ $\langle \bv_\LEFT, \bM_k + \bm_\bvare \rangle
    \ase \langle \bv_\LEFT, \bvare_k \rangle$
 \ and, by \eqref{rec_V}, 
 \begin{equation}\label{rec_V_m}
  V_k \ase \lambda V_{k-1} + \langle \bv_\LEFT, \bvare_k \rangle ,
  \qquad k \in \NN .
 \end{equation}
The recursion \eqref{regr} has the solution
 \[
   \bX_k = \sum_{j=1}^k \bm_\bxi^{k-j} (\bm_\bvare + \bM_j) , \qquad k \in \NN .
 \]
Consequently, using \eqref{Putzer},
 \begin{align*}
  \bX_k &= \sum_{j=1}^k
            \left( \bu_\RIGHT \bu_\LEFT^\top
                   + \lambda^{k-j} \bv_\RIGHT \bv_\LEFT^\top \right)
            (\bm_\bvare + \bM_j) \\
       &= \bu_\RIGHT \bu_\LEFT^\top \sum_{j=1}^k (\bX_j - \bm_\bxi \bX_{j-1})
          + \bv_\RIGHT \bv_\LEFT^\top
            \sum_{j=1}^k
             \lambda^{k-j} (\bX_j - \bm_\bxi \bX_{j-1}) \\
       &= \bu_\RIGHT \bu_\LEFT^\top \sum_{j=1}^k (\bX_j - \bX_{j-1})
          + \bv_\RIGHT \bv_\LEFT^\top
            \sum_{j=1}^k
             \left[ \lambda^{k-j} \bX_j
                    - \lambda^{k-j+1} \bX_{j-1} \right] \\
       &= \bu_\RIGHT \bu_\LEFT^\top \bX_k + \bv_\RIGHT \bv_\LEFT^\top \bX_k
        = U_k \bu_\RIGHT + V_k \bv_\RIGHT ,
 \end{align*}
 hence
 \begin{equation}\label{XUV}
  \bX_k
  = \begin{bmatrix} X_{k,1} \\ X_{k,2} \end{bmatrix}
  = \begin{bmatrix} \bu_\RIGHT & \bv_\RIGHT \end{bmatrix}
    \begin{bmatrix} U_k \\ V_k \end{bmatrix}
  = \begin{bmatrix}
     \frac{\beta}{\beta+1-\alpha} U_k - \frac{\beta+1-\alpha}{1-\lambda} V_k \\
     \frac{1-\alpha}{\beta+1-\alpha} U_k + \frac{\gamma+1-\delta}{1-\lambda} V_k
    \end{bmatrix} ,
  \qquad k \in \ZZ_+ .
 \end{equation}
This decomposition yields
 \begin{equation}\label{A_n_deco}
  \det(\bA_n)
  = \left( \sum_{k=1}^{n-1} U_k^2 \right)
    \left( \sum_{k=1}^{n-1} V_k^2 \right)
    - \left( \sum_{k=1}^{n-1} U_k V_k \right)^2 ,
 \end{equation}
 since
 \begin{align*}
  \det(\bA_n)
  &= \det\left( \sum_{k=1}^n \bX_{k-1} \bX_{k-1}^\top \right) \\
  &= \det\left( \begin{bmatrix} \bu_\RIGHT & \bv_\RIGHT \end{bmatrix}
                \sum_{k=1}^n
                 \begin{bmatrix} U_{k-1} \\ V_{k-1} \end{bmatrix}
                 \begin{bmatrix} U_{k-1} \\ V_{k-1} \end{bmatrix}^\top
                \begin{bmatrix} \bu_\RIGHT & \bv_\RIGHT \end{bmatrix}^\top
         \right) \\
  &= \det\left( \sum_{k=1}^n
                 \begin{bmatrix} U_{k-1} \\ V_{k-1} \end{bmatrix}
                 \begin{bmatrix} U_{k-1} \\ V_{k-1} \end{bmatrix}^\top \right)
     \left[ \det\left(\begin{bmatrix} \bu_\RIGHT & \bv_\RIGHT \end{bmatrix}
                \right) \right]^2 ,
 \end{align*}
 where
 \begin{equation}\label{C}
  \det\left(\begin{bmatrix} \bu_\RIGHT & \bv_\RIGHT \end{bmatrix}\right)
  = 1 .
 \end{equation}

Theorem \ref{main} will follow from the following statements by the continuous
 mapping theorem and by Slutsky's lemma.

\begin{Thm}\label{main_Ad}
Suppose that the assumptions of Theorem \ref{main} hold. 
If \ $\langle \obV_{\!\!\bxi} \bv_\LEFT, \bv_\LEFT \rangle > 0$, \ then
 \begin{gather*}
  \sum_{k=1}^n n^{-5/2} U_{k-1} V_{k-1} \stoch 0
  \qquad \text{as \ $n \to \infty$,} \\
  \sum_{k=1}^n
   \begin{bmatrix}
    n^{-3} U_{k-1}^2 \\
    n^{-2} V_{k-1}^2 \\
    n^{-2} \bM_k U_{k-1} \\
    n^{-3/2} \bM_k V_{k-1}
   \end{bmatrix}
  \distr
  \begin{bmatrix}
   \int_0^1 \cY_t^2 \, \dd t \\[1mm]
   \frac{\langle \obV_{\!\!\bxi} \bv_\LEFT, \bv_\LEFT \rangle}
        {1-\lambda^2}
   \int_0^1 \cY_t \, \dd t \\[1mm]
   \int_0^1 \cY_t \, \dd \bcM_t \\[1mm]
   \frac{\langle \obV_{\!\!\bxi} \bv_\LEFT, \bv_\LEFT \rangle^{1/2}}
        {(1-\lambda^2)^{1/2}}
   \obV_{\!\!\bxi}^{1/2} \!\!
   \int_0^1 \cY_t \, \dd\tbcW_t
  \end{bmatrix}
  \qquad \text{as \ $n \to \infty$.}
 \end{gather*}
\end{Thm}

In case of \ $\langle \obV_{\!\!\bxi} \bv_\LEFT, \bv_\LEFT \rangle = 0$ \ the
 second and fourth coordinates of the limit vector of the second convergence
 is 0 in Theorem \ref{main_Ad}, thus other scaling factors should be chosen
 for these coordinates, described in the following theorem.

\begin{Thm}\label{maint_Ad}
Suppose that the assumptions of Theorem \ref{main} hold. 
If \ $\langle \obV_{\!\!\bxi} \bv_\LEFT, \bv_\LEFT \rangle = 0$, \ then
 \begin{gather*}
  \sum_{k=1}^n n^{-1} V_{k-1}^2
  \stoch \frac{\langle \bv_\LEFT, \bm_\bvare \rangle^2}{(1-\lambda)^2}
         + \frac{\langle \bV_{\!\!\bvare} \bv_\LEFT, \bv_\LEFT \rangle}
                {1-\lambda^2}
         =: M
  \qquad \text{as \ $n \to \infty$,} \\
   \sum_{k=1}^n
    \begin{bmatrix}
     n^{-3} U_{k-1}^2 \\
     n^{-2} U_{k-1} V_{k-1} \\
     n^{-2} \bM_k U_{k-1} \\
     n^{-1} \bM_k V_{k-1}
    \end{bmatrix}
   \distr
   \begin{bmatrix}
    \int_0^1 \cY_t^2 \, \dd t \\[1mm]
    \frac{\langle\bv_\LEFT,\bm_\bvare\rangle}{1-\lambda} \,
    \int_0^1 \cY_t \, \dd t \\[1mm]
    \int_0^1 \cY_t \, \dd \bcM_t \\[1mm]
    \frac{\langle \bv_\LEFT, \bm_\bvare \rangle}{1-\lambda} \, \bcM_1
    + \frac{\langle \bV_{\!\!\bvare} \bv_\LEFT, \bv_\LEFT \rangle^{1/2}}
           {(1-\lambda^2)^{1/2}}
      \obV_{\!\!\bxi}^{1/2} \!\!
      \int_0^1 \cY_t \, \dd\tbcW_t
   \end{bmatrix}
   \qquad \text{as \ $n \to \infty$.}
 \end{gather*}
\end{Thm}

\noindent
\textbf{Proof of Theorem \ref{main}.}
In order to derive the statements, we can use the continuous mapping theorem
 and Slutsky's lemma.

Theorem \ref{main_Ad} implies \eqref{m_xi}.
Indeed, we can use the representation \eqref{hbmbxideco}, where the adjugate
 \ $\tbA_n$ \ can be written in the form
 \[
   \tbA_n = \begin{bmatrix} 0 & 1 \\ -1 & 0 \end{bmatrix}
            \sum_{\ell=1}^n \bX_{\ell-1} \bX_{\ell-1}^\top
            \begin{bmatrix} 0 & -1 \\ 1 & 0 \end{bmatrix} ,
   \qquad n \in \NN .
 \]
Using \eqref{XUV}, we have
 \begin{align*}
  \bD_n \tbA_n
  &= \sum_{k=1}^n
      \bM_k \begin{bmatrix} U_{k-1} \\ V_{k-1} \end{bmatrix}^\top
     \begin{bmatrix} \bu_\RIGHT^\top \\ \bv_\RIGHT^\top \end{bmatrix}
     \begin{bmatrix} 0 & 1 \\ -1 & 0 \end{bmatrix}
     \begin{bmatrix} \bu_\RIGHT^\top \\ \bv_\RIGHT^\top \end{bmatrix}^\top
     \sum_{\ell=1}^n
      \begin{bmatrix} U_{\ell-1} \\ V_{\ell-1} \end{bmatrix}
      \begin{bmatrix} U_{\ell-1} \\ V_{\ell-1} \end{bmatrix}^\top
     \begin{bmatrix} \bu_\RIGHT^\top \\ \bv_\RIGHT^\top \end{bmatrix}
     \begin{bmatrix} 0 & -1 \\ 1 & 0 \end{bmatrix} .
 \end{align*}
Here we have
 \[
   \begin{bmatrix} \bu_\RIGHT^\top \\ \bv_\RIGHT^\top \end{bmatrix}
   \begin{bmatrix} 0 & 1 \\ -1 & 0 \end{bmatrix}
   \begin{bmatrix} \bu_\RIGHT^\top \\ \bv_\RIGHT^\top \end{bmatrix}^\top
   = \begin{bmatrix} 0 & 1 \\ -1 & 0 \end{bmatrix} , \qquad
   \begin{bmatrix} \bu_\RIGHT^\top \\ \bv_\RIGHT^\top \end{bmatrix}
   \begin{bmatrix} 0 & -1 \\ 1 & 0 \end{bmatrix}
   = \begin{bmatrix} - \bv_\LEFT^\top \\ \bu_\LEFT^\top \end{bmatrix} .
 \]
Theorem \ref{main_Ad} implies asymptotic expansions
 \begin{align*}
  \sum_{k=1}^n \bM_k \begin{bmatrix} U_{k-1} \\ V_{k-1} \end{bmatrix}^\top
  &= n^2 \bD_{n,1} + n^{3/2} \bD_{n,2} , \\
     \sum_{\ell=1}^n
      \begin{bmatrix} U_{\ell-1} \\ V_{\ell-1} \end{bmatrix}
      \begin{bmatrix} U_{\ell-1} \\ V_{\ell-1} \end{bmatrix}^\top
  &= n^3 \bA_{n,1} + n^{5/2} \bA_{n,2} + n^2 \bA_{n,3} , 
 \end{align*}
 where
 \begin{align*}
  \bD_{n,1}
  &:= n^{-2} \sum_{k=1}^n \bM_k U_{k-1} \be_1^\top
   \distr
   \int_0^1 \cY_t \, \dd \bcM_t \, \be_1^\top
   =: \bcD_1 , \\
  \bD_{n,2}
  &:= n^{-3/2}
      \sum_{k=1}^n \bM_k V_{k-1} \be_2^\top
   \distr
   \frac{\langle \obV_{\!\!\bxi} \bv_\LEFT, \bv_\LEFT \rangle^{1/2}}
        {(1-\lambda^2)^{1/2}}
   \obV_{\!\!\bxi}^{1/2} \!\!
   \int_0^1 \cY_t \, \dd\tbcW_t \, \be_2^\top
   =: \bcD_2 , \\
  \bA_{n,1}
  &:= n^{-3} \sum_{\ell=1}^n \begin{bmatrix} U_{\ell-1}^2 & 0 \\ 0 & 0 \end{bmatrix}
   \distr
   \int_0^1 \cY_t^2 \, \dd t
   \begin{bmatrix} 1 & 0 \\ 0 & 0 \end{bmatrix}
   =: \bcA_1 , \\
  \bA_{n,2}
  &:= n^{-5/2}
      \sum_{\ell=1}^n
       \begin{bmatrix}
        0 & U_{\ell-1} V_{\ell-1} \\
        U_{\ell-1} V_{\ell-1} & 0
       \end{bmatrix}
   \distr
   \bzero , \\
  \bA_{n,3}
  &:= n^{-2}
      \sum_{\ell=1}^n \begin{bmatrix} 0 & 0 \\ 0 & V_{\ell-1}^2 \end{bmatrix}
   \distr
   \frac{\langle \obV_{\!\!\bxi} \bv_\LEFT, \bv_\LEFT \rangle}{1-\lambda^2}
   \int_0^1 \cY_t \, \dd t
   \begin{bmatrix} 0 & 0 \\ 0 & 1 \end{bmatrix}
   =: \bcA_3
 \end{align*}
 jointly as \ $n \to \infty$.
\ Consequently, we obtain an asymptotic expansion
 \begin{align*}
  \bD_n \tbA_n
  &= (n^2 \bD_{n,1} + n^{3/2} \bD_{n,2})
     \begin{bmatrix} 0 & 1 \\ -1 & 0 \end{bmatrix}
     (n^3 \bA_{n,1} + n^{5/2} \bA_{n,2} + n^2 \bA_{n,3})
     \begin{bmatrix} -\bv_\LEFT^\top \\ \bu_\LEFT^\top \end{bmatrix} \\
  &= (n^5 \bC_{n,1} + n^{9/2} \bC_{n,2} + n^4 \bC_{n,3} + n^{7/2} \bC_{n,4})
     \begin{bmatrix} -\bv_\LEFT^\top \\ \bu_\LEFT^\top \end{bmatrix} ,
 \end{align*}
 where
 \[
   \bC_{n,1}
   := \bD_{n,1} \begin{bmatrix} 0 & 1 \\ -1 & 0 \end{bmatrix} \bA_{n,1}
    = n^{-5}
      \sum_{k=1}^n
       \sum_{\ell=1}^n
        \bM_k U_{k-1} U_{\ell-1}^2 \be_1^\top
        \begin{bmatrix} 0 & 1 \\ -1 & 0 \end{bmatrix}
        \begin{bmatrix} 1 & 0 \\ 0 & 0 \end{bmatrix}
    = \bzero
 \]
 for all \ $n \in \NN$, \ and
 \begin{align*}
  \bC_{n,2}
  &:= \bD_{n,1} \begin{bmatrix} 0 & 1 \\ -1 & 0 \end{bmatrix} \bA_{n,2}
      + \bD_{n,2} \begin{bmatrix} 0 & 1 \\ -1 & 0 \end{bmatrix} \bA_{n,1}
   \distr \bcD_2 \begin{bmatrix} 0 & 1 \\ -1 & 0 \end{bmatrix} \bcA_1 , \\
  \bC_{n,3}
  &:= \bD_{n,1} \begin{bmatrix} 0 & 1 \\ -1 & 0 \end{bmatrix} \bA_{n,3}
      + \bD_{n,2} \begin{bmatrix} 0 & 1 \\ -1 & 0 \end{bmatrix} \bA_{n,2}
   \distr \bcD_1 \begin{bmatrix} 0 & 1 \\ -1 & 0 \end{bmatrix} \bcA_3 , \\
  \bC_{n,4}
  &:= \bD_{n,2} \begin{bmatrix} 0 & 1 \\ -1 & 0 \end{bmatrix} \bA_{n,3}
   \distr \bcD_2 \begin{bmatrix} 0 & 1 \\ -1 & 0 \end{bmatrix} \bcA_3
 \end{align*}
 as \ $n \to \infty$.
\ Using again Theorem \ref{main_Ad} and \eqref{A_n_deco}, we conclude
 \[
   \begin{bmatrix} n^{-5} \det(\bA_n) \\ n^{-9/2} \bD_n \tbA_n \end{bmatrix}
   \distr
   \begin{bmatrix} 
    \frac{\langle\obV_{\!\!\bxi}\bv_\LEFT,\bv_\LEFT\rangle}{1-\lambda^2}
    \int_0^1 \cY_t^2 \, \dd t
    \int_0^1 \cY_t \, \dd t \\[2mm] 
    \bcD_2 \begin{bmatrix} 0 & 1 \\ -1 & 0 \end{bmatrix} \bcA_1
    \begin{bmatrix} -\bv_\LEFT^\top \\ \bu_\LEFT^\top \end{bmatrix}
   \end{bmatrix}
   \qquad \text{as \ $n \to \infty$.}
 \]
Here
 \begin{align*}
  &\bcD_2 \begin{bmatrix} 0 & 1 \\ -1 & 0 \end{bmatrix} \bcA_1
   \begin{bmatrix} -\bv_\LEFT^\top \\ \bu_\LEFT^\top \end{bmatrix} \\
  &= \frac{\langle \obV_{\!\!\bxi} \bv_\LEFT, \bv_\LEFT \rangle^{1/2}}
          {(1-\lambda^2)^{1/2}}
     \int_0^1 \cY_t^2 \, \dd t \,
     \obV_{\!\!\bxi}^{1/2} \!\!
     \int_0^1 \cY_t \, \dd\tbcW_t \, \be_2^\top
     \begin{bmatrix} 0 & 1 \\ -1 & 0 \end{bmatrix}
     \begin{bmatrix} 1 & 0 \\ 0 & 0 \end{bmatrix}
     \begin{bmatrix} -\bv_\LEFT^\top \\ \bu_\LEFT^\top \end{bmatrix} \\
  &= \frac{\langle \obV_{\!\!\bxi} \bv_\LEFT, \bv_\LEFT \rangle^{1/2}}
          {(1-\lambda^2)^{1/2}}
     \int_0^1 \cY_t^2 \, \dd t \,
     \obV_{\!\!\bxi}^{1/2} \!\!
     \int_0^1 \cY_t \, \dd\tbcW_t \, \bv_\LEFT^\top .
 \end{align*}
Since \ $\bm_\bvare \ne \bzero$, \ by the SDE \eqref{SDE_Y}, we have
 \ $\PP(\text{$\cY_t = 0$ \ for all \ $t \in [0,1]$}) = 0$, \ which implies
 that
 \ $\PP\bigl( \int_0^1 \cY_t^2 \, \dd t \int_0^1 \cY_t \, \dd t > 0 \bigr) = 1$,
 \ hence the continuous mapping theorem implies \eqref{m_xi}.

For the proof of \eqref{varrho_n}, we can write \ $\hvarrho_n - 1$ \ in the
 form
 \begin{align*}
  \hvarrho_n - 1
  &= \hvarrho_n - \varrho
   = \frac{(\halpha_n - \alpha) + (\hdelta_n - \delta)
           + \sqrt{(\halpha_n - \hdelta_n)^2 + 4 \hbeta_n \hgamma_n}
           - \sqrt{(\alpha - \delta)^2 + 4 \beta \gamma}}
          {2} \\
  &= \frac{(\halpha_n - \alpha) + (\hdelta_n - \delta)}{2} \\
  &\quad
     + \frac{\bigl[(\halpha_n - \alpha) - (\hdelta_n - \delta)\bigr]
             \bigl[(\halpha_n + \alpha) - (\hdelta_n + \delta)\bigr]
             + 4 (\hbeta_n - \beta) \hgamma_n + 4 (\hgamma_n - \gamma) \beta}
            {2 \Bigl(\sqrt{(\halpha_n - \hdelta_n)^2 + 4 \hbeta_n \hgamma_n}
                     + \sqrt{(\alpha - \delta)^2 + 4 \beta \gamma}\Bigr)} \\
  &= \frac{c_n}
          {2 \Bigl(\sqrt{(\halpha_n - \hdelta_n)^2 + 4 \hbeta_n \hgamma_n}
                   + (1 - \lambda)\Bigr)} ,
 \end{align*}
 where
 \begin{align*}
  c_n &:= (\halpha_n - \alpha)
          \Bigl[\sqrt{(\halpha_n - \hdelta_n)^2 + 4 \hbeta_n \hgamma_n}
                + (1 - \lambda) + (\halpha_n + \alpha) - (\hdelta_n + \delta)\Bigr]
          + 4 (\hbeta_n - \beta) \hgamma_n \\
      &\quad
          + 4 (\hgamma_n - \gamma) \beta
          + (\hdelta_n - \delta)
            \Bigl[\sqrt{(\halpha_n - \hdelta_n)^2 + 4 \hbeta_n \hgamma_n}
                  + (1 - \lambda) - (\halpha_n + \alpha)
                  + (\hdelta_n + \delta)\Bigr] .
 \end{align*} 
Slutsky's lemma and \eqref{m_xi} imply
 \ $\hbmbxi^{(n)} - \bm_\bxi \distr \bzero$, \ and hence
 \ $\hbmbxi^{(n)} - \bm_\bxi \stoch \bzero$ \ as \ $n \to \infty$.
\ Thus \ $\hgamma_n \stoch \gamma$, \ and
 \begin{align*}
  (\halpha_n + \alpha) - (\hdelta_n + \delta)
  &\stoch
   2 (\alpha - \delta) , \\
  (\halpha_n - \hdelta_n)^2 + 4 \hbeta_n \hgamma_n
  &\stoch
   (\alpha - \delta)^2 + 4 \beta \gamma 
   = (2 - \alpha - \delta)^2
   = (1 - \lambda)^2
 \end{align*}
 as \ $n \to \infty$.
\ The aim of the following discussion is to show \ $n (c_n - d_n) \distr 0$ \ as
 \ $n \to \infty$, \ where
 \begin{align*}
  d_n &:= (\halpha_n - \alpha) \bigl[ 2(1 - \lambda) + 2 \alpha - 2 \delta\bigr]
          + 4 (\hbeta_n - \beta) \gamma \\
      &\:\quad
          + 4 (\hgamma_n - \gamma) \beta
          + (\hdelta_n - \delta)
            \bigl[ 2(1 - \lambda) - 2 \alpha + 2 \delta\bigr] \\
      &\:= 4 (1 - \delta) (\halpha_n - \alpha) + 4 \gamma (\hbeta_n - \beta)
           + 4 \beta (\hgamma_n - \gamma) + 4 (1 - \alpha) (\hdelta_n - \delta) .
 \end{align*}
We have
 \begin{align*}
  n (c_n - d_n)
  &= n \Bigl[ \sqrt{(\halpha_n - \hdelta_n)^2 + 4 \hbeta_n \hgamma_n}
              - (1 - \lambda) \Bigr]
       \bigl[ (\halpha_n - \alpha) + (\hdelta_n - \delta) \bigr] \\
  &\quad
     + n \bigl[ (\halpha_n - \alpha) - (\hdelta_n - \delta) \bigr]^2
     + 4 n (\hbeta_n - \beta) (\hgamma_n - \gamma) ,
 \end{align*}
 where
 \begin{multline*}
  \sqrt{(\halpha_n - \hdelta_n)^2 + 4 \hbeta_n \hgamma_n} - (1 - \lambda) \\
  = \frac{\bigl[ (\halpha_n - \alpha) - (\hdelta_n - \delta) \bigr]
          \bigl[ (\halpha_n + \alpha) - (\hdelta_n + \delta) \bigr]
          + 4 (\hbeta_n - \beta) \hgamma_n + 4 (\hgamma_n - \gamma) \beta}
         {\sqrt{(\halpha_n - \hdelta_n)^2 + 4 \hbeta_n \hgamma_n}
          + (1 - \lambda)} ,
 \end{multline*}
 hence \ $n (c_n - d_n) \distr 0$ \ as \ $n \to \infty$ \ will follow from
 \begin{equation}\label{c-d}
  \begin{aligned}
   &n \bigl\{\bigl[(\halpha_n - \alpha) - (\hdelta_n - \delta)\bigr]
             \bigl[(\halpha_n + \alpha) - (\hdelta_n + \delta)\bigr]
             + 4 (\hbeta_n - \beta) \hgamma_n
             + 4 (\hgamma_n - \gamma)\beta\bigr\} \\
   &\times
      \bigl[(\halpha_n - \alpha) + (\hdelta_n - \delta)\bigr] \\
   & + n \Bigl[\sqrt{(\halpha_n - \hdelta_n)^2 + 4 \hbeta_n \hgamma_n}
               + (1 - \lambda)\Bigr] \\
   &\phantom{+}\times
         \Bigl\{\bigl[(\halpha_n - \alpha) - (\hdelta_n - \delta)\bigr]^2
                + 4 (\hbeta_n - \beta) (\hgamma_n - \gamma)\Bigr\}
     \distr 0 \qquad \text{as \ $n \to \infty$.}
  \end{aligned}
 \end{equation}
By \eqref{m_xi} we have
 \[
   n^{1/2} (\hbmbxi^{(n)} - \bm_\bxi)
   \distr
   \cC \obV_{\!\!\bxi}^{1/2} \int_0^1 \cY_t \, \dd\tbcW_t \, \bv_\LEFT^\top \qquad
   \text{as \ $n \to \infty$}
 \]
 with
 \[
   \cC := \frac{(1-\lambda)^{1/2}}
               {\langle\obV_{\!\!\bxi}\bv_\LEFT,\bv_\LEFT\rangle^{1/2}
                \int_0^1 \cY_t \, \dd t} .
 \]
Consequently,
 \[
   \begin{bmatrix}
    n^{1/2} (\halpha_n - \alpha) \\
    n^{1/2} (\hbeta_n - \beta) \\
    n^{1/2} (\hgamma_n - \gamma) \\
    n^{1/2} (\hdelta_n - \delta)
   \end{bmatrix}
   \distr
   \cC \begin{bmatrix}
        \be_1^\top \obV_{\!\!\bxi}^{1/2} \int_0^1 \cY_t \, \dd\tbcW_t \, \bv_\LEFT^\top
        \be_1 \\
        \be_1^\top \obV_{\!\!\bxi}^{1/2} \int_0^1 \cY_t \, \dd\tbcW_t \, \bv_\LEFT^\top
        \be_2 \\
        \be_2^\top \obV_{\!\!\bxi}^{1/2} \int_0^1 \cY_t \, \dd\tbcW_t \, \bv_\LEFT^\top
        \be_1 \\
        \be_2^\top \obV_{\!\!\bxi}^{1/2} \int_0^1 \cY_t \, \dd\tbcW_t \, \bv_\LEFT^\top
        \be_2
      \end{bmatrix}
   = \frac{\cC}{\beta+1-\alpha}
     \begin{bmatrix}
      - (1-\alpha) \be_1^\top \bcI \\
      \beta \be_1^\top \bcI \\
      - (1-\alpha) \be_2^\top \bcI \\
      \beta \be_2^\top \bcI
     \end{bmatrix}
 \]
 as \ $n \to \infty$ \ with
 \ $\bcI := \obV_\bxi^{1/2} \int_0^1 \cY_t \, \dd\tbcW_t$.
\ By continuous mapping theorem,
 \[
   \begin{bmatrix}
    n (\halpha_n - \alpha)^2 \\
    n (\halpha_n - \alpha) (\hbeta_n - \beta) \\
    n (\halpha_n - \alpha) (\hgamma_n - \gamma) \\
    n (\halpha_n - \alpha) (\hdelta_n - \delta) \\
    n (\hbeta_n - \beta) (\hgamma_n - \gamma) \\
    n (\hbeta_n - \beta) (\hdelta_n - \delta) \\
    n (\hgamma_n - \gamma) (\hdelta_n - \delta) \\
    n (\hdelta_n - \delta)^2
   \end{bmatrix}
   \distr
   \frac{\cC^2}{(\beta+1-\alpha)^2}
   \begin{bmatrix}
    (1-\alpha)^2 \bcI^\top \be_1 \be_1^\top \bcI \\
    - (1-\alpha) \beta \bcI^\top \be_1 \be_1^\top \bcI \\
    (1-\alpha)^2 \bcI^\top \be_1 \be_2^\top \bcI \\
    - (1-\alpha) \beta \bcI^\top \be_1 \be_2^\top \bcI \\
    - (1-\alpha) \beta \bcI^\top \be_1 \be_2^\top \bcI \\
    \beta^2 \bcI^\top \be_1 \be_2^\top \bcI \\
    - (1-\alpha) \beta \bcI^\top \be_2 \be_2^\top \bcI \\
    \beta^2 \bcI^\top \be_2 \be_2^\top \bcI
   \end{bmatrix}
 \]
 as \ $n \to \infty$, \ and by continuous mapping theorem, Slutsky’s lemma and
 \ $\hbmbxi^{(n)} \distr \bm_\bxi$ \ as \ $n \to \infty$, \ we conclude
 \eqref{c-d}, since
 \begin{align*}
  &2 (\alpha - \delta)
   \bigl[(1-\alpha)^2 \be_1 \be_1^\top - \beta^2 \be_2 \be_2^\top\bigr]
   + 4 \gamma 
       \bigl[- (1-\alpha) \beta \be_1 \be_1^\top + \beta^2 \be_1 \be_2^\top\bigr] \\
  &+ 4 \beta
       \bigl[(1-\alpha)^2 \be_1 \be_2^\top
             - (1-\alpha) \beta \be_2 \be_2^\top\bigr] \\
  &+ 2 (1 - \lambda)
     \bigl[(1-\alpha)^2 \be_1 \be_1^\top + 2 (1-\alpha) \beta \be_1 \be_2^\top
           + \beta^2 \be_2 \be_2^\top
           - 4 (1 -\alpha) \beta \be_1 \be_2^\top\bigr] \\
  &= \be_1 \be_1^\top
     \bigl[2 (\alpha - \delta) (1-\alpha)^2 - 4 \beta \gamma (1-\alpha)
           + 2 (1 - \lambda) (1-\alpha)^2\bigr] \\
  &\quad
     + \be_1 \be_2^\top
       \bigl[4 \beta^2 \gamma + 4 \beta (1-\alpha)^2
             - 4 \beta (1-\alpha) (1 - \lambda)\bigr] \\
  &\quad
     + \be_2 \be_2^\top
       \bigl[- 2 \beta^2 (\alpha - \delta) - 4 \beta^2 (1-\alpha)
             + 2 \beta^2 (1 - \lambda)\bigr]
   = 0 .
 \end{align*}
Consequently, \eqref{varrho_n} will follow from
 \begin{equation}\label{tvarrho_n}
  \frac{nd_n}{2 \Bigl[\sqrt{(\halpha_n - \hdelta_n)^2 + 4 \hbeta_n \hgamma_n}
                      + (1 - \lambda)\Bigr]}
  \distr
  \frac{\int_0^1 \cY_t \, \dd(\cY_t - t \langle \bu_\LEFT, \bm_\bvare \rangle)}
       {\int_0^1 \cY_t^2 \, \dd t}
 \end{equation}
 as \ $n \to \infty$.
\ We can write
 \begin{align*}
  d_n
  &= 4 (1-\delta) \be_1^\top (\hbmbxi^{(n)}-\bm_\bxi) \be_1
     + 4 \gamma \be_1^\top (\hbmbxi^{(n)}-\bm_\bxi) \be_2 \\
  &\quad
     + 4 \beta \be_2^\top (\hbmbxi^{(n)}-\bm_\bxi) \be_1
     + 4 (1-\alpha) \be_2^\top (\hbmbxi^{(n)}-\bm_\bxi) \be_2 .
 \end{align*}
We use again the representation \eqref{hbmbxideco} and the asymptotic expansion of
 \ $\bD_n \tbA_n$. 
\ We have
 \begin{gather*}
  \bD_{n,1} \begin{bmatrix} 0 & 1 \\ -1 & 0 \end{bmatrix} \bA_{n,2}
  \begin{bmatrix} -\bv_\LEFT^\top \\ \bu_\LEFT^\top \end{bmatrix}
  = - n^{-9/2}
       \sum_{k=1}^n 
        \sum_{\ell=1}^n
         \bM_k U_{k-1} U_{\ell-1} V_{\ell-1}
       \bv_\LEFT^\top , \\
  \bD_{n,2} \begin{bmatrix} 0 & 1 \\ -1 & 0 \end{bmatrix} \bA_{n,1}
  \begin{bmatrix} -\bv_\LEFT^\top \\ \bu_\LEFT^\top \end{bmatrix}
  = n^{-9/2}
     \sum_{k=1}^n 
      \sum_{\ell=1}^n
       \bM_k V_{k-1} U_{\ell-1}^2
       \bv_\LEFT^\top ,
 \end{gather*}
 implying
 \begin{align*}
  &\Biggl((1-\delta)
          \be_1^\top \bD_{n,1} \begin{bmatrix} 0 & 1 \\ -1 & 0 \end{bmatrix}
          \bA_{n,2} \be_1
          + \gamma \be_1^\top \bD_{n,1} \begin{bmatrix} 0 & 1 \\ -1 & 0 \end{bmatrix}
            \bA_{n,2} \be_2 \\
  &\phantom{\biggl(}
          + \beta \be_2^\top \bD_{n,1} \begin{bmatrix} 0 & 1 \\ -1 & 0 \end{bmatrix}
            \bA_{n,2} \be_1
          + (1-\alpha)
            \be_2^\top \bD_{n,1} \begin{bmatrix} 0 & 1 \\ -1 & 0 \end{bmatrix}
            \bA_{n,2} \be_2\Biggr)
   \begin{bmatrix} -\bv_\LEFT^\top \\ \bu_\LEFT^\top \end{bmatrix}
   = \bzero , \\
  &\Biggl((1-\delta) \be_1^\top \bD_{n,2} \begin{bmatrix} 0 & 1 \\ -1 & 0 \end{bmatrix}
          \bA_{n,1} \be_1
          + \gamma \be_1^\top \bD_{n,2} \begin{bmatrix} 0 & 1 \\ -1 & 0 \end{bmatrix}
            \bA_{n,1} \be_2 \\
  &\phantom{\biggl(}
          + \beta \be_2^\top \bD_{n,2} \begin{bmatrix} 0 & 1 \\ -1 & 0 \end{bmatrix}
            \bA_{n,1} \be_1
          + (1-\alpha)
            \be_2^\top \bD_{n,2} \begin{bmatrix} 0 & 1 \\ -1 & 0 \end{bmatrix}
            \bA_{n,1} \be_2\Biggr)
   \begin{bmatrix} -\bv_\LEFT^\top \\ \bu_\LEFT^\top \end{bmatrix}
   = \bzero
 \end{align*}
 for all \ $n \in \NN$, \ and
 \begin{align*}
  &(1-\delta) \be_1^\top \bD_{n,2} \begin{bmatrix} 0 & 1 \\ -1 & 0 \end{bmatrix}
   \bA_{n,2} \be_1
   + \gamma \be_1^\top \bD_{n,2} \begin{bmatrix} 0 & 1 \\ -1 & 0 \end{bmatrix}
     \bA_{n,2} \be_2 \\
  &+ \beta \be_2^\top \bD_{n,2} \begin{bmatrix} 0 & 1 \\ -1 & 0 \end{bmatrix}
     \bA_{n,2} \be_1
   + (1-\alpha) \be_2^\top \bD_{n,2} \begin{bmatrix} 0 & 1 \\ -1 & 0 \end{bmatrix}
     \bA_{n,2} \be_2
   \stoch \bzero
 \end{align*}
 as \ $n \to \infty$.
\ Moreover,
 \[
   \bcD_1 \begin{bmatrix} 0 & 1 \\ -1 & 0 \end{bmatrix} \bcA_3
   \begin{bmatrix} -\bv_\LEFT^\top \\ \bu_\LEFT^\top \end{bmatrix} \\
   = \frac{\langle \obV_{\!\!\bxi} \bv_\LEFT, \bv_\LEFT \rangle}{1-\lambda^2}
     \int_0^1 \cY_t \, \dd t
     \int_0^1 \cY_t \, \dd\bcM_t \,
     \bu_\LEFT^\top .
 \]
Using again Theorem \ref{main_Ad} and \eqref{A_n_deco}, we conclude
 \begin{align*}
  &\frac{nd_n}{2 \Bigl[\sqrt{(\halpha_n - \hdelta_n)^2 + 4 \hbeta_n \hgamma_n}
               + (1 - \lambda)\Bigr]} \\
  &\distr
   \frac{1}{(1-\lambda)\int_0^1 \cY_t^2 \, \dd t}
   \biggl((1-\delta) \be_1^\top \int_0^1 \cY_t \, \dd\bcM_t \, \bu_\LEFT^\top \be_1
          + \gamma \be_1^\top \int_0^1 \cY_t \, \dd\bcM_t \, \bu_\LEFT^\top
            \be_2 \\
  &\phantom{\distr\frac{1}{(1-\lambda)\int_0^1 \cY_t^2 \, \dd t}\biggl(}
   + \beta \be_2^\top \int_0^1 \cY_t \, \dd\bcM_t \, \bu_\LEFT^\top \be_1
   + (1-\alpha) \be_2^\top \int_0^1 \cY_t \, \dd\bcM_t \, \bu_\LEFT^\top
     \be_2\biggr) \\
  &= \frac{1}{(1-\lambda)^2\int_0^1 \cY_t^2 \, \dd t}
     \biggl(\bigl[(1-\delta) (\gamma+1-\delta) + \gamma (\beta+1-\alpha)\bigr]
            \be_1^\top \int_0^1 \cY_t \, \dd\bcM_t \\
  &\phantom{\distr\frac{1}{(1-\lambda)\int_0^1 \cY_t^2 \, \dd t}\biggl(}
            + \bigl[\beta (\gamma+1-\delta) + (1-\delta) (\beta+1-\alpha)\bigr] 
              \be_2^\top \int_0^1 \cY_t \, \dd\bcM_t \biggr) \\
  &= \frac{1}{(1-\lambda)\int_0^1 \cY_t^2 \, \dd t}
     \biggl((\gamma+1-\delta) \be_1^\top \int_0^1 \cY_t \, \dd\bcM_t
             + (\beta+1-\alpha) \be_2^\top \int_0^1 \cY_t \, \dd\bcM_t\biggr) \\
  &= \frac{\int_0^1 \cY_t \, \dd\langle\bu_\LEFT,\bcM_t\rangle}
          {\int_0^1 \cY_t^2 \, \dd t}
   = \frac{\int_0^1 \cY_t \, \dd(\cY_t - t \langle \bu_\LEFT, \bm_\bvare \rangle)}
          {\int_0^1 \cY_t^2 \, \dd t}
 \end{align*}
 as \ $n \to \infty$.
 
Next we show that Theorem \ref{maint_Ad} yields \eqref{m_xi_m}.
Theorem \ref{maint_Ad} implies asymptotic expansions
 \begin{align*}
  \sum_{k=1}^n \bM_k \begin{bmatrix} U_{k-1} \\ V_{k-1} \end{bmatrix}^\top
  &= n^2 \bD_{n,1} + n \bD_{n,2} , \\
     \sum_{\ell=1}^n
      \begin{bmatrix} U_{\ell-1} \\ V_{\ell-1} \end{bmatrix}
      \begin{bmatrix} U_{\ell-1} \\ V_{\ell-1} \end{bmatrix}^\top
  &= n^3 \bA_{n,1} + n^2 \bA_{n,2} + n \bA_{n,3} , 
 \end{align*}
 where
 \begin{align*}
  \bD_{n,1}
  &:= n^{-2} \sum_{k=1}^n \bM_k U_{k-1} \be_1^\top
   \distr
   \int_0^1 \cY_t \, \dd \bcM_t \, \be_1^\top
   =: \bcD_1 , \\
  \bD_{n,2}
  &:= n^{-1} \sum_{k=1}^n \bM_k V_{k-1} \be_2^\top \\
  &\distr
   \biggl(\frac{\langle\bv_\LEFT,\bm_\bvare\rangle}{1-\lambda} \bcM_1
          + \frac{\langle \bV_{\!\!\bvare} \bv_\LEFT, \bv_\LEFT \rangle^{1/2}}
                 {(1-\lambda^2)^{1/2}}
          \obV_{\!\!\bxi}^{1/2} \!\! \int_0^1 \cY_t \, \dd\tbcW_t \biggr)
   \be_2^\top
   =: \bcD_2 , \\
  \bA_{n,1}
  &:= n^{-3}
      \sum_{\ell=1}^n \begin{bmatrix} U_{\ell-1}^2 & 0 \\ 0 & 0 \end{bmatrix}
   \distr
   \int_0^1 \cY_t^2 \, \dd t
   \begin{bmatrix} 1 & 0 \\ 0 & 0 \end{bmatrix}
   =: \bcA_1 , \\
  \bA_{n,2}
  &:= n^{-2}
      \sum_{\ell=1}^n
       \begin{bmatrix}
        0 & U_{\ell-1} V_{\ell-1} \\
        U_{\ell-1} V_{\ell-1} & 0
       \end{bmatrix}
   \distr
   \frac{\langle\bv_\LEFT,\bm_\bvare\rangle}{1-\lambda}
   \int_0^1 \cY_t^2 \, \dd t
   \begin{bmatrix} 0 & 1 \\ 1 & 0 \end{bmatrix}
   =: \bcA_2 , \\
  \bA_{n,3}
  &:= n^{-1}
      \sum_{\ell=1}^n \begin{bmatrix} 0 & 0 \\ 0 & V_{\ell-1}^2 \end{bmatrix}
   \distr
   \biggl(\frac{\langle\bv_\LEFT,\bm_\bvare\rangle^2}{(1-\lambda)^2}
          + \frac{\langle\bV_{\!\!\bvare}\bv_\LEFT,\bv_\LEFT\rangle}
                 {1-\lambda^2}\biggr)
   \begin{bmatrix}
    0 & 0 \\
    0 & 1
   \end{bmatrix}
   =: \bcA_3
 \end{align*}
 jointly as \ $n \to \infty$.
\ Consequently, we obtain an asymptotic expansion
 \begin{align*}
  \bD_n \tbA_n
  &= (n^2 \bD_{n,1} + n \bD_{n,2})
     \begin{bmatrix} 0 & 1 \\ -1 & 0 \end{bmatrix}
     (n^3 \bA_{n,1} + n^2 \bA_{n,2} + n \bA_{n,3})
     \begin{bmatrix} -\bv_\LEFT^\top \\ \bu_\LEFT^\top \end{bmatrix} \\
  &= (n^5 \bC_{n,1} + n^4 \bC_{n,2} + n^3 \bC_{n,3} + n^2 \bC_{n,4})
     \begin{bmatrix} -\bv_\LEFT^\top \\ \bu_\LEFT^\top \end{bmatrix} ,
 \end{align*}
 where
 \[
   \bC_{n,1}
   := \bD_{n,1} \begin{bmatrix} 0 & 1 \\ -1 & 0 \end{bmatrix} \bA_{n,1} = \bzero
   \qquad \text{for all \ $n \in \NN$,}
 \]
 and
 \begin{align*}
  \bC_{n,2}
  &:= \bD_{n,1} \begin{bmatrix} 0 & 1 \\ -1 & 0 \end{bmatrix} \bA_{n,2}
      + \bD_{n,2} \begin{bmatrix} 0 & 1 \\ -1 & 0 \end{bmatrix} \bA_{n,1}
   \distr
   \bcD_1 \begin{bmatrix} 0 & 1 \\ -1 & 0 \end{bmatrix} \bcA_2
   + \bcD_2 \begin{bmatrix} 0 & 1 \\ -1 & 0 \end{bmatrix} \bcA_1 , \\
  \bC_{n,3}
  &:= \bD_{n,1} \begin{bmatrix} 0 & 1 \\ -1 & 0 \end{bmatrix} \bA_{n,3}
      + \bD_{n,2} \begin{bmatrix} 0 & 1 \\ -1 & 0 \end{bmatrix} \bA_{n,2}
   \distr
   \bcD_1 \begin{bmatrix} 0 & 1 \\ -1 & 0 \end{bmatrix} \bcA_3
   + \bcD_2 \begin{bmatrix} 0 & 1 \\ -1 & 0 \end{bmatrix} \bcA_2 , \\
  \bC_{n,4}
  &:= \bD_{n,2} \begin{bmatrix} 0 & 1 \\ -1 & 0 \end{bmatrix} \bA_{n,3}
   \distr \bcD_2 \begin{bmatrix} 0 & 1 \\ -1 & 0 \end{bmatrix} \bcA_3 
 \end{align*}
 as \ $n \to \infty$.
\ Using again Theorem \ref{main_Ad} and \eqref{A_n_deco}, we conclude
 \[
   \begin{bmatrix} n^{-4} \det(\bA_n) \\ n^{-4} \bD_n \tbA_n \end{bmatrix}
   \distr
   \begin{bmatrix} 
    \int_0^1 \cY_t^2 \, \dd t
    \left(\frac{\langle\bv_\LEFT,\bm_\bvare\rangle^2}{(1-\lambda)^2}
          + \frac{\langle\bV_{\!\!\bvare}\bv_\LEFT,\bv_\LEFT\rangle}{1-\lambda^2}\right)
    - \frac{\langle\bv_\LEFT,\bm_\bvare\rangle^2}{(1-\lambda)^2}
    \left(\int_0^1 \cY_t \, \dd t\right)^2 \\[2mm] 
    \bcD_1 \begin{bmatrix} 0 & 1 \\ -1 & 0 \end{bmatrix} \bcA_2
    \begin{bmatrix} -\bv_\LEFT^\top \\ \bu_\LEFT^\top \end{bmatrix}
    + \bcD_2 \begin{bmatrix} 0 & 1 \\ -1 & 0 \end{bmatrix} \bcA_1
      \begin{bmatrix} -\bv_\LEFT^\top \\ \bu_\LEFT^\top \end{bmatrix}
   \end{bmatrix}
 \]
 as \ $n \to \infty$.
\ Here
 \[
   \int_0^1 \cY_t^2 \, \dd t
   \left(\frac{\langle\bv_\LEFT,\bm_\bvare\rangle^2}{(1-\lambda)^2}
         + \frac{\langle\bV_{\!\!\bvare}\bv_\LEFT,\bv_\LEFT\rangle}{1-\lambda^2}\right)
   - \frac{\langle\bv_\LEFT,\bm_\bvare\rangle^2}{(1-\lambda)^2}
   \left(\int_0^1 \cY_t \, \dd t\right)^2
   = \cI_1 + \cI_2 ,
 \]
 \begin{align*}
  \bcD_1 \begin{bmatrix} 0 & 1 \\ -1 & 0 \end{bmatrix} \bcA_2
   \begin{bmatrix} -\bv_\LEFT^\top \\ \bu_\LEFT^\top \end{bmatrix}
  &= \frac{\langle\bv_\LEFT,\bm_\bvare\rangle}{1-\lambda}
     \int_0^1 \cY_t \, \dd t
     \int_0^1 \cY_t \, \dd\bcM_t \, \be_1^\top
     \begin{bmatrix} 0 & 1 \\ -1 & 0 \end{bmatrix}
     \begin{bmatrix} 0 & 1 \\ 1 & 0 \end{bmatrix}
     \begin{bmatrix} -\bv_\LEFT^\top \\ \bu_\LEFT^\top \end{bmatrix} \\
  &= \frac{\langle\bv_\LEFT,\bm_\bvare\rangle}{1-\lambda}
     \int_0^1 \cY_t \, \dd t
     \int_0^1 \cY_t \, \dd\bcM_t \, \bv_\LEFT^\top
 \end{align*}
 and
 \begin{align*}
  &\bcD_2 \begin{bmatrix} 0 & 1 \\ -1 & 0 \end{bmatrix} \bcA_1
   \begin{bmatrix} -\bv_\LEFT^\top \\ \bu_\LEFT^\top \end{bmatrix} \\
  &= \int_0^1 \cY_t^2 \, \dd t
     \left(\frac{\langle\bv_\LEFT,\bm_\bvare\rangle}{1-\lambda} \bcM_1
           + \frac{\langle \bV_{\!\!\bvare} \bv_\LEFT, \bv_\LEFT \rangle^{1/2}}
                  {(1-\lambda^2)^{1/2}}
             \obV_{\!\!\bxi}^{1/2} \!\! \int_0^1 \cY_t \, \dd\tbcW_t\right)
     \be_2^\top
     \begin{bmatrix} 0 & 1 \\ -1 & 0 \end{bmatrix}
     \begin{bmatrix} 1 & 0 \\ 0 & 0 \end{bmatrix}
     \begin{bmatrix} -\bv_\LEFT^\top \\ \bu_\LEFT^\top \end{bmatrix} \\
  &= \int_0^1 \cY_t^2 \, \dd t
     \left(\frac{\langle\bv_\LEFT,\bm_\bvare\rangle}{1-\lambda} \bcM_1
           + \frac{\langle \bV_{\!\!\bvare} \bv_\LEFT, \bv_\LEFT \rangle^{1/2}}
                  {(1-\lambda^2)^{1/2}}
             \obV_{\!\!\bxi}^{1/2} \!\! \int_0^1 \cY_t \, \dd\tbcW_t\right)
     \bv_\LEFT^\top .
 \end{align*}
Since \ $\bm_\bvare \ne \bzero$, \ by the SDE \eqref{SDE_Y}, we have
 \ $\PP(\text{$\cY_t = 0$ \ for all \ $t \in [0,1]$}) = 0$, \ which implies
 that
 \ $\PP\bigl( \int_0^1 \cY_t^2 \, \dd t \int_0^1 \cY_t \, \dd t > 0 \bigr) = 1$,
 \ hence the continuous mapping theorem implies \eqref{m_xi_m}.
\proofend

\section{Proof of Theorem \ref{main_Ad}}
\label{section_proof_main}

The first convergence in Theorem \ref{main_Ad} follows from Lemma
 \ref{limsupUV}.

For the second convergence in Theorem \ref{main_Ad}, consider the sequence of
 stochastic processes
 \[
   \bcZ^{(n)}_t
   := \begin{bmatrix}
       \bcM_t^{(n)} \\
       \bcN_t^{(n)} \\
       \bcP_t^{(n)}
      \end{bmatrix}
   := \sum_{k=1}^\nt
       \bZ^{(n)}_k
   \qquad \text{with} \qquad
   \bZ^{(n)}_k
   := \begin{bmatrix}
       n^{-1} \bM_k \\
       n^{-2} \bM_k U_{k-1} \\
       n^{-3/2} \bM_k V_{k-1}
      \end{bmatrix} 
    = \begin{bmatrix}
       n^{-1} \\
       n^{-2} U_{k-1} \\
       n^{-3/2} V_{k-1}
      \end{bmatrix} 
      \otimes \bM_k
 \]
 for \ $t \in \RR_+$ \ and \ $k, n \in \NN$, \ where \ $\otimes$ \ denotes
 Kronecker product of matrices.
The second convergence in Theorem \ref{main_Ad} follows from Lemma
 \ref{main_VV} and the following theorem (this will be explained after Theorem
 \ref{main_conv}).

\begin{Thm}\label{main_conv}
Suppose that the assumptions of Theorem \ref{main_Ad} hold.
Then we have
 \begin{equation}\label{conv_Z}
   \bcZ^{(n)} \distr \bcZ \qquad \text{as \ $n\to\infty$,}
 \end{equation}
 where the process \ $(\bcZ_t)_{t \in \RR_+}$ \ with values in \ $(\RR^2)^3$ \ is
 the unique strong solution of the SDE
 \begin{equation}\label{ZSDE}
  \dd \bcZ_t
  = \gamma(t, \bcZ_t) \begin{bmatrix} \dd \bcW_t \\ \dd \tbcW_t \end{bmatrix} ,
  \qquad t \in \RR_+ ,
 \end{equation}
 with initial value \ $\bcZ_0 = \bzero$, \ where \ $(\bcW_t)_{t \in \RR_+}$ \ and
 \ $(\tbcW_t)_{t \in \RR_+}$ \ are independent 2-dimensional standard Wiener
 processes, and \ $\gamma : \RR_+ \times (\RR^2)^3 \to (\RR^{2\times2})^{3\times2}$
 \ is defined by
 \[
   \gamma(t, \bx)
   := \begin{bmatrix}
       {(\langle \bu_\LEFT, \bx_1 + t \bm_\bvare \rangle^+)^{1/2}} & 0 \\
       {(\langle \bu_\LEFT, \bx_1 + t \bm_\bvare \rangle^+)^{3/2}} & 0 \\
       0 & \frac{\langle \obV_{\!\!\bxi} \bv_\LEFT, \bv_\LEFT \rangle^{1/2}}
                {(1-\lambda^2)^{1/2}}
           \langle \bu_\LEFT, \bx_1 + t \bm_\bvare \rangle
      \end{bmatrix}
      \otimes \obV_{\!\!\bxi}^{1/2}
 \]
 for \ $t \in \RR_+$ \ and \ $\bx = (\bx_1 , \bx_2 , \bx_3) \in (\RR^2)^3$.
\end{Thm}
(Note that the statement of Theorem \ref{main_conv} holds even if
 \ $\langle \obV_{\!\!\bxi} \bv_\LEFT, \bv_\LEFT \rangle = 0$, \ when the last
 2-dimensional coordinate process of the unique strong solution
 \ $(\bcZ_t)_{t\in\RR_+}$ \ is \ $\bzero$.)

The SDE \eqref{ZSDE} has the form
 \begin{align}\label{MNPSDE}
  \dd \bcZ_t
  = \begin{bmatrix}
     \dd \bcM_t \\
     \dd \bcN_t \\
     \dd \bcP_t
    \end{bmatrix}
  = \begin{bmatrix}
     {(\langle \bu_\LEFT, \bcM_t + t \bm_\bvare \rangle^+)^{1/2}
     \, \obV_{\!\!\bxi}^{1/2} \, \dd \bcW_t} \\[1mm]
     {(\langle \bu_\LEFT, \bcM_t + t \bm_\bvare \rangle^+)^{3/2}
     \, \obV_{\!\!\bxi}^{1/2} \, \dd \bcW_t} \\[1mm]
     \frac{\langle \obV_{\!\!\bxi} \bv_\LEFT, \bv_\LEFT \rangle^{1/2}}
          {(1-\lambda^2)^{1/2}}
     \langle \bu_\LEFT, \bcM_t + t \bm_\bvare \rangle
     \, \obV_{\!\!\bxi}^{1/2} \, \dd \tbcW_t
    \end{bmatrix} , \qquad t\in\RR_+ .
 \end{align}
One can prove that the first 2-dimensional equation of the SDE
 \eqref{MNPSDE} has a pathwise unique strong solution
 \ $(\bcM_t^{(\by_0)})_{t\in\RR_+}$ \ with arbitrary initial value
 \ $\bcM_0^{(\by_0)} = \by_0 \in \RR^2$.
\ Indeed, it is equivalent to the existence of a pathwise unique strong
 solution of the SDE
 \begin{equation}\label{SDE_P_Q}
  \begin{cases}
   \dd \cS_t
   = \langle \bu_\LEFT, \bm_\bvare \rangle \, \dd t
     + (\cS_t^+)^{1/2} \, \bu_\LEFT^\top \obV_{\!\!\bxi}^{1/2} \, \dd \bcW_t , \\[2mm]
   \dd \bcQ_t
   = - \bPi \bm_\bvare \, \dd t
     + (\cS_t^+)^{1/2} \, \bigl( \bI_2 - \bPi \bigr) \obV_{\!\!\bxi}^{1/2}
       \, \dd \bcW_t ,
  \end{cases}
  \qquad t \in \RR_+ ,
 \end{equation}
 with initial value
 \ $\bigl(\cS_0^{(\by_0)}, \, \bcQ_0^{(\by_0)}\bigr)
    =\big(\langle \bu_\LEFT, \by_0 \rangle, \, (\bI_2 - \bPi)\by_0\bigr)
    \in \RR \times \RR^2$,
 \ where \ $\bI_2$ \ denotes the $2$-dimensional unit matrix and
 \ $\bPi := \bu_\RIGHT \bu_\LEFT^\top$, \ since we have the correspondences
 \begin{gather*}
  \cS_t^{(\by_0)} = \bu_\LEFT^\top (\bcM_t^{(\by_0)} + t \bm_\bvare) , \qquad
  \bcQ_t^{(\by_0)} = \bcM_t^{(\by_0)} - \cS_t^{(\by_0)} \bu_\RIGHT \\
  \bcM_t^{(\by_0)} = \bcQ_t^{(\by_0)} + \cS_t^{(\by_0)} \bu_\RIGHT ,
 \end{gather*}
 see the proof of Isp\'any and Pap \cite[Theorem 3.1]{IspPap2}.
By Remark \ref{REMARK_SDE}, \ $\cS_t^+$ \  may be replaced by
 \ $\cS_t$ \ for all \ $t \in \RR_+$ \ in the first equation of
 \eqref{SDE_P_Q} provided that \ $\langle \bu_\LEFT, \by_0 \rangle \in \RR_+$,
 \ hence \ $\langle \bu_\LEFT, \bcM_t + t \bm_\bvare \rangle^+$ \ may be
 replaced by \ $\langle \bu_\LEFT, \bcM_t + t \bm_\bvare \rangle$ \ for
 all \ $t \in \RR_+$ \ in \eqref{MNPSDE}.
Thus the SDE \eqref{ZSDE} has a pathwise unique strong solution with initial
 value \ $\bcZ_0 = \bzero$, \ and we have
 \[
   \bcZ_t
   = \begin{bmatrix}
      \bcM_t \\
      \bcN_t \\
      \bcP_t
     \end{bmatrix}
   = \begin{bmatrix}
      \int_0^t
       \langle \bu_\LEFT, \bcM_s + s \bm_\bvare \rangle^{1/2} \, \obV_{\!\!\bxi}^{1/2}
       \, \dd \bcW_s \\[1mm]
      \int_0^t
       \langle \bu_\LEFT, \bcM_s + s \bm_\bvare \rangle \, \dd \bcM_s \\[1mm]
       \frac{\langle \obV_{\!\!\bxi} \bv_\LEFT, \bv_\LEFT \rangle^{1/2}}
            {(1-\lambda^2)^{1/2}}
      \int_0^t
       \langle \bu_\LEFT, \bcM_s + s \bm_\bvare \rangle \, \obV_{\!\!\bxi}^{1/2}
       \, \dd \tbcW_s
     \end{bmatrix} , \qquad t\in\RR_+ .
 \]
By the method of the proof of \ $\cX^{(n)} \distr \cX$ \ in Theorem 3.1 in
 Barczy et al.\ \cite{BarIspPap0}, applying Lemma \ref{Conv2Funct}, one can
 easily derive
 \begin{align}\label{convXZ}
  \begin{bmatrix} \bcX^{(n)} \\ \bcZ^{(n)} \end{bmatrix}
  \distr \begin{bmatrix} \tbcX \\ \bcZ \end{bmatrix} \qquad
  \text{as \ $n \to \infty$,}
 \end{align}
 where
 \[
   \bcX^{(n)}_t := n^{-1} \bX_\nt , \qquad
   \tbcX_t: = \langle \bu_\LEFT, \bcM_t + t \bm_\bvare \rangle \bu_\RIGHT ,
   \qquad t \in \RR_+ , \qquad n \in \NN .
 \]
More precisely, using that
 \[
   \bX_k = \sum_{j=1}^k \bm_\bxi^{k-j} (\bM_j + \bm_\bvare) , \qquad k \in \NN ,
 \]
 we have
 \[
   \begin{bmatrix}
    \bcX^{(n)} \\
    \bcZ^{(n)} \\
   \end{bmatrix}
   = \psi_n(\bcZ^{(n)}), \qquad n \in \NN ,
 \]
 where the mapping \ $\psi_n : \DD(\RR_+,(\RR^2)^3) \to \DD(\RR_+,(\RR^2)^4)$
 \ is given by
 \[
   \psi_n(f_1,f_2,f_3)(t)
   := \begin{bmatrix}
       \sum_{j=1}^\nt
        \bm_\bxi^{\nt-j}
        \left(f_1\left(\frac{j}{n}\right) - f_1\left(\frac{j-1}{n}\right)
              + \frac{\bm_\bvare}{n}\right) \\
       f_1(t) \\
       f_2(t) \\
       f_3(t) \\
      \end{bmatrix}
 \]
 for \ $f_1, f_2, f_3 \in \DD(\RR_+,\RR^2)$, \ $t \in \RR_+$, \ $n \in \NN$.
\ Further, we have
 \[
   \begin{bmatrix}
    \tbcX \\
    \bcZ \\
   \end{bmatrix}
    = \psi(\bcZ) ,
 \]
 where the mapping \ $\psi : \DD(\RR_+,(\RR^2)^3) \to \DD(\RR_+,(\RR^2)^4)$ \ is
 given by
 \[
   \psi(f_1,f_2,f_3)(t)
   := \begin{bmatrix}
       \langle \bu_\LEFT, f_1(t) + t \bm_\bvare \rangle \bu_\RIGHT \\
       f_1(t) \\
       f_2(t) \\
       f_3(t) \\
      \end{bmatrix}
 \]
 for \ $f_1, f_2, f_3 \in \DD(\RR_+,\RR^2)$ \ and \ $t \in \RR_+$.
\ By page 603 in Barczy et al. \cite{BarIspPap0}, the mappings \ $\psi_n$,
 \ $n \in \NN$, \ and \ $\psi$ \ are measurable (the latter one is continuous
 too), since the coordinate functions are measurable.
Using page 604 in Barczy et al. \cite{BarIspPap0}, we obtain that the set
 \[
   C := \big\{f \in \CC(\RR_+,(\RR^2)^3) : f(0) = \bzero \in (\RR^2)^3 \big\}
 \]
 has the properties \ $C \subseteq C_{\psi,(\psi_n)_{n\in\NN}}$ \ with
 \ $C \in \cB(\DD(\RR_+,(\RR^2)^3))$ \ and \ $\PP(\bcZ\in C) = 1$, \ where
 \ $C_{\psi,(\psi_n)_{n\in\NN}}$ \ is defined in Appendix \ref{app_C}.
Hence, by \eqref{conv_Z} and Lemma \ref{Conv2Funct}, we have
 \[
   \begin{bmatrix}
    \bcX^{(n)} \\
    \bcZ^{(n)} \\
   \end{bmatrix}
   =\psi_n(\bcZ^{(n)})
   \distr
   \psi(\bcZ)
   =\begin{bmatrix}
     \tbcX \\
     \bcZ \\
    \end{bmatrix}
   \qquad \text{as \ $n \to \infty$,}
 \]
 as desired.

Next, similarly to the proof of \eqref{seged2}, by Lemma \ref{Marci},
 convergence \eqref{convXZ} with
 \ $U_{k-1} = \langle \bu_\LEFT, \bX_{k-1} \rangle$ \ and Lemma \ref{main_VV}
 implies
 \[
   \sum_{k=1}^n
    \begin{bmatrix}
     n^{-3} U_{k-1}^2 \\
     n^{-2} V_{k-1}^2 \\
     n^{-2} \bM_k U_{k-1} \\
     n^{-3/2} \bM_k V_{k-1}
    \end{bmatrix}
   \distr \begin{bmatrix}
           \int_0^1 \langle \bu_\LEFT, \tbcX_t \rangle^2 \, \dd t \\[1mm]
           \frac{\langle \obV_{\!\!\bxi} \bv_\LEFT, \bv_\LEFT \rangle}
                {1-\lambda^2}
           \int_0^1 \langle \bu_\LEFT, \tbcX_t \rangle \, \dd t \\[1mm]
           \int_0^1 \cY_t \, \dd \bcM_t \\[1mm]
           \frac{\langle \obV_{\!\!\bxi} \bv_\LEFT, \bv_\LEFT \rangle^{1/2}}
                {(1-\lambda^2)^{1/2}}
           \int_0^1 \cY_t \obV_{\!\!\bxi}^{1/2} \, \dd \tbcW_t
          \end{bmatrix}
  \qquad \text{as \ $n \to \infty$.}
 \]
This limiting random vector can be written in the form as given in Theorem
 \ref{main_Ad}, since \ $\langle \bu_\LEFT, \tbcX_t \rangle = \cY_t$
 \ for all \ $t \in \RR_+$.

\noindent
\textbf{Proof of Theorem \ref{main_conv}.}
In order to show convergence \ $\bcZ^{(n)} \distr \bcZ$, \ we apply Theorem
 \ref{Conv2DiffThm} with the special choices \ $\bcU := \bcZ$,
 \ $\bU^{(n)}_k := \bZ^{(n)}_k$, \ $n, k \in \NN$,
 \ $(\cF_k^{(n)})_{k\in\ZZ_+} := (\cF_k)_{k\in\ZZ_+}$ \ and the function \ $\gamma$
 \ which is defined in Theorem \ref{main_conv}.
Note that the discussion after Theorem \ref{main_conv} shows that the SDE
 \eqref{ZSDE} admits a unique strong solution \ $(\bcZ_t^\bz)_{t\in\RR_+}$ \ for
 all initial values \ $\bcZ_0^\bz = \bz \in (\RR^2)^3$.

Now we show that conditions (i) and (ii) of Theorem \ref{Conv2DiffThm} hold.
The conditional variance has the form
 \[
   \var\bigl(\bZ^{(n)}_k \mid \cF_{k-1}\bigr) =
   \begin{bmatrix}
    n^{-2}
    & n^{-3} U_{k-1}
    & n^{-5/2} V_{k-1} \\
    n^{-3} U_{k-1}
    & n^{-4} U_{k-1}^2
    & n^{-7/2} U_{k-1} V_{k-1} \\
    n^{-5/2} V_{k-1}
    & n^{-7/2} U_{k-1} V_{k-1}
    & n^{-3} V_{k-1}^2
   \end{bmatrix}
   \otimes \bV_{\!\!\bM_k}
 \]
 for \ $n \in \NN$, \ $k \in \{1, \ldots, n\}$, \ with
 \ $\bV_{\!\!\bM_k} := \var(\bM_k \mid \cF_{k-1})$, \ and
 \ $\gamma(s,\bcZ_s^{(n)}) \gamma(s,\bcZ_s^{(n)})^\top$ \ has the form
 \[
   \begin{bmatrix}
    \langle \bu_\LEFT, \bcM_s^{(n)} + s \bm_\bvare \rangle
     & \langle \bu_\LEFT, \bcM_s^{(n)} + s \bm_\bvare \rangle^2
     & \bzero \\
    \langle \bu_\LEFT, \bcM_s^{(n)} + s \bm_\bvare \rangle^2
     & \langle \bu_\LEFT, \bcM_s^{(n)} + s \bm_\bvare \rangle^3
     & \bzero \\
    \bzero
     & \bzero
     & \frac{\langle \obV_{\!\!\bxi} \bv_\LEFT, \bv_\LEFT \rangle}
            {1-\lambda^2}
       \langle \bu_\LEFT, \bcM_s^{(n)} + s \bm_\bvare \rangle^2
   \end{bmatrix}
   \otimes \obV_{\!\!\xi}
 \]
 for \ $s\in\RR_+$, \ where we used that
 \ $\langle \bu_\LEFT, \bcM_s^{(n)} + s \bm_\bvare \rangle^+
    = \langle \bu_\LEFT, \bcM_s^{(n)} + s \bm_\bvare \rangle$, \ $s \in \RR_+$,
 \ $n \in \NN$.
\ Indeed, by \eqref{Mk}, we get
 \begin{align}\label{M+}
  \begin{split}
   &\langle \bu_\LEFT, \bcM_s^{(n)} + s \bm_\bvare \rangle
    = \frac{1}{n}
      \sum_{k=1}^\ns
       \langle \bu_\LEFT, \bX_k - \bm_\bxi \bX_{k-1} - \bm_\bvare \rangle
      + \langle \bu_\LEFT, s \bm_\bvare \rangle \\
   &= \frac{1}{n}
      \sum_{k=1}^\ns
       \langle \bu_\LEFT, \bX_k - \bX_{k-1} - \bm_\bvare \rangle
      + s \langle \bu_\LEFT, \bm_\bvare \rangle \\
   &= \frac{1}{n}
      \langle \bu_\LEFT, \bX_\ns \rangle
      + \frac{ns - \ns}{n} \langle \bu_\LEFT, \bm_\bvare \rangle
    = \frac{1}{n} U_\ns
      + \frac{ns - \ns}{n} \langle \bu_\LEFT, \bm_\bvare \rangle
    \in \RR_+
  \end{split}
 \end{align}
 for \ $s \in \RR_+$, \ $n \in \NN$, \ since
 \ $\bu_\LEFT^\top \bm_\bxi = \bu_\LEFT^\top$ \ implies
 \ $\langle \bu_\LEFT, \bm_\bxi \bX_{k-1} \rangle = \bu_\LEFT^\top \bm_\bxi \bX_{k-1}
    = \bu_\LEFT^\top \bX_{k-1} = \langle \bu_\LEFT, \bX_{k-1} \rangle$.

In order to check condition (i) of Theorem \ref{Conv2DiffThm}, we need  to
 prove that for each \ $T > 0$,
 \begin{gather}
  \sup_{t\in[0,T]}
   \bigg\| \frac{1}{n^2} \sum_{k=1}^{\nt} \bV_{\!\!\bM_k}
           - \int_0^t
              \langle \bu_\LEFT, \bcM_s^{(n)} + s \bm_\bvare \rangle
              \, \obV_{\!\!\bxi} \,
              \dd s \bigg\|
  \stoch 0 , \label{Zcond1} \\
  \sup_{t\in[0,T]}
   \bigg\| \frac{1}{n^3}
           \sum_{k=1}^{\nt} U_{k-1} \bV_{\!\!\bM_k}
           - \int_0^t
              \langle \bu_\LEFT, \bcM_s^{(n)} + s \bm_\bvare \rangle^2 \,
              \obV_{\!\!\bxi} \,
              \dd s \bigg\|
  \stoch 0 , \label{Zcond2} \\
  \sup_{t\in[0,T]}
   \bigg\| \frac{1}{n^4}
           \sum_{k=1}^{\nt} U_{k-1}^2 \bV_{\!\!\bM_k}
           - \int_0^t
              \langle \bu_\LEFT, \bcM_s^{(n)} + s \bm_\bvare \rangle^3 \,
              \obV_{\!\!\bxi} \,
              \dd s \bigg\|
  \stoch 0 , \label{Zcond3} \\
  \sup_{t\in[0,T]}
   \bigg\| \frac{1}{n^3}
           \sum_{k=1}^{\nt} V_{k-1}^2 \bV_{\!\!\bM_k}
           - \frac{\langle \obV_{\!\!\bxi} \bv_\LEFT, \bv_\LEFT \rangle}
                  {1-\lambda^2}
             \int_0^t
              \langle \bu_\LEFT, \bcM_s^{(n)} + s \bm_\bvare \rangle^2 \,
              \obV_{\!\!\bxi} \,
              \dd s \bigg\|
  \stoch 0 , \label{Zcond4} \\
  \sup_{t\in[0,T]}
   \bigg\| \frac{1}{n^{5/2}}
           \sum_{k=1}^{\nt} V_{k-1} \bV_{\!\!\bM_k} \bigg\|
  \stoch 0 , \label{Zcond5} \\
  \sup_{t\in[0,T]}
   \bigg\| \frac{1}{n^{7/2}}
           \sum_{k=1}^{\nt}
            U_{k-1} V_{k-1} \bV_{\!\!\bM_k} \bigg\|
  \stoch 0  \label{Zcond6}
 \end{gather}
 as \ $n \to \infty$.

First we show \eqref{Zcond1}.
By \eqref{M+},
 \[
   \int_0^t \langle \bu_\LEFT, \bcM_s^{(n)} + s \bm_\bvare \rangle \, \dd s
   = \frac{1}{n^2} \sum_{k=1}^{\nt-1} U_k
     + \frac{nt - \nt}{n^2} U_\nt
     + \frac{\nt + (nt - \nt)^2}{2 n^2} \langle \bu_\LEFT, \bm_\bvare \rangle .
 \]
Using Lemma \ref{Moments}, we have
 \ $\bV_{\!\!\bM_k} = U_{k-1} \obV_{\!\!\bxi} + V_{k-1} \tbV_\bxi + \bV_{\!\!\bvare}$,
 \ thus, in order to show \eqref{Zcond1}, it suffices to prove
 \begin{gather}\label{1supsumV_supU}
  n^{-2} \sum_{k=1}^{\nT} |V_k| \stoch 0 , \qquad
  n^{-2} \sup_{t \in [0,T]} U_\nt \stoch 0 , \\
  n^{-2} \sup_{t \in [0,T]} \left[ \nt + (nt - \nt)^2 \right] \to 0 \label{1supnt}
 \end{gather}
 as \ $n \to \infty$.
\ Using \eqref{seged_UV_UNIFORM1} with \ $(\ell, i, j) = (2, 0, 1)$ \ and
 \eqref{seged_UV_UNIFORM2} with \ $(\ell, i , j) = (2, 1, 0)$, \ we have
 \eqref{1supsumV_supU}.
Clearly, \eqref{1supnt} follows from \ $|nt - \nt| \leq 1$, \ $n \in \NN$,
 \ $t \in \RR_+$, \ thus we conclude \eqref{Zcond1}.

Next we turn to prove \eqref{Zcond2}.
By \eqref{M+},
 \begin{align*}
  \int_0^t \langle \bu_\LEFT, \bcM_s^{(n)} + s \bm_\bvare \rangle^2 \, \dd s
  &= \frac{1}{n^3} \sum_{k=1}^{\nt-1} U_k^2
     + \frac{1}{n^3} \langle \bu_\LEFT, \bm_\bvare \rangle \sum_{k=1}^{\nt-1} U_k
     + \frac{nt - \nt}{n^3} U_\nt^2 \\
  &\phantom{\quad}
     + \frac{(nt - \nt)^2}{n^3} \langle \bu_\LEFT, \bm_\bvare \rangle U_\nt
     + \frac{\nt + (nt - \nt)^3}{3n^3} \langle \bu_\LEFT, \bm_\bvare \rangle^2 .
 \end{align*}
Using Lemma \ref{Moments}, we obtain
 \begin{align}\label{UM2F}
  \sum_{k=1}^{\nt} U_{k-1} \bV_{\!\!\bM_k}
  = \sum_{k=1}^{\nt} U_{k-1}^2 \obV_{\!\!\bxi} 
    + \sum_{k=1}^{\nt} U_{k-1} V_{k-1} \tbV_{\!\!\bxi}
    + \sum_{k=1}^{\nt} U_{k-1} \bV_{\!\!\bvare} .
 \end{align}
Thus, in order to show \eqref{Zcond2}, it suffices to prove
 \begin{gather}\label{2supsumUV_2supsumU_2supU}
  n^{-3} \sum_{k=1}^{\nT} |U_k V_k| \stoch 0 , \qquad
  n^{-3} \sum_{k=1}^{\nT} U_k \stoch 0 , \qquad
  n^{-3/2} \sup_{t \in [0,T]} U_\nt \stoch 0 , \\
  n^{-3} \sup_{t \in [0,T]} \left[ \nt + (nt - \nt)^3 \right] \to 0
    \label{2supnt}
 \end{gather}
 as \ $n \to \infty$.
\ By \eqref{seged_UV_UNIFORM1} with \ $(\ell, i, j) = (2, 1, 1)$ \ and
 \ $(\ell, i, j) = (2, 1, 0)$, \ and by \eqref{seged_UV_UNIFORM5}, we have
 \eqref{2supsumUV_2supsumU_2supU}.
Clearly, \eqref{2supnt} follows from \ $|nt - \nt| \leq 1$, \ $n \in \NN$,
 \ $t \in \RR_+$, \ thus we conclude \eqref{Zcond2}.

Now we turn to check \eqref{Zcond3}.
Again by \eqref{M+}, we have
 \begin{align*}
  \int_0^t \langle \bu_\LEFT, \bcM_s^{(n)} + s \bm_\bvare \rangle^3 \, \dd s
  &= \frac{1}{n^4} \sum_{k=1}^{\nt-1} U_k^3
     + \frac{3}{2n^4} \langle \bu_\LEFT, \bm_\bvare \rangle \sum_{k=1}^{\nt-1} U_k^2
     + \frac{1}{n^4} \langle \bu_\LEFT, \bm_\bvare \rangle^2
       \sum_{k=1}^{\nt-1} U_k \\
  &\phantom{\quad}
     + \frac{nt - \nt}{n^4} U_\nt^3
     + \frac{3 (nt - \nt)^2}{2n^4} \langle \bu_\LEFT, \bm_\bvare \rangle U_\nt^2 \\
  &\phantom{\quad}
     + \frac{(nt - \nt)^3}{n^4} \langle \bu_\LEFT, \bm_\bvare \rangle^2 \, U_\nt
     + \frac{\nt + (nt - \nt)^4}{4n^4} \langle \bu_\LEFT, \bm_\bvare \rangle^3 .
 \end{align*}
Using Lemma \ref{Moments}, we obtain
 \begin{equation}\label{U2M2F}
  \sum_{k=1}^{\nt} U_{k-1}^2 \bV_{\!\!\bM_k}
  = \sum_{k=1}^{\nt} U_{k-1}^3 \obV_{\!\!\bxi} 
    + \sum_{k=1}^{\nt} U_{k-1}^2 V_{k-1} \tbV_{\!\!\bxi}
    + \sum_{k=1}^{\nt} U_{k-1}^2 \bV_{\!\!\bvare} .
 \end{equation}
Thus, in order to show \eqref{Zcond3}, it suffices to prove
 \begin{gather}
  n^{-4} \sum_{k=1}^{\nT} | U_k^2 V_k | \stoch 0 , \qquad
  n^{-4} \sum_{k=1}^{\nT} U_k^2 \stoch 0 , \label{3supsumUUV_3supsumUU} \\
  n^{-4} \sum_{k=1}^{\nT} U_k \stoch 0 , \qquad
  n^{-4/3} \sup_{t \in [0,T]} U_\nt \stoch 0 , \label{3supsumU_3supU} \\
  n^{-4} \sup_{t \in [0,T]} \left[ \nt + (nt - \nt)^4 \right] \to 0
    \label{3supnt}
 \end{gather}
 as \ $n \to \infty$.
\ By \eqref{seged_UV_UNIFORM1} with \ $(\ell, i , j) = (4, 2, 1)$,
 \ $(\ell, i, j) = (4, 2, 0)$, \ and \ $(\ell, i , j) = (2, 1, 0)$, \ and by
 \eqref{seged_UV_UNIFORM5}, \ we have \eqref{3supsumUUV_3supsumUU} and
 \eqref{3supsumU_3supU}.
Clearly, \eqref{3supnt} follows again from \ $|nt - \nt| \leq 1$,
 \ $n \in \NN$, \ $t \in \RR_+$, \ thus we conclude \eqref{Zcond3}.

Next we turn to prove \eqref{Zcond4}.
First we show that
 \begin{align}\label{Zcond4a}
  n^{-3}
  \sup_{t \in [0,T]}
   \left\| \sum_{k=1}^{\nt} V_{k-1}^2 \bV_{\!\!\bM_k}
           - \frac{\langle \obV_{\!\!\bxi} \bv_\LEFT, \bv_\LEFT \rangle}
                  {1-\lambda^2}
             \sum_{k=1}^{\nt} U_{k-1}^2 \obV_{\!\!\bxi} \right\|
  \stoch 0
 \end{align}
 as \ $n \to \infty$ \ for all \ $T > 0$.
\ Using Lemma \ref{Moments}, we obtain
 \begin{align}\label{V2M2F}
  \sum_{k=1}^{\nt} V_{k-1}^2 \bV_{\!\!\bM_k}
  = \sum_{k=1}^{\nt} U_{k-1} V_{k-1}^2 \obV_{\!\!\bxi} 
    + \sum_{k=1}^{\nt} V_{k-1}^3 \tbV_{\!\!\bxi}
    + \sum_{k=1}^{\nt} V_{k-1}^2 \bV_{\!\!\bvare} .
 \end{align}
Using \eqref{seged_UV_UNIFORM1} with \ $(\ell, i , j) = (6, 0, 3)$ \ and
 \ $(\ell, i , j) = (4, 0, 2)$, \ we have
 \begin{align*}
  n^{-3} \sum_{k=1}^{\nT} | V_k |^3 \stoch 0 , \qquad
  n^{-3} \sum_{k=1}^{\nT} V_k^2 \stoch 0 \qquad \text{as \ $n\to\infty$,}
 \end{align*}
 hence \eqref{Zcond4a} will follow from
 \begin{align}\label{Zcond4b}
  n^{-3}
  \sup_{t \in [0,T]}
   \left\| \sum_{k=1}^{\nt} U_{k-1} V_{k-1}^2
           - \frac{\langle \obV_{\!\!\bxi} \bv_\LEFT, \bv_\LEFT \rangle}
                  {1-\lambda^2}
             \sum_{k=1}^{\nt} U_{k-1}^2 \right\|
  \stoch 0 \qquad \text{as \ $n \to \infty$}
 \end{align}
 for all \ $T > 0$.
\ The aim of the following discussion is to decompose
 \ $\sum_{k=1}^{\nt} U_{k-1} V_{k-1}^2$ \ as a sum of a martingale and some
 other terms.
Using recursions  \eqref{rec_V}, \eqref{rec_U} and formulas \eqref{Mcond} and
 \eqref{M3cond}, we obtain
 \begin{align*}
  \EE(U_{k-1} V_{k-1}^2 \mid \cF_{k-2})
  &= \EE\Bigl((U_{k-2} + \langle \bu_\LEFT, \bM_{k-1} + \bm_\bvare \rangle)
              \bigl(\lambda V_{k-2}
                    + \langle \bv_\LEFT, \bM_{k-1} + \bm_\bvare \rangle\bigr)^2
              \,\Big|\, \cF_{k-2}\Bigr) \\
  &= \lambda^2 U_{k-2} V_{k-2}^2
     + \bv_\LEFT^\top \EE(\bM_{k-1} \bM_{k-1}^\top \mid \cF_{k-2}) \bv_\LEFT \,
       U_{k-2} \\
  &\quad 
     + \text{constant}
     + \text{linear combination of \ $U_{k-2} V_{k-2}$, \ $V_{k-2}^2$, \ $U_{k-2}$
             \ and \ $V_{k-2}$} \\
  &= \lambda^2 U_{k-2} V_{k-2}^2
     + \langle \obV_{\!\!\bxi} \bv_\LEFT, \bv_\LEFT \rangle \, U_{k-2}^2
     + \text{constant} \\
  &\quad 
     + \text{linear combination of \ $U_{k-2} V_{k-2}$, \ $V_{k-2}^2$, \ $U_{k-2}$
             \ and \ $V_{k-2}$.}
 \end{align*}
Thus
 \begin{align*}
  &\sum_{k=1}^{\nt} U_{k-1} V_{k-1}^2
   = \sum_{k=2}^{\nt}
      \big[U_{k-1} V_{k-1}^2 - \EE(U_{k-1} V_{k-1}^2 \mid \cF_{k-2}) \big]
     + \sum_{k=2}^{\nt} \EE(U_{k-1} V_{k-1}^2 \mid \cF_{k-2}) \\
  &= \sum_{k=2}^{\nt}
      \big[U_{k-1} V_{k-1}^2 - \EE(U_{k-1} V_{k-1}^2 \mid \cF_{k-2}) \big]
     + \lambda^2 \sum_{k=2}^{\nt} U_{k-2} V_{k-2}^2
     + \langle \obV_{\!\!\bxi} \bv_\LEFT, \bv_\LEFT \rangle \,
       \sum_{k=2}^{\nt} U_{k-2}^2 \\
  &\quad
     + \OO(n)
     + \text{linear combination of \ $\sum_{k=2}^{\nt} U_{k-2} V_{k-2}$,
             \ $\sum_{k=2}^{\nt} V_{k-2}^2$, \ $\sum_{k=2}^{\nt} U_{k-2}$
             \ and \ $\sum_{k=2}^{\nt} V_{k-2}$.}
 \end{align*}
Consequently,
 \begin{align*}
  \sum_{k=1}^{\nt} U_{k-1} V_{k-1}^2
  &= \frac{1}{1-\lambda^2}
     \sum_{k=2}^{\nt}
      \big[U_{k-1} V_{k-1}^2 - \EE(U_{k-1} V_{k-1}^2 \mid \cF_{k-2}) \big] \\
  &\quad
     + \frac{\langle \obV_{\!\!\bxi} \bv_\LEFT, \bv_\LEFT \rangle}
            {1-\lambda^2}
       \sum_{k=2}^{\nt} U_{k-2}^2
     - \frac{\lambda^2}{1-\lambda^2}
       U_{\nt - 1} V_{\nt - 1}^2
     + \OO(n) \\
  &\quad
     + \text{linear combination of \ $\sum_{k=2}^{\nt} U_{k-2} V_{k-2}$,
             \ $\sum_{k=2}^{\nt} V_{k-2}^2$, \ $\sum_{k=2}^{\nt} U_{k-2}$
             \ and \ $\sum_{k=2}^{\nt} V_{k-2}$.}
 \end{align*}
Using \eqref{seged_UV_UNIFORM4} with \ $(\ell, i , j) = (8, 1, 2)$ \ we have
 \begin{align*}
  n^{-3} \sup_{t \in [0,T]} \,
         \Biggl| \sum_{k=2}^\nt
                  \big[ U_{k-1} V_{k-1}^2
                        - \EE(U_{k-1} V_{k-1}^2 \mid \cF_{k-2}) \big] \Biggr|
  \stoch 0 \qquad \text{as \ $n\to\infty$.}
 \end{align*}
Thus, in order to show \eqref{Zcond4b}, it suffices to prove
 \begin{gather}
  n^{-3} \sum_{k=1}^{\nT} | U_k V_k | \stoch 0 , \qquad
  n^{-3} \sum_{k=1}^{\nT} V_k^2 \stoch 0 , \label{4supsumUV_4supsumVV} \\
  n^{-3} \sum_{k=1}^{\nT} U_k \stoch 0 , \qquad
  n^{-3} \sum_{k=1}^{\nT} |V_k| \stoch 0 , \label{4supsumU_4supsumV} \\
  n^{-3} \sup_{t \in [0,T]} U_\nt V_\nt^2 \stoch 0 , \qquad
  n^{-3/2} \sup_{t \in [0,T]} U_{\nt} \stoch 0 \label{4supUVV_4supUU}
   \end{gather}
 as \ $n \to \infty$.
\ Using \eqref{seged_UV_UNIFORM1} with \ $(\ell, i , j) = (2, 1, 1)$;
 \ $(\ell, i , j) = (4, 0, 2)$; \ $(\ell, i , j) = (2, 1, 0)$, \ and
 \ $(\ell, i , j) = (2, 0, 1)$, \ we have \eqref{4supsumUV_4supsumVV} and
 \eqref{4supsumU_4supsumV}.
By \eqref{seged_UV_UNIFORM2} with \ $(\ell, i , j) = (4, 1, 2)$, \ and by
 \eqref{seged_UV_UNIFORM5}, we have \eqref{4supUVV_4supUU}.
Thus we conclude \eqref{Zcond4b}, and hence \eqref{Zcond4a}.
By Lemma \ref{Moments} and \eqref{seged_UV_UNIFORM1} with
 \ $(\ell, i , j) = (2, 1, 1)$ \ and \ $(\ell, i , j) = (2, 1, 0)$, we get
 \begin{align}\label{Zcond2a}
  n^{-3}
  \sup_{t \in [0,T]}
   \left\| \sum_{k=1}^{\nt} U_{k-1} \bV_{\!\!\bM_k}
           - \sum_{k=1}^{\nt} U_{k-1}^2 \obV_{\!\!\bxi} \right\|
  \stoch 0
 \end{align}
 as \ $n \to \infty$ \ for all \ $T > 0$.
\ As a last step, using \eqref{Zcond2}, we obtain \eqref{Zcond4}.

For \eqref{Zcond5}, consider
 \begin{align}\label{VM2F}
  \sum_{k=1}^{\nt} V_{k-1} \bV_{\!\!\bM_k}
  = \sum_{k=1}^{\nt} U_{k-1} V_{k-1} \obV_{\!\!\bxi} 
    + \sum_{k=1}^{\nt} V_{k-1}^2 \tbV_{\!\!\bxi}
    + \sum_{k=1}^{\nt} V_{k-1} \bV_{\!\!\bvare}  ,
 \end{align}
 where we used Lemma \ref{Moments}.
Using \eqref{seged_UV_UNIFORM1} with \ $(\ell, i , j) = (4, 0, 2)$, \ and
 \ $(\ell, i , j) = (2, 0, 1)$, \ we have
 \begin{align*}
  n^{-5/2} \sum_{k=1}^{\nT} V_k^2 \stoch 0 , \qquad
  n^{-5/2} \sum_{k=1}^{\nT} |V_k| \stoch 0 \qquad \text{as \ $n \to \infty$,}
 \end{align*}
 hence \eqref{Zcond5} follows from Lemma \ref{limsupUV}.

Convergence \eqref{Zcond6} can be handled in the same way as \eqref{Zcond5}.
For completeness we present all of the details.
By Lemma \ref{Moments}, we have
 \begin{equation}\label{UVM2F}
  \sum_{k=1}^{\nt} U_{k-1} V_{k-1} \bV_{\!\!\bM_k}
  = \sum_{k=1}^{\nt} U_{k-1}^2 V_{k-1} \obV_{\!\!\bxi} 
    + \sum_{k=1}^{\nt} U_{k-1} V_{k-1}^2 \tbV_{\!\!\bxi}
    + \sum_{k=1}^{\nt} U_{k-1} V_{k-1} \bV_{\!\!\bvare}  .
 \end{equation}
Using \eqref{seged_UV_UNIFORM1} with \ $(\ell, i , j) = (4, 1, 2)$, \ and
 \ $(\ell, i , j) =(2, 1, 1)$, \ we have
 \begin{align*}
  n^{-7/2} \sum_{k=1}^{\nT} U_{k-1}V_{k-1}^2 \stoch 0 , \qquad
  n^{-7/2} \sum_{k=1}^{\nT} |U_{k-1} V_{k-1}| \stoch 0 \qquad
  \text{as \ $n \to \infty$,}
 \end{align*}
 hence \eqref{Zcond6} will follow from
 \begin{align}\label{Zcond6a}
  n^{-7/2} \sup_{t \in [0,T]} \left| \sum_{k=1}^{\nt} U_{k-1}^2 V_{k-1} \right|
  \stoch 0 \qquad \text{as \ $n \to \infty$.}
 \end{align}
The aim of the following discussion is to decompose
 \ $\sum_{k=1}^{\nt} U_{k-1}^2 V_{k-1}$ \ as a sum of a martingale and some
 other terms.
Using recursions \eqref{rec_V}, \eqref{rec_U} and Lemma \ref{Moments}, we
 obtain
 \begin{align*}
  &\EE(U_{k-1}^2V_{k-1} \mid \cF_{k-2})
   = \EE\Bigl((U_{k-2} + \langle \bu_\LEFT, \bM_{k-1} + \bm_\bvare \rangle)^2
              \bigl(\lambda V_{k-2}
                    + \langle \bv_\LEFT, \bM_{k-1} + \bm_\bvare \rangle\bigr)
              \,\Big|\, \cF_{k-2}\Bigr) \\
  &\qquad
   = \lambda U_{k-2}^2 V_{k-2} + \text{constant}  
     + \text{linear combination of \ $U_{k-2}$, $V_{k-2}$, $U_{k-2}^2$, $V_{k-2}^2$
             \ and \ $U_{k-2}V_{k-2}$.}
\end{align*}
Thus
 \begin{multline*}
  \sum_{k=1}^{\nt} U_{k-1}^2 V_{k-1}
   = \sum_{k=2}^{\nt}
      \big[ U_{k-1}^2V_{k-1} - \EE(U_{k-1}^2V_{k-1} \mid \cF_{k-2}) \big]
     + \sum_{k=2}^{\nt} \EE(U_{k-1}^2V_{k-1} \mid \cF_{k-2}) \\
  \begin{aligned}
   &= \sum_{k=2}^{\nt}
       \big[ U_{k-1}^2 V_{k-1} - \EE(U_{k-1}^2V_{k-1} \mid \cF_{k-2}) \big]
      + \lambda \sum_{k=2}^{\nt} U_{k-2}^2 V_{k-2} + \OO(n) \\
   &\quad
     + \text{linear combination of \ $\sum_{k=1}^{\nt} U_{k-2}$,
             \ $\sum_{k=1}^{\nt} V_{k-2}$, \ $\sum_{k=1}^{\nt} U_{k-2}^2$,
             \ $\sum_{k=1}^{\nt} V_{k-2}^2$ \ and
             \ $\sum_{k=1}^{\nt} U_{k-2} V_{k-2}$.}
  \end{aligned}
 \end{multline*}
Consequently
 \begin{align*}
  &\sum_{k=1}^{\nt} U_{k-1}^2 V_{k-1}
   =\frac{1}{1-\lambda}
    \sum_{k=2}^{\nt}
     \big[ U_{k-1}^2 V_{k-1} - \EE(U_{k-1}^2V_{k-1} \mid \cF_{k-2}) \big]
    - \frac{\lambda}{1-\lambda} U_{\nt-1}^2 V_{\nt-1} \\
  &+ \OO(n)
   + \text{linear combination of \ $\sum_{k=1}^{\nt} U_{k-2}$,
          \ $\sum_{k=1}^{\nt} V_{k-2}$, \ $\sum_{k=1}^{\nt} U_{k-2}^2$,
          \ $\sum_{k=1}^{\nt} V_{k-2}^2$ \ and
          \ $\sum_{k=1}^{\nt} U_{k-2}V_{k-2}$.}
 \end{align*}
Using \eqref{seged_UV_UNIFORM4} with \ $(\ell, i , j) = (8, 2, 1)$ \ we have
 \begin{align*}
    n^{-7/2}\sup_{t \in [0,T]}\,
           \Biggl\vert \sum_{k=2}^\nt
                    \big[U_{k-1}^2 V_{k-1}
                         - \EE(U_{k-1}^2 V_{k-1} \mid \cF_{k-2}) \big]
           \Biggr\vert
\stoch 0 \qquad \text{as \ $n\to\infty$.}
 \end{align*}
Thus,  in order to show \eqref{Zcond6a}, it suffices to prove
 \begin{gather}
  n^{-7/2} \sum_{k=1}^{\nT} U_k \stoch 0 , \qquad
  n^{-7/2} \sum_{k=1}^{\nT} U_k^2 \stoch 0 , \qquad
  n^{-7/2} \sum_{k=1}^{\nT} |V_k| \stoch 0 , \label{6supsumU_6supsumUU_6supsumV} \\
  n^{-7/2} \sum_{k=1}^{\nT} V_k^2 \stoch 0 , \qquad
  n^{-7/2} \sum_{k=1}^{\nT} \vert U_k V_k\vert \stoch 0 , \qquad
  n^{-7/2} \sup_{t \in [0,T]} | U_\nt^2 V_\nt | \stoch 0 
  \label{6supsumVV_6supsumUV_6supUUV}
 \end{gather}
 as \ $n \to \infty$.
\ Here \eqref{6supsumU_6supsumUU_6supsumV} and
 \eqref{6supsumVV_6supsumUV_6supUUV} follow by \eqref{seged_UV_UNIFORM1} and
 \eqref{seged_UV_UNIFORM2}, thus we conclude \eqref{Zcond6}.

Finally, we check condition (ii) of Theorem \ref{Conv2DiffThm}, i.e., the
 conditional Lindeberg condition
 \begin{equation}\label{Zcond3_new}
   \sum_{k=1}^{\lfloor nT \rfloor}
     \EE \big( \|\bZ^{(n)}_k\|^2 \bbone_{\{\|\bZ^{(n)}_k\| > \theta\}}
               \bmid \cF_{k-1} \big)
    \stoch 0 \qquad \text{for all \ $\theta>0$ \ and \ $T>0$.}
 \end{equation}
We have
 \ $\EE \big( \|\bZ^{(n)}_k\|^2 \bbone_{\{\|\bZ^{(n)}_k\| > \theta\}}
              \bmid \cF_{k-1} \big)
    \leq \theta^{-2} \EE \big( \|\bZ^{(n)}_k\|^4 \bmid \cF_{k-1} \big)$
 \ and
 \[
   \|\bZ^{(n)}_k\|^4
   \leq 3 \left( n^{-4} + n^{-8} U_{k-1}^4 + n^{-6} V_{k-1}^4 \right)
        \|\bM_{k-1}\|^4 .
 \]
Hence
 \[
    \sum_{k=1}^{\nT}
       \EE \big( \|\bZ^{(n)}_k\|^2 \bbone_{\{\|\bZ^{(n)}_k\| > \theta\}} \big)
         \to 0
    \qquad \text{as \ $n\to\infty$ \ for all \ $\theta>0$ \ and \ $T>0$,}
 \]
 since \ $\EE( \|\bM_k\|^4 ) = \OO(k^2)$,
 \ $\EE( \|\bM_k\|^4 U_{k-1}^4 ) \leq \sqrt{\EE(\|\bM_k\|^8) \EE(U_{k-1}^8)}
    = \OO(k^6)$
 \ and
 \ $\EE( \|\bM_k\|^4 V_{k-1}^4 ) \leq \sqrt{\EE(\|\bM_k\|^8) \EE(V_{k-1}^8)}
    = \OO(k^4)$
 \ by Corollary \ref{EEX_EEU_EEV}.
This yields \eqref{Zcond3_new}.
\proofend

We call the attention to the fact that our eighth order moment conditions
 \ $\EE(\|\bxi_{1,1,1}\|^8) < \infty$, \ $\EE(\|\bxi_{1,1,2}\|^8) < \infty$ \ and
 \ $\EE(\|\bvare_1\|^8) < \infty$ \ are used for applying Corollary
 \ref{EEX_EEU_EEV}.

\section{Proof of Theorem \ref{maint_Ad}}
\label{section_proof_maint}

The first convergence in Theorem \ref{maint_Ad} follows from Lemma
 \ref{limsupUVt}.

For the second convergence in Theorem \ref{maint_Ad}, consider the sequence of
 stochastic processes
 \[
   \bcZ^{(n)}_t
   := \begin{bmatrix}
       \bcM_t^{(n)} \\
       \bcN_t^{(n)} \\
       \bcP_t^{(n)}
      \end{bmatrix}
   := \sum_{k=1}^\nt
       \bZ^{(n)}_k
   \qquad \text{with} \qquad
   \bZ^{(n)}_k
   := \begin{bmatrix}
       n^{-1} \bM_k \\
       n^{-2} \bM_k U_{k-1} \\
       n^{-1} \bM_k V_{k-1}
      \end{bmatrix} 
 \]
 for \ $t \in \RR_+$ \ and \ $k, n \in \NN$.
\ Theorem \ref{maint_Ad} follows from Lemma \ref{main_VVt} and the following
 theorem (this will be explained after Theorem \ref{maint_conv}).

\begin{Thm}\label{maint_conv}
Suppose that the assumptions of Theorem \ref{maint_Ad} hold.
If \ $\langle \obV_{\!\!\bxi} \bv_\LEFT, \bv_\LEFT \rangle = 0$ \ then
 \begin{equation}\label{conv_Zt}
  \bcZ^{(n)} \distr \bcZ \qquad \text{as \ $n\to\infty$,}
 \end{equation}
 where the process \ $(\bcZ_t)_{t \in \RR_+}$ \ with values in
 \ $\RR^2 \times (\RR^2)^3$ \ is the unique strong solution of the SDE
 \begin{equation}\label{ZSDEt}
  \dd \bcZ_t
  = \gamma(t, \bcZ_t) \begin{bmatrix} \dd \bcW_t \\ \dd \tbcW_t \end{bmatrix} ,
  \qquad t \in \RR_+ ,
 \end{equation}
 with initial value \ $\bcZ_0 = \bzero$, \ where \ $(\bcW_t)_{t \in \RR_+}$ \ and
 \ $(\tbcW_t)_{t \in \RR_+}$ \ are independent 2-dimensional standard Wiener
 processes, and \ $\gamma : \RR_+ \times (\RR^2)^3 \to (\RR^{2\times2})^{3\times2}$
 \ is defined by
 \[
   \gamma(t, \bx)
   := \begin{bmatrix}
       (\langle \bu_\LEFT, \bx_1 + t \bm_\bvare \rangle^+)^{1/2} & \bzero \\
       (\langle \bu_\LEFT, \bx_1 + t \bm_\bvare \rangle^+)^{3/2} & \bzero \\
       \frac{\langle\bv_\LEFT,\bm_\bvare\rangle}{1-\lambda}
       (\langle \bu_\LEFT, \bx_1 + t \bm_\bvare \rangle^+)^{1/2}
        & \frac{\langle \bV_{\!\!\bvare} \bv_\LEFT, \bv_\LEFT \rangle^{1/2}}
               {(1-\lambda^2)^{1/2}}
          (\langle \bu_\LEFT, \bx_1 + t \bm_\bvare \rangle^+)^{1/2}
      \end{bmatrix}
      \otimes \obV_{\!\!\bxi}^{1/2}
 \]
 for \ $t \in \RR_+$ \ and \ $\bx = (\bx_1 , \bx_2 , x_3) \in (\RR^2)^3$.
\end{Thm}

As in the case of Theorem \ref{main_Ad}, the SDE \eqref{ZSDEt} has a unique
 strong solution with initial value \ $\bcZ_0 = \bzero$, \ for which we have
 \[
   \bcZ_t
   = \begin{bmatrix}
      \bcM_t \\
      \bcN_t \\
      \bcP_t
     \end{bmatrix}
   = \begin{bmatrix}
      \int_0^t \cY_s^{1/2} \, \obV_{\!\!\bxi}^{1/2} \, \dd \bcW_s \\[1mm]
      \int_0^t \cY_s \, \dd \bcM_s \\[1mm]
      \frac{\langle\bv_\LEFT,\bm_\bvare\rangle}{1-\lambda}
      \int_0^t \cY_s^{1/2} \, \obV_{\!\!\bxi}^{1/2} \, \dd \bcW_s
      + \frac{\langle \bV_{\!\!\bvare} \bv_\LEFT, \bv_\LEFT \rangle^{1/2}}
             {(1-\lambda^2)^{1/2}} \!
        \int_0^t \cY_s^{1/2} \, \obV_{\!\!\bxi}^{1/2} \, \dd \tbcW_s
     \end{bmatrix} , \qquad t \in \RR_+ .
 \]
One can again easily derive
 \begin{align}\label{convXZt}
  \begin{bmatrix} \bcX^{(n)} \\ \bcZ^{(n)} \end{bmatrix}
  \distr \begin{bmatrix} \tbcX \\ \bcZ \end{bmatrix} \qquad
  \text{as \ $n \to \infty$,}
 \end{align}
 where
 \[
   \bcX^{(n)}_t := n^{-1} \bX_\nt , \qquad
   \tbcX_t: = \langle \bu_\LEFT, \bcM_t + t \bm_\bvare \rangle \bu_\RIGHT ,
   \qquad t \in \RR_+ , \qquad n \in \NN .
 \]
Next, similarly to the proof of \eqref{seged2}, by Lemma \ref{Marci},
 convergence \eqref{convXZt} and Lemma \ref{main_VVt} imply
 \[
   \sum_{k=1}^n
    \begin{bmatrix}
     n^{-3} U_{k-1}^2 \\
     n^{-1} V_{k-1}^2 \\
     n^{-2} \bM_k U_{k-1} \\
     n^{-1/2} \bM_k V_{k-1}
    \end{bmatrix}
   \distr \begin{bmatrix}
           \int_0^1 \langle \bu_\LEFT, \tbcX_t \rangle^2 \, \dd t \\[1mm]
           \frac{\langle \bv_\LEFT, \bm_\bvare \rangle^2}{(1-\lambda)^2}
           + \frac{\langle \bV_{\!\!\bvare} \bv_\LEFT, \bv_\LEFT \rangle}
                  {1-\lambda^2} \\[1mm]
           \int_0^1 \cY_t \, \dd \bcM_t \\[1mm]
           \frac{\langle\bv_\LEFT,\bm_\bvare\rangle}{1-\lambda} \bcM_1
           + \frac{\langle \bV_{\!\!\bvare} \bv_\LEFT, \bv_\LEFT \rangle^{1/2}}
                  {(1-\lambda^2)^{1/2}}
             \int_0^1 \cY_t^{1/2} \, \obV_{\!\!\bxi}^{1/2} \, \dd \tbcW_t
          \end{bmatrix}
   \qquad \text{as \ $n \to \infty$.}
 \]
The limiting random vector can be written in the form as given in Theorem
 \ref{maint_Ad}, since \ $\langle \bu_\LEFT, \tbcX_t \rangle = \cY_t$ \ for all
 \ $t \in \RR_+$.

\noindent
\textbf{Proof of Theorem \ref{maint_conv}.}
Similar to the proof of Theorem \ref{main_conv}.
The conditional variance has the form
 \[
   \var\bigl(\bZ^{(n)}_k \mid \cF_{k-1}\bigr)
   = \begin{bmatrix}
      n^{-2} & n^{-3} U_{k-1} & n^{-2} V_{k-1} \\
      n^{-3} U_{k-1} & n^{-4} U_{k-1}^2 & n^{-3} U_{k-1} V_{k-1} \\
      n^{-2} V_{k-1} & n^{-3} U_{k-1} V_{k-1} & n^{-2} V_{k-1}^2
     \end{bmatrix}
     \otimes \bV_{\!\!\bM_k}
 \]
 for \ $n \in \NN$, \ $k \in \{1, \ldots, n\}$, \ with
 \ $\bV_{\bM_k} := \var(\bM_k \mid \cF_{k-1})$, \ and
 \ $\gamma(s, \bcZ_s^{(n)}) \gamma(s, \bcZ_s^{(n)})^\top$ \ has the form
 \[
   \begin{bmatrix}
    \langle \bu_\LEFT, \bcM_s^{(n)} + s \bm_\bvare \rangle
     & \langle \bu_\LEFT, \bcM_s^{(n)} + s \bm_\bvare \rangle^2
     & \frac{\langle\bv_\LEFT,\bm_\bvare\rangle}{1-\lambda}
       \langle \bu_\LEFT, \bcM_s^{(n)} + s \bm_\bvare \rangle \\
    \langle \bu_\LEFT, \bcM_s^{(n)} + s \bm_\bvare \rangle^2
     & \langle \bu_\LEFT, \bcM_s^{(n)} + s \bm_\bvare \rangle^3
     & \frac{\langle\bv_\LEFT,\bm_\bvare\rangle}{1-\lambda}
       \langle \bu_\LEFT, \bcM_s^{(n)} + s \bm_\bvare \rangle^2 \\
    \frac{\langle\bv_\LEFT,\bm_\bvare\rangle}{1-\lambda}
    \langle \bu_\LEFT, \bcM_s^{(n)} + s \bm_\bvare \rangle
     & \frac{\langle\bv_\LEFT,\bm_\bvare\rangle}{1-\lambda}
       \langle \bu_\LEFT, \bcM_s^{(n)} + s \bm_\bvare \rangle^2
     & M \langle \bu_\LEFT, \bcM_s^{(n)} + s \bm_\bvare \rangle
   \end{bmatrix}
   \otimes \obV_\xi
 \]
 for \ $s \in \RR_+$.
\ In order to check condition (i) of Theorem \ref{Conv2DiffThm}, we need to
 prove only that for each \ $T > 0$,
 \begin{gather}
  \sup_{t\in[0,T]}
   \bigg\| \frac{1}{n^2}
           \sum_{k=1}^{\nt} V_{k-1}^2 \bV_{\!\!\bM_k}
           - M \int_0^t
                \langle \bu_\LEFT, \bcM_s^{(n)} + s \bm_\bvare \rangle \obV_{\!\!\xi}
                \, \dd s \bigg\|
  \stoch 0 , \label{tZcond4} \\
  \sup_{t\in[0,T]}
   \bigg\| \frac{1}{n^2}
           \sum_{k=1}^{\nt} V_{k-1} \bV_{\!\!\bM_k} 
           - \frac{\langle\bv_\LEFT,\bm_\bvare\rangle}{1-\lambda}
             \int_0^t
              \langle \bu_\LEFT, \bcM_s^{(n)} + s \bm_\bvare \rangle \obV_{\!\!\xi}
              \, \dd s \bigg\|
  \stoch 0 , \label{tZcond5} \\
  \sup_{t\in[0,T]}
   \bigg\| \frac{1}{n^3}
           \sum_{k=1}^{\nt} U_{k-1} V_{k-1} \bV_{\!\!\bM_k}
           - \frac{\langle\bv_\LEFT,\bm_\bvare\rangle}{1-\lambda}
             \int_0^t
              \langle \bu_\LEFT, \bcM_s^{(n)} + s \bm_\bvare \rangle^2 \obV_{\!\!\xi}
              \, \dd s  \bigg\|
  \stoch 0  \label{tZcond6}
 \end{gather}
 as \ $n \to \infty$, \ since the rest, namely, \eqref{Zcond1}, \eqref{Zcond2} 
 and \eqref{Zcond3}, have already been proved.

We turn to prove \eqref{tZcond5}.
First we show that
 \begin{align}\label{tZcond5a}
  n^{-2}
  \sup_{t \in [0,T]}
   \Biggl\| \sum_{k=1}^{\nt} V_{k-1} \bV_{\!\!\bM_k}
            - \frac{\langle\bv_\LEFT,\bm_\bvare\rangle}{1-\lambda}
              \sum_{k=1}^{\nt} U_{k-1} \obV_{\!\!\bxi} \Biggr\|
  \stoch 0
 \end{align}
 as \ $n \to \infty$ \ for all \ $T > 0$.
\ We use the decomposition \eqref{VM2F}.
Using \eqref{seged_UV_UNIFORM1_mod} with \ $(\ell, i , j) = (4, 0, 2)$ \ and
 \ $(\ell, i , j) = (2, 0, 1)$, \ we have
 \begin{align*}
  n^{-2} \sum_{k=1}^{\nT} V_k^2 \stoch 0 , \qquad
  n^{-2} \sum_{k=1}^{\nT} |V_k| \stoch 0 \qquad \text{as \ $n\to\infty$,}
 \end{align*}
 hence \eqref{tZcond5a} will follow from
 \begin{align}\label{tZcond5b}
  n^{-2}
  \sup_{t \in [0,T]}
   \left\| \sum_{k=1}^{\nt} U_{k-1} V_{k-1}
           - \frac{\langle\bv_\LEFT,\bm_\bvare\rangle}{1-\lambda}
             \sum_{k=1}^{\nt} U_{k-1} \right\|
  \stoch 0 \qquad \text{as \ $n \to \infty$}
 \end{align}
 for all \ $T > 0$.
\ The aim of the following discussion is to decompose
 \ $\sum_{k=1}^{\nt} U_{k-1} V_{k-1}$ \ as a sum of a martingale and some
 other terms.
Using recursions \eqref{rec_U}, \eqref{rec_V_m} and formula \eqref{Mcond}, we
 obtain
 \begin{align*}
  \EE(U_{k-1} V_{k-1} \mid \cF_{k-2})
  &= \EE\Bigl((U_{k-2} + \langle \bu_\LEFT, \bM_{k-1} + \bm_\bvare \rangle)
              \bigl(\lambda V_{k-2}
                    + \langle \bv_\LEFT, \bvare_{k-1} \rangle\bigr)
              \,\Big|\, \cF_{k-2}\Bigr) \\
  &= \lambda U_{k-2} V_{k-2}
     + \langle \bv_\LEFT, \bm_\bvare \rangle U_{k-2}
     + \text{constant}
     + \text{constant$\times V_{k-2}$.}
 \end{align*}
Thus
 \begin{align*}
  \sum_{k=1}^{\nt} U_{k-1} V_{k-1}
  &= \sum_{k=2}^{\nt}
      \big[U_{k-1} V_{k-1} - \EE(U_{k-1} V_{k-1} \mid \cF_{k-2}) \big]
     + \sum_{k=2}^{\nt} \EE(U_{k-1} V_{k-1}^2 \mid \cF_{k-2}) \\
  &= \sum_{k=2}^{\nt}
      \big[U_{k-1} V_{k-1} - \EE(U_{k-1} V_{k-1} \mid \cF_{k-2}) \big]
     + \lambda \sum_{k=2}^{\nt} U_{k-2} V_{k-2} \\
  &\quad
     + \langle \bv_\LEFT, \bm_\bvare \rangle \sum_{k=2}^{\nt} U_{k-2}
     + \OO(n)
     + \text{constant$\times \sum_{k=2}^{\nt} V_{k-2}$.}
 \end{align*}
Consequently,
 \begin{align*}
  \sum_{k=1}^{\nt} U_{k-1} V_{k-1}
  &= \frac{1}{1-\lambda}
     \sum_{k=2}^{\nt}
      \big[U_{k-1} V_{k-1} - \EE(U_{k-1} V_{k-1} \mid \cF_{k-2}) \big]
     + \frac{\langle\bv_\LEFT,\bm_\bvare\rangle}
            {1-\lambda}
       \sum_{k=2}^{\nt} U_{k-2} \\
  &\quad
     - \frac{\lambda}{1-\lambda} U_{\nt - 1} V_{\nt - 1}
     + \OO(n)
     + \text{constant$\times \sum_{k=2}^{\nt} V_{k-2}$.}
 \end{align*}
Using \eqref{seged_UV_UNIFORM4_mod} with \ $(\ell, i , j) = (4, 1, 1)$ \ we
 have \begin{align*}
  n^{-2} \sup_{t \in [0,T]} \,
         \Biggl| \sum_{k=2}^\nt
                  \big[ U_{k-1} V_{k-1}
                        - \EE(U_{k-1} V_{k-1} \mid \cF_{k-2}) \big] \Biggr|
  \stoch 0 \qquad \text{as \ $n\to\infty$.}
 \end{align*}
Thus, in order to show \eqref{tZcond5b}, it suffices to prove
 \begin{gather}\label{5supsumV_mod_5supUV_mod_5supU_mod} 
  n^{-2} \sum_{k=1}^{\nT} |V_k| \stoch 0 , \qquad
  n^{-2} \sup_{t \in [0,T]} U_\nt V_\nt \stoch 0 , \qquad
  n^{-2} \sup_{t \in [0,T]} U_\nt \stoch 0
 \end{gather}
 as \ $n \to \infty$.
\ By \eqref{seged_UV_UNIFORM1_mod} with \ $(\ell, i , j) = (2, 0, 1)$, \ and by
 \eqref{seged_UV_UNIFORM2_mod} with \ $(\ell, i , j) = (2, 1, 1)$ \ and
 \ $(\ell, i , j) = (2, 1, 0)$, \ we have
 \eqref{5supsumV_mod_5supUV_mod_5supU_mod}.
Thus we conclude \eqref{tZcond5b}, and hence \eqref{tZcond5a}.
By Lemma \ref{Moments} and \eqref{seged_UV_UNIFORM1_mod} with
 \ $(\ell, i , j) = (2, 1, 1)$ \ and \ $(\ell, i , j) = (2, 1, 0)$, \ we get
 \begin{align}\label{tZcond2a_mod}
  n^{-2}
  \sup_{t \in [0,T]}
   \left\| \sum_{k=1}^{\nt} \bV_{\!\!\bM_k}
           - \sum_{k=1}^{\nt} U_{k-1} \obV_{\!\!\bxi} \right\|
  \stoch 0
 \end{align}
 as \ $n \to \infty$ \ for all \ $T > 0$.
\ As a last step, using \eqref{Zcond1} and \eqref{tZcond5a}, we obtain
 \eqref{tZcond5}.

Next we turn to prove \eqref{tZcond4}.
First we show that
 \begin{align}\label{tZcond4a}
  n^{-2}
  \sup_{t \in [0,T]}
   \Biggl\| \sum_{k=1}^{\nt} V_{k-1}^2 \bV_{\!\!\bM_k}
            - M \sum_{k=1}^{\nt} U_{k-1} \obV_{\!\!\bxi} \Biggr\|
  \stoch 0
 \end{align}
 as \ $n \to \infty$ \ for all \ $T > 0$.
\ We use the decomposition \eqref{V2M2F}.
Using \eqref{seged_UV_UNIFORM1_mod} with \ $(\ell, i , j) = (6, 0, 3)$ \ and
 \ $(\ell, i , j) = (4, 0, 2)$, \ we have
 \begin{align*}
  n^{-2} \sum_{k=1}^{\nT} | V_k |^3 \stoch 0 , \qquad
  n^{-2} \sum_{k=1}^{\nT} V_k^2 \stoch 0 \qquad \text{as \ $n\to\infty$,}
 \end{align*}
 hence \eqref{tZcond4a} will follow from
 \begin{align}\label{tZcond4b}
  n^{-2}
  \sup_{t \in [0,T]}
   \left\| \sum_{k=1}^{\nt} U_{k-1} V_{k-1}^2 - M \sum_{k=1}^{\nt} U_{k-1} \right\|
  \stoch 0 \qquad \text{as \ $n \to \infty$}
 \end{align}
 for all \ $T > 0$.
\ The aim of the following discussion is to decompose
 \ $\sum_{k=1}^{\nt} U_{k-1} V_{k-1}^2$ \ as a sum of a martingale and some
 other terms.
Using recursions \eqref{rec_U}, \eqref{rec_V_m} and formula \eqref{Mcond}, we
 obtain
 \begin{align*}
  \EE(U_{k-1} V_{k-1}^2 \mid \cF_{k-2})
  &= \EE\Bigl((U_{k-2} + \langle \bu_\LEFT, \bM_{k-1} + \bm_\bvare \rangle)
              \bigl(\lambda V_{k-2}
                    + \langle \bv_\LEFT, \bvare_{k-1} \rangle\bigr)^2
              \,\Big|\, \cF_{k-2}\Bigr) \\
  &= \lambda^2 U_{k-2} V_{k-2}^2
     + 2 \lambda \langle \bv_\LEFT, \bm_\bvare \rangle U_{k-2} V_{k-2}
     + \EE(\langle \bv_\LEFT, \bvare_{k-1} \rangle^2) \, U_{k-2} \\
  &\quad 
     + \text{constant}
     + \text{constant $\times$ $V_{k-2}$.}
 \end{align*}
Thus
 \begin{align*}
  &\sum_{k=1}^{\nt} U_{k-1} V_{k-1}^2
   = \sum_{k=2}^{\nt}
      \big[U_{k-1} V_{k-1}^2 - \EE(U_{k-1} V_{k-1}^2 \mid \cF_{k-2}) \big]
     + \sum_{k=2}^{\nt} \EE(U_{k-1} V_{k-1}^2 \mid \cF_{k-2}) \\
  &= \sum_{k=2}^{\nt}
      \big[U_{k-1} V_{k-1}^2 - \EE(U_{k-1} V_{k-1}^2 \mid \cF_{k-2}) \big]
     + \lambda^2 \sum_{k=2}^{\nt} U_{k-2} V_{k-2}^2
     + 2 \lambda \langle \bv_\LEFT, \bm_\bvare \rangle
       \sum_{k=2}^{\nt} U_{k-2} V_{k-2} \\
  &\quad
     + \EE(\langle \bv_\LEFT, \bvare_{k-1} \rangle^2) \sum_{k=2}^{\nt} U_{k-2}
     + \OO(n)
     + \text{constant $\times$ $\sum_{k=2}^{\nt} V_{k-2}$.}
 \end{align*}
Consequently,
 \begin{align*}
  &\sum_{k=1}^{\nt} U_{k-1} V_{k-1}^2
   = \frac{1}{1-\lambda^2}
     \sum_{k=2}^{\nt}
      \big[U_{k-1} V_{k-1}^2 - \EE(U_{k-1} V_{k-1}^2 \mid \cF_{k-2}) \big]
     + \frac{2\lambda\langle\bv_\LEFT,\bm_\bvare\rangle}
            {1-\lambda^2}
       \sum_{k=2}^{\nt} U_{k-2} V_{k-2} \\
  &\qquad\quad
     + \frac{\EE(\langle \bv_\LEFT, \bvare_{k-1} \rangle^2)}
            {1-\lambda^2}
       \sum_{k=2}^{\nt} U_{k-2} 
     - \frac{\lambda^2}{1-\lambda^2} U_{\nt - 1} V_{\nt - 1}^2
     + \OO(n)
     + \text{constant $\times$ $\sum_{k=2}^{\nt} V_{k-2}$.}
 \end{align*}
Using \eqref{seged_UV_UNIFORM4_mod} with \ $(\ell, i , j) = (8, 1, 2)$ \ we
 have \begin{align*}
  n^{-2} \sup_{t \in [0,T]} \,
         \Biggl| \sum_{k=2}^\nt
                  \big[ U_{k-1} V_{k-1}^2
                        - \EE(U_{k-1} V_{k-1}^2 \mid \cF_{k-2}) \big] \Biggr|
  \stoch 0 \qquad \text{as \ $n \to \infty$.}
 \end{align*}
By \eqref{seged_UV_UNIFORM1_mod} with \ $(\ell, i , j) = (2, 0, 1)$, \ and by
 \eqref{seged_UV_UNIFORM2_mod} with \ $(\ell, i , j) = (4, 1, 2)$, \ we obtain
 \begin{gather}
  n^{-2} \sum_{k=1}^{\nT} |V_k| \stoch 0 , \qquad
  n^{-2} \sup_{t \in [0,T]} U_\nt V_\nt^2 \stoch 0 \label{4supsumV_mod_4supUVV_mod}
   \end{gather}
 as \ $n \to \infty$, \ hence
 \begin{align*}
  n^{-2}
  \sup_{t \in [0,T]}
   \left\| \sum_{k=1}^{\nt} U_{k-1} V_{k-1}^2
           - \frac{2\lambda\langle\bv_\LEFT,\bm_\bvare\rangle}
                  {1-\lambda^2}
             \sum_{k=2}^{\nt} U_{k-2} V_{k-2}
           - \frac{\EE(\langle \bv_\LEFT, \bvare_{k-1} \rangle^2)}
                  {1-\lambda^2}
             \sum_{k=2}^{\nt} U_{k-2} \right\|
  \stoch 0
 \end{align*}
 as \ $n \to \infty$ \ for all \ $T > 0$.
\ Thus, taking into account \eqref{tZcond5b}, we conclude \eqref{tZcond4b},
 and hence \eqref{tZcond4a}, since
 \[
   \frac{\langle \bv_\LEFT, \bm_\bvare \rangle}{1-\lambda}
   \frac{2\lambda\langle\bv_\LEFT,\bm_\bvare\rangle}
        {1-\lambda^2}
   + \frac{\EE(\langle \bv_\LEFT, \bvare_{k-1} \rangle^2)}
          {1-\lambda^2}
   = \frac{\langle \bv_\LEFT, \bm_\bvare \rangle^2}{(1-\lambda)^2}
     + \frac{\langle \bV_{\!\!\bvare} \bv_\LEFT, \bv_\LEFT \rangle}
            {1-\lambda^2} .
 \]
As a last step, using \eqref{tZcond2a_mod} and \eqref{Zcond1}, we obtain
 \eqref{tZcond4}.

Finally we turn to prove \eqref{tZcond6}.
First we show that
 \begin{align}\label{tZcond6a}
  n^{-3}
  \sup_{t \in [0,T]}
   \Biggl\| \sum_{k=1}^{\nt} U_{k-1} V_{k-1} \bV_{\!\!\bM_k}
            - \frac{\langle\bv_\LEFT,\bm_\bvare\rangle}{1-\lambda}
              \sum_{k=1}^{\nt} U_{k-1}^2 \obV_{\!\!\bxi} \Biggr\|
  \stoch 0
 \end{align}
 as \ $n \to \infty$ \ for all \ $T > 0$.
\ We use the decomposition \eqref{UVM2F}.
Using \eqref{seged_UV_UNIFORM1_mod} with \ $(\ell, i , j) = (4, 1, 2)$ \ and
 \ $(\ell, i , j) = (2, 1, 1)$, \ we have
 \begin{align*}
  n^{-3} \sum_{k=1}^{\nT} U_k V_k^2 \stoch 0 , \qquad
  n^{-3} \sum_{k=1}^{\nT} U_k |V_k| \stoch 0 \qquad \text{as \ $n\to\infty$,}
 \end{align*}
 hence \eqref{tZcond6a} will follow from
 \begin{align}\label{tZcond6b}
  n^{-3}
  \sup_{t \in [0,T]}
   \left\| \sum_{k=1}^{\nt} U_{k-1}^2 V_{k-1}
           - \frac{\langle\bv_\LEFT,\bm_\bvare\rangle}{1-\lambda}
             \sum_{k=1}^{\nt} U_{k-1}^2 \right\|
  \stoch 0 \qquad \text{as \ $n \to \infty$}
 \end{align}
 for all \ $T > 0$.
\ The aim of the following discussion is to decompose
 \ $\sum_{k=1}^{\nt} U_{k-1}^2 V_{k-1}$ \ as a sum of a martingale and some
 other terms.
Using recursions \eqref{rec_U}, \eqref{rec_V_m} and formula \eqref{Mcond}, we
 obtain
 \begin{align*}
  \EE(U_{k-1}^2 V_{k-1} \mid \cF_{k-2})
  &= \EE\Bigl((U_{k-2} + \langle \bu_\LEFT, \bM_{k-1} + \bm_\bvare \rangle)^2
              \bigl(\lambda V_{k-2}
                    + \langle \bv_\LEFT, \bvare_{k-1} \rangle\bigr)
              \,\Big|\, \cF_{k-2}\Bigr) \\
  &= \lambda U_{k-2}^2 V_{k-2}
     + \langle \bv_\LEFT, \bm_\bvare \rangle U_{k-2}^2
     + \text{constant} \\
  &\quad 
     + \text{linear combinations of \ $U_{k-2} V_{k-2}$, \ $U_{k-2}$ \ and
             \ $V_{k-2}$.}
 \end{align*}
Thus
 \begin{align*}
  &\sum_{k=1}^{\nt} U_{k-1}^2 V_{k-1}
   = \sum_{k=2}^{\nt}
      \big[U_{k-1}^2 V_{k-1} - \EE(U_{k-1}^2 V_{k-1} \mid \cF_{k-2}) \big]
     + \sum_{k=2}^{\nt} \EE(U_{k-1}^2 V_{k-1} \mid \cF_{k-2}) \\
  &= \sum_{k=2}^{\nt}
      \big[U_{k-1}^2 V_{k-1} - \EE(U_{k-1}^2 V_{k-1} \mid \cF_{k-2}) \big]
     + \lambda \sum_{k=2}^{\nt} U_{k-2}^2 V_{k-2}
     + \langle \bv_\LEFT, \bm_\bvare \rangle \sum_{k=2}^{\nt} U_{k-2}^2 \\
  &\quad
     + \OO(n)
     + \text{linear combinations of \ $\sum_{k=2}^{\nt} U_{k-2} V_{k-2}$,
             \ $\sum_{k=2}^{\nt} U_{k-2}$ \ and \ $\sum_{k=2}^{\nt} V_{k-2}$.}
 \end{align*}
Consequently,
 \begin{align*}
  \sum_{k=1}^{\nt} U_{k-1}^2 V_{k-1}
  &= \frac{1}{1-\lambda}
     \sum_{k=2}^{\nt}
      \big[U_{k-1}^2 V_{k-1} - \EE(U_{k-1}^2 V_{k-1} \mid \cF_{k-2}) \big] \\
  &\quad
      + \frac{\langle\bv_\LEFT,\bm_\bvare\rangle}
            {1-\lambda}
       \sum_{k=2}^{\nt} U_{k-2}^2
     - \frac{\lambda}{1-\lambda} U_{\nt-1}^2 V_{\nt-1} \\
  &\quad
     + \OO(n)
     + \text{linear combinations of \ $\sum_{k=2}^{\nt} U_{k-2} V_{k-2}$,
             \ $\sum_{k=2}^{\nt} U_{k-2}$ \ and \ $\sum_{k=2}^{\nt} V_{k-2}$.}
 \end{align*}
Using \eqref{seged_UV_UNIFORM4_mod} with \ $(\ell, i , j) = (8, 2, 1)$ \ we
 have \begin{align*}
  n^{-3} \sup_{t \in [0,T]} \,
         \Biggl| \sum_{k=2}^\nt
                  \big[ U_{k-1}^2 V_{k-1}
                        - \EE(U_{k-1}^2 V_{k-1} \mid \cF_{k-2}) \big] \Biggr|
  \stoch 0 \qquad \text{as \ $n \to \infty$.}
 \end{align*}
Thus, in order to show \eqref{tZcond6b}, it suffices to prove
 \begin{gather}\label{6supsumUV_mod_6supsumU_mod}
  n^{-3} \sum_{k=1}^{\nT} U_k |V_k| \stoch 0 , \qquad
  n^{-3} \sum_{k=1}^{\nT} U_k \stoch 0, \\
  n^{-3} \sum_{k=1}^{\nT} |V_k| \stoch 0, \qquad
  n^{-3} \sup_{t \in [0,T]} U_\nt^2 |V_\nt| \stoch 0
  \label{6supsumV_mod_6supUUV_mod}
 \end{gather}
 as \ $n \to \infty$. 
\ By \eqref{seged_UV_UNIFORM1_mod} with \ $(\ell, i , j) = (4, 1, 1)$ \ and
  \ $(\ell, i , j) = (4, 1, 0)$ \ we obtain
 \eqref{6supsumUV_mod_6supsumU_mod}.
By \eqref{seged_UV_UNIFORM1_mod} with \ $(\ell, i , j) = (4, 0, 1)$ \ and
 \eqref{seged_UV_UNIFORM2_mod} with \ $(\ell, i , j) = (4, 2, 1)$, \ we obtain
 \eqref{6supsumV_mod_6supUUV_mod}.
Thus we conclude \eqref{tZcond6b}, and hence \eqref{tZcond6a}.
By Lemma \ref{Moments} and \eqref{seged_UV_UNIFORM1_mod} with
 \ $(\ell, i , j) = (4, 2, 1)$ \ and \ $(\ell, i , j) = (4, 2, 0)$, \ we get
 \begin{align}\label{tZcond6a_mod}
  n^{-3}
  \sup_{t \in [0,T]}
   \left\| \sum_{k=1}^{\nt} U_{k-1} \bV_{\!\!\bM_k}
           - \sum_{k=1}^{\nt} U_{k-1}^2 \obV_{\!\!\bxi} \right\|
  \stoch 0
 \end{align}
 as \ $n \to \infty$ \ for all \ $T > 0$.
\ As a last step, using \eqref{tZcond6a} and \eqref{Zcond2}, we obtain
 \eqref{tZcond6}.

Condition (ii) of Theorem \ref{Conv2DiffThm} can be checked again as in case of
 Theorem \ref{main_conv}.
\proofend

\section{Proof of Theorem \ref{main_sub}}
\label{section_proof_sub}
 
By Quine and Durham \cite{QD}, the Markov chain \ $(\bX_k)_{k\in\ZZ_+}$ \ admits a
 unique stationary distribution, and we have
 \begin{equation}\label{SLLN}
  \frac{1}{n} \sum_{k=1}^n f(\bX_{k-1}, \bX_k)
  \as
  \EE\biggl(f\biggl(\tbX,
                    \sum_{j=1}^{\tX_1} \bxi_{1,j,1}
                    + \sum_{j=1}^{\tX_2} \bxi_{1,j,2}
                    + \bvare_1\biggr)\biggr) ,
  \qquad \text{as \ $n \to \infty$,}
 \end{equation}
 for all Borel measurable functions \ $f : \RR^2 \times \RR^2 \to \RR$ \ with
 \[
   \EE\biggl(\biggl|f\biggl(\tbX,
                            \sum_{j=1}^{\tX_1} \bxi_{1,j,1}
                            + \sum_{j=1}^{\tX_2} \bxi_{1,j,2}
                            + \bvare_1\biggr)\biggr|\biggr)
   < \infty ,
 \]
 see (2.1) in Quine and Durham \cite{QD}.
By Quine \cite{Q}, \ we have \ $\EE(\|\tbX\|^2) < \infty$, \ hence
 \begin{equation}\label{SLLN1}
  \frac{1}{n} \sum_{k=1}^n \bX_{k-1} \bX_{k-1}^\top
  \as
  \EE\bigl(\tbX \tbX^\top\bigr) ,
  \qquad \text{as \ $n \to \infty$.}
 \end{equation}
The aim of the following discussion is to show that the matrix
 \ $\EE\bigl(\tbX \tbX^\top\bigr)$ \ is invertible.
By Quine \cite{Q}, \ we have
 \[
   \EE\bigl(\tbX \tbX^\top\bigr)
   = \sum_{i=0}^\infty
      \bm_\bxi^i
      \bigl(\EE(\tX_1) \bV_{\!\!\bxi_1} + \EE(\tX_2) \bV_{\!\!\bxi_2}
            + \bV_{\!\!\bvare}\bigr)
      (\bm_\bxi^\top)^i
     + \EE(\tbX) \EE(\tbX^\top) .
 \]
If \ $\EE\bigl(\tbX \tbX^\top\bigr)$ \ is not invertible, then there exists
 \ $\bw \in \RR^2 \setminus \{\bzero\}$ \ with
 \ $\bw^\top \EE\bigl(\tbX \tbX^\top\bigr) \bw = 0$.
\ But then \ $\EE(\tX_1) \bw^\top \bV_{\!\!\bxi_1} \bw = 0$,
 \ $\EE(\tX_2) \bw^\top \bV_{\!\!\bxi_2} \bw = 0$ \ and
 \ $\bw^\top \bV_{\!\!\bvare} \bw = 0$.
\ We have \ $\EE(\tX_1) > 0$ \ and \ $\EE(\tX_2) > 0$, \ since, by Quine \cite{Q}, \ we
 have \ $\EE(\tbX) = \sum_{i=0}^\infty \bm_\bxi^i \bm_\bvare$, \ where
 \ $\bm_\bvare \ne \bzero$ \ and there exists \ $i \in \NN$ \ such that the entries of
 \ $\bm_\bxi^i$ \ are positive.
Consequently, we obtain \ $\bw^\top \bV_{\!\!\bxi_1} \bw = 0$,
 \ $\bw^\top \bV_{\!\!\bxi_2} \bw = 0$ \ and \ $\bw^\top \bV_{\!\!\bvare} \bw = 0$,
 \ which is impossible, since at least one of the matrices \ $\bV_{\!\!\bxi_1}$,
 \ $\bV_{\!\!\bxi_2}$, \ $\bV_{\!\!\bvare}$ \ is invertible, hence
 \ $\EE\bigl(\tbX \tbX^\top\bigr)$ \ has to be invertible, and we conclude
 \begin{equation}\label{SLLN2}
  \biggl(\frac{1}{n} \sum_{k=1}^n \bX_{k-1} \bX_{k-1}^\top\biggl)^{-1}
  \as
  \bigl[\EE\bigl(\tbX \tbX^\top\bigr)\bigr]^{-1} ,
  \qquad \text{as \ $n \to \infty$.}
 \end{equation}
Applying again \eqref{SLLN} and using the decomposition \eqref{Mdeco},
 \begin{align*}
  &\frac{1}{n} \sum_{k=1}^n (\bX_k - \bm_\bxi \bX_{k-1} - \bm_\bxi) \bX_{k-1}^\top \\
  &\as \EE\Biggl(\biggl(\sum_{j=1}^{\tX_1} (\bxi_{1,j,1} - \EE(\bxi_{1,j,1}))
                        + \sum_{j=1}^{\tX_2} (\bxi_{1,j,2} - \EE(\bxi_{1,j,2}))
                        + (\bvare_1 - \EE(\bvare_1))\biggr) \tbX^\top\Biggr) \\
  &\phantom{\as}
   = \EE\Biggl[\EE\biggl(\sum_{j=1}^{\tX_1} (\bxi_{1,j,1} - \EE(\bxi_{1,j,1}))
                         + \sum_{j=1}^{\tX_2} (\bxi_{1,j,2} - \EE(\bxi_{1,j,2}))
                         + (\bvare_1 - \EE(\bvare_1))
                         \, \bigg| \, \tbX\biggr) \tbX^\top \Biggr]
   = \bzero
 \end{align*}
 as \ $n \to \infty$.
\ Consequently,
 \[
   \hbmbxi^{(n)} - \bm_\bxi
   = \biggl(\frac{1}{n}
            \sum_{k=1}^n (\bX_k - \bm_\bxi \bX_{k-1} - \bm_\bxi) \bX_{k-1}^\top\biggr)
     \biggl(\frac{1}{n} \sum_{k=1}^n \bX_{k-1} \bX_{k-1}^\top\biggr)^{-1}
   \as \bzero
 \]
 as \ $n \to \infty$, \ hence we obtain the strong consistency of \ $\hbmbxi^{(n)}$.
\ By the continuity of the function \ $r$, \ this implies the strong consistency of
 \ $\hvarrho_n = r(\hbmbxi^{(n)})$.

The asymptotic normality \eqref{main_sub} can be proved by the martingale central
 limit theorem. 
We can write
 \[
   n^{1/2} (\hbmbxi^{(n)} - \bm_\bxi)
   = \biggl(\frac{1}{n^{1/2}} \sum_{k=1}^n \bZ_k \biggr)
     \biggl(\frac{1}{n} \sum_{k=1}^n \bX_{k-1} \bX_{k-1}^\top\biggr)^{-1}
 \]
 with \ $\bZ_k := \bM_k \bX_{k-1}^\top$. 
\ We have
 \ $\EE\bigl(\bZ_k^{\otimes2} \,\big|\, \cF_{k-1}\bigr)
    = \EE\bigl(\bM_k^{\otimes2} \,\big|\, \cF_{k-1}\bigr)
      (\bX_{k-1}^\top)^{\otimes2}$.
\ By the decomposition \eqref{Mdeco},
 \begin{equation}\label{M2cond}
  \EE\bigl(\bM_k^{\otimes2} \,\big|\, \cF_{k-1}\bigr)
   = \sum_{i=1}^2 
      X_{k-1,i} \EE\bigl[(\bxi_{1,1,i} - \EE(\bxi_{1,1,i}))^{\otimes2}\bigr]
     + \EE\bigl[(\bvare_1 - \EE(\bvare_1))^{\otimes2}\bigr] , \qquad k \in \NN ,
 \end{equation}
 thus by \eqref{SLLN}, the asymptotic covariances have the form
 \begin{align*}
  \frac{1}{n} \sum_{k=1}^n \EE\bigl(\bZ_k^{\otimes2} \,\big|\, \cF_{k-1}\bigr)
  &\as \sum_{i=1}^2 
        \EE\bigl[(\bxi_{1,1,i} - \EE(\bxi_{1,1,i}))^{\otimes2}\bigr]
        \EE\biggl[\tX_i \Bigl(\tbX^\top\Bigr)^{\otimes2}\biggr] \\
  &\phantom{\as}
       + \EE\bigl[(\bvare_1 - \EE(\bvare_1))^{\otimes2}\bigr]
         \EE\biggl[\Bigl(\tbX^\top\Bigr)^{\otimes2}\biggr] \qquad
  \text{as \ $n \to \infty$.}
 \end{align*}
The aim of the following discussion is to check the conditional Lindeberg condition
 \begin{equation}\label{CLC}
  \frac{1}{n}
  \sum_{k=1}^n
   \EE(\|\bZ_k\|^2 \bbone_{\{\|\bZ_k\|>\theta\sqrt{n}\}} \mid \cF_{k-1})
  \stoch 0 \qquad \text{as \ $n \to \infty$ \ for all \ $\theta > 0$.}
 \end{equation}
By the decomposition \eqref{Mdeco},
 \[
   \EE\bigl(\bM_k^{\otimes4} \,\big|\, \cF_{k-1}\bigr) = \bP(X_{k-1,1}, X_{k-1,2})
 \]
 with \ $\bP = (P_1, \ldots, P_{16}) : \RR^2 \to \RR^{16}$, \ where
 \ $P_1$, \ldots, $P_{16}$ \ are polynomials having degree at most 2, and their
 coefficients depend on the moments
 \ $\EE[(\bxi_{1,1,1} - \EE(\bxi_{1,1,1})^{\otimes4}]$,
 \ $\EE[(\bxi_{1,1,2} - \EE(\bxi_{1,1,2})^{\otimes4}]$,
 \ $\EE[(\bvare_1 - \EE(\bvare_1)^{\otimes4}]$, \ $\bV_{\!\!\bxi_1}$,
 \ $\bV_{\!\!\bxi_2}$ \ and $\bV_{\!\!\bvare}$.
\ Thus
 \begin{align*}
  &\frac{1}{n}
   \sum_{k=1}^n
    \EE(\|\bZ_k\|^2 \bbone_{\{\|\bZ_k\|>\theta\sqrt{n}\}} \mid \cF_{k-1})
   \leq
   \frac{1}{n^2\theta^2} \sum_{k=1}^n \EE(\|\bZ_k\|^4 \mid \cF_{k-1}) \\
  &\leq
   \frac{1}{n^2\theta^2} \sum_{k=1}^n \|\bX_k\|^4 \EE(\|\bM_k\|^4 \mid \cF_{k-1})
   = \frac{1}{n^2\theta^2} \sum_{k=1}^n \|\bX_k\|^4 \bP(X_{k-1,1}, X_{k-1,2})
   \stoch 0
 \end{align*}
 as \ $n \to \infty$ \ for all \ $\theta > 0$, \ since, by Lemma
 \ref{LEM_moments_X_sub}, \ $\EE(\|\bX_k\|^4 |\bP(X_{k-1,1}, X_{k-1,2})|) = \OO(1)$,
 \ hence \ $n^{-2} \sum_{k=1}^n \EE(\|\bX_k\|^4 |\bP(X_{k-1,1}, X_{k-1,2})|) \to 0$ \ as
 \ $n \to \infty$.
\ Consequently, by the martingale central limit theorem, we obtain
 \[
   \frac{1}{n^{1/2}} \sum_{k=1}^n \bZ_k \distr \tbZ \qquad \text{as \ $n \to \infty$,}
 \]
 where \ $\tbZ$ \ is a \ $2 \times 2$ \ random matrix having a normal distribution
 with zero mean and with
 \begin{align*}
  \EE\Bigl(\tbZ^{\otimes2}\Bigr)
  = \sum_{i=1}^2 
      \EE\bigl[(\bxi_{1,1,i} - \EE(\bxi_{1,1,i}))^{\otimes2}\bigr]
      \EE\biggl[\tX_i \Bigl(\tbX^\top\Bigr)^{\otimes2}\biggr] 
     + \EE\bigl[(\bvare_1 - \EE(\bvare_1))^{\otimes2}\bigr]
       \EE\biggl[\Bigl(\tbX^\top\Bigr)^{\otimes2}\biggr] .
 \end{align*}
Using \eqref{SLLN2} and applying Slutsky's lemma, we obtain \eqref{main_sub}.

The convergence \eqref{varrho_n_sub} follows from \eqref{main_sub} by the so called
 Delta Method with the function \ $r$, \ see, e.g., Lehmann and Romano
 \cite[Theorem 11.2.14]{LR}.

\vspace*{5mm}


\appendix


\vspace*{5mm}

\noindent{\bf\Large Appendices}

\section{Estimations of moments}
\label{section_moments}

In the proof of Theorem \ref{main}, good bounds for moments of the random
 vectors and variables \ $(\bM_k)_{k\in\ZZ_+}$, \ $(\bX_k)_{k\in\ZZ_+}$,
 \ $(U_k)_{k\in\ZZ_+}$ \ and \ $(V_k)_{k\in\ZZ_+}$ \ are extensively used.
First note that, for all \ $k \in \NN$, \ $\EE( \bM_k \mid \cF_{k-1} ) = \bzero$
 \ and \ $\EE(\bM_k) = \bzero$, \ since
 \ $\bM_k = \bX_k - \EE(\bX_k \mid \cF_{k-1})$.

\begin{Lem}\label{Moments}
Let \ $(\bX_k)_{k\in\ZZ_+}$ \ be a 2-type Galton--Watson process with immigration
 and with \ $\bX_0 = \bzero$.
\ If \ $\EE(\|\bxi_{1,1,1}\|^2) < \infty$, \ $\EE(\|\bxi_{1,1,2}\|^2) < \infty$
 \ and \ $\EE(\|\bvare_1\|^2) < \infty$ \ then
 \begin{equation}\label{Mcond}
  \var( \bM_k \mid \cF_{k-1} )
  = X_{k-1,1} \bV_{\!\!\bxi_1} + X_{k-1,2} \bV_{\!\!\bxi_2} + \bV_{\!\!\bvare}
  = U_{k-1} \obV_{\!\!\bxi} + V_{k-1} \tbV_{\!\!\bxi} + \bV_{\!\!\bvare}
 \end{equation}
 for all \ $k \in \NN$, \ where
 \[
   \tbV_{\!\!\bxi}
   := \sum_{i=1}^2 \langle \be_i, \bv_\RIGHT \rangle \bV_{\!\!\bxi_i}
   = \frac{\beta\bV_{\!\!\bxi_1}-(1-\delta)\bV_{\!\!\bxi_2}}{\beta+1-\delta} .
 \]
If \ $\EE(\|\bxi_{1,1,1}\|^3) < \infty$, \ $\EE(\|\bxi_{1,1,2}\|^3) < \infty$
 \ and \ $\EE(\|\bvare_1\|^3) < \infty$ \ then, for all \ $k \in \NN$,
 \begin{align}\label{M3cond}
  \begin{split}
   \EE( \bM_k^{\otimes3} \mid \cF_{k-1} )
   &= X_{k-1,1} \EE[(\bxi_{1,1,1} - \EE(\bxi_{1,1,1})^{\otimes3}] \\
   &\quad
      + X_{k-1,2} \EE[(\bxi_{1,1,2} - \EE(\bxi_{1,1,2})^{\otimes3}]
      + \EE[(\bvare_1 - \EE(\bvare_1)^{\otimes3}] .
  \end{split}
 \end{align}
\end{Lem}

\noindent
\textbf{Proof.}
Using the decomposition \eqref{Mdeco}, where, for all \ $k \in \NN$, \ the
 random vectors
 \ $\big\{\bxi_{k,j,1} - \EE(\bxi_{k,j,1}), \, \bxi_{k,j.2} - \EE(\bxi_{k,j,2}), \,
          \bvare_k - \EE(\bvare_k)
          : j \in \NN \big\}$
 \ are independent of each other, independent of \ $\cF_{k-1}$, \ and have
 zero mean vector, we conclude \eqref{Mcond} and \eqref{M3cond}.
\proofend

\begin{Lem}\label{LEM_moments_seged0}
Let \ $(\bzeta_k)_{k\in\NN}$ \ be independent and identically distributed
 random vectors with values in \ $\RR^d$ \ such that
 \ $\EE(\|\bzeta_1\|^\ell) < \infty$ \ with some \ $\ell \in \NN$.
\begin{enumerate}
 \item[\textup{(i)}]
  Then there exists \ $\bQ = (Q_1, \ldots, Q_{d^\ell}) : \RR \to \RR^{d^\ell}$,
   \ where \ $Q_1$, \ldots, $Q_{d^\ell}$ \ are polynomials having degree at most
   \ $\ell-1$ \ such that
   \[
     \EE\bigl((\bzeta_1 + \cdots + \bzeta_N)^{\otimes \ell}\bigr)
     = N^\ell \bigl[\EE(\bzeta_1)\bigr]^{\otimes \ell} + \bQ(N) ,
     \qquad N \in \NN , \quad N \geq \ell .
   \]
 \item[\textup{(ii)}]
  If \ $\EE(\bzeta_1) = \bzero$ \ then there exists
   \ $\bR = (R_1, \ldots, R_{d^\ell}) : \RR \to \RR^{d^\ell}$, \ where
   \ $R_1$, \ldots, $R_{d^\ell}$ \ are polynomials having degree at most
   \ $\lfloor \ell/2 \rfloor$ \ such that
   \[
     \EE\bigl((\bzeta_1 + \cdots + \bzeta_N)^{\otimes \ell}\bigr) = \bR(N) ,
     \qquad N \in \NN , \quad N \geq \ell .
   \]
\end{enumerate}
The coefficients of the polynomials \ $\bQ$ \ and \ $\bR$ \ depend on the
 moments \ $\EE(\bzeta_{i_1} \otimes \cdots \otimes \bzeta_{i_\ell})$,
 \ $i_1, \ldots, i_\ell \in \{1, \ldots, N\}$.
\end{Lem}

\noindent
\textbf{Proof.}
(i) We have
 \begin{align*}
  \EE\bigl((\bzeta_1 + \cdots + \bzeta_N)^{\otimes\ell}\bigr)
  = \sum_{\underset{\SC k_1 + 2 k_2 + \cdots + s k_s = \ell, \; k_s \ne 0}
                   {s \in \{1,\ldots,\ell\}, \; k_1, \ldots, k_s \in \ZZ_+,}}
    &\binom{N}{k_1} \binom{N - k_1}{k_2} \cdots
     \binom{N - k_1 - \cdots - k_{s-1}}{k_s} \\
  &\times \sum_{(i_1,\ldots,i_\ell) \in P_{k_1,\ldots,k_s}^{(N,\ell)}}
           \EE(\bzeta_{i_1} \otimes \cdots \otimes \bzeta_{i_\ell}) ,
 \end{align*}
 where the set \ $P_{k_1,\ldots,k_s}^{(N,\ell)}$ \ consists of permutations of all
 the multisets containing pairwise different elements
 \ $j_{k_1}, \ldots, j_{k_s}$ \ of the set \ $\{1, \ldots, N\}$ \ with
 multiplicities \ $k_1, \ldots, k_s$, \ respectively.
Since
 \[
   \binom{N}{k_1}\binom{N - k_1}{k_2} \cdots
   \binom{N - k_1 - \cdots - k_{s-1}}{k_s}
   = \frac{N (N-1) \cdots (N - k_1 - k_2 - \cdots - k_s + 1)}
          {k_1! k_2! \cdots k_s!}
 \]
 is a polynomial of the variable \ $N$ \ having degree
 \ $k_1 + \cdots + k_s \leq \ell$, \ there exists
 \ $\bP = (P_1, \ldots, P_{d^\ell}) : \RR \to \RR^{d^\ell}$, \ where
 \ $P_1$, \ldots, $P_{d^\ell}$ \ are polynomials having degree at most
 \ $\ell$ \ such that
 \ $\EE\bigl((\bzeta_1 + \cdots + \bzeta_N)^{\otimes \ell}\bigr) = \bP(N)$.
\ A term of degree \ $\ell$ \ can occur only in case 
 \ $k_1 + \cdots + k_s = \ell$, \ when \ $k_1 + 2k_2 + \cdots + sk_s = \ell$
 \ implies \ $s = 1$ \ and \ $k_1 = \ell$, \ thus the corresponding term of
 degree \ $\ell$ \ is
 \ $N (N-1) \cdots (N-\ell+1) \bigl[\EE(\bzeta_1)\bigr]^{\otimes\ell}$, \ hence
 we obtain the statement.

(ii) Using the same decomposition, we have
 \begin{align*}
  \EE\bigl((\bzeta_1 + \cdots + \bzeta_N)^{\otimes\ell}\bigr)
  = \sum_{\underset{\SC 2 k_2 + 3 k_3 + \cdots + s k_s = \ell, \; k_s \ne 0}
                   {s \in \{2,\ldots,\ell\}, \; k_2, \ldots, k_s \in \ZZ_+,}}
     &\binom{N}{k_2} \binom{N - k_2}{k_3} \cdots
      \binom{N - k_2 - \cdots - k_{s-1}}{k_s} \\
  &\times \sum_{(i_1,\ldots,i_\ell) \in P_{0,k_2,\ldots,k_s}^{(N,\ell)}}
           \EE(\bzeta_{i_1} \otimes \cdots \otimes \bzeta_{i_\ell}) .
 \end{align*}
Here
 \[
   \binom{N}{k_2}\binom{N - k_2}{k_3} \cdots
   \binom{N - k_2 - \cdots - k_{s-1}}{k_s}
   = \frac{N (N-1) \cdots (N - k_2 - k_3 - \cdots - k_s + 1)}
          {k_2! k_3! \cdots k_s!}
 \]
 is a polynomial of the variable \ $N$ \ having degree \ $k_2 + \cdots + k_s$.
\ Since
 \[
   \ell = 2k_2 + 3k_3 + \cdots + sk_s \geq 2(k_2 + k_3 + \cdots + k_s) ,
 \]
 we have \ $k_2 + \cdots + k_s \leq \ell/2$ \ yielding part (ii).
\proofend

\begin{Rem}
In what follows, using the proof of Lemma \ref{LEM_moments_seged},
 we give a bit more explicit form of the polynomial \ $R_\ell$ \ in
 part (ii) of Lemma \ref{LEM_moments_seged} for the special cases
 \ $\ell = 1, 2, 3, 4, 5, 6$.
 \[
   \EE(\bzeta_1 + \cdots + \bzeta_N) = \bzero
 \]
 \[
   \EE((\bzeta_1 + \cdots + \bzeta_N)^{\otimes 2}) = N\EE(\bzeta_1^{\otimes 2}) .
 \]
 \[
   \EE((\bzeta_1 + \cdots + \bzeta_N)^{\otimes 3}) = N\EE(\bzeta_1^{\otimes 3}) .
 \] 
 \[
   \EE((\bzeta_1 + \cdots + \bzeta_N)^{\otimes 4})
   = N\EE(\bzeta_1^{\otimes 4})
     + \frac{N(N-1)}{2!}
       \sum_{(i_1,i_2,i_3,i_4) \in P_{0,2}^{(N,4)}}
        \EE(\bzeta_{i_1} \otimes \bzeta_{i_2} \otimes \bzeta_{i_3}
            \otimes \bzeta_{i_4}) .
 \]
 \[
   \EE((\bzeta_1 + \cdots + \bzeta_N)^{\otimes 5})
   = N\EE(\bzeta_1^{\otimes 5})
     + N(N-1)
       \sum_{(i_1,i_2,i_3,i_4,i_5) \in P_{0,1,1}^{(N,5)}}
        \EE(\bzeta_{i_1} \otimes \bzeta_{i_2} \otimes \bzeta_{i_3}
            \otimes \bzeta_{i_4} \otimes \bzeta_{i_5}) .
 \]
 \begin{multline*}
  \EE((\bzeta_1 + \cdots + \bzeta_N)^{\otimes 6}) \\
   \begin{aligned}
   &= N\EE(\bzeta_1^{\otimes 6})
      + N(N-1)
        \sum_{(i_1,i_2,i_3,i_4,i_5,i_6) \in P_{0,1,0,1}^{(N,6)}}
         \EE(\bzeta_{i_1} \otimes \bzeta_{i_2} \otimes \bzeta_{i_3}
             \otimes \bzeta_{i_4} \otimes \bzeta_{i_5} \otimes \bzeta_{i_6}) \\
   &\quad
      + \frac{N(N-1)}{2!}
        \sum_{(i_1,i_2,i_3,i_4,i_5,i_6) \in P_{0,0,2}^{(N,6)}}
         \EE(\bzeta_{i_1} \otimes \bzeta_{i_2} \otimes \bzeta_{i_3}
             \otimes \bzeta_{i_4} \otimes \bzeta_{i_5} \otimes \bzeta_{i_6}) \\
   &\quad
      + \frac{N(N-1)(N-2)}{3!}
        \sum_{(i_1,i_2,i_3,i_4,i_5,i_6) \in P_{0,3}^{(N,6)}}
         \EE(\bzeta_{i_1} \otimes \bzeta_{i_2} \otimes \bzeta_{i_3}
             \otimes \bzeta_{i_4} \otimes \bzeta_{i_5} \otimes \bzeta_{i_6}) .
  \end{aligned}
 \end{multline*}
\proofend
\end{Rem}

Lemma \ref{LEM_moments_seged0} can be generalized in the following way.

\begin{Lem}\label{LEM_moments_seged}
For each \ $i \in \NN$, \ let \ $(\bzeta_{i,k})_{k\in\NN}$ \ be independent and
 identically distributed random vectors with values in \ $\RR^d$ \ such that
 \ $\EE(\|\bzeta_{i,1}\|^\ell) < \infty$ \ with some \ $\ell \in \NN$.
\ Let \ $j_1, \ldots, j_\ell \in \NN$.
\begin{enumerate}
 \item[\textup{(i)}]
  Then there exists
   \ $\bQ = (Q_1, \ldots, Q_{d^\ell}) : \RR^\ell \to \RR^{d^\ell}$,
   \ where \ $Q_1$, \ldots, $Q_{d^\ell}$ \ are polynomials of \ $\ell$
   \ variables having degree at most \ $\ell-1$ \ such that
   \begin{multline*}
    \EE\bigl((\bzeta_{j_1,1} + \cdots + \bzeta_{j_1,N_1}) \otimes \cdots \otimes
             (\bzeta_{j_\ell,1} + \cdots + \bzeta_{j_\ell,N_\ell})\bigr) \\
    = N_1 \ldots N_\ell
      \EE(\bzeta_{j_1,1}) \otimes \cdots \otimes \EE(\bzeta_{j_\ell,1})
      + \bQ(N_1, \ldots, N_\ell)
   \end{multline*}
   for \ $N_1, \ldots, N_\ell \in \NN$ \ with \ $N_1 \geq \ell$, \ldots,
   \ $N_\ell \geq \ell$.
 \item[\textup{(ii)}]
  If \ $\EE(\bzeta_{j_1,1}) = \ldots = \EE(\bzeta_{j_\ell,1}) = \bzero$ \ then
   there exists \ $\bR = (R_1, \ldots, R_{d^\ell}) : \RR^\ell \to \RR^{d^\ell}$,
   \ where \ $R_1$, \ldots, $R_{d^\ell}$ \ are polynomials of \ $\ell$
   \ variables having degree at most \ $\lfloor \ell/2 \rfloor$ \ such that
   \[
     \EE\bigl((\bzeta_{j_1,1} + \cdots + \bzeta_{j_1,N_1}) \otimes \cdots \otimes
              (\bzeta_{j_\ell,1} + \cdots + \bzeta_{j_\ell,N_\ell})\bigr) \\
     = \bR(N_1, \ldots, N_\ell)
   \]
   for \ $N_1, \ldots, N_\ell \in \NN$ \ with \ $N_1 \geq \ell$, \ldots,
   \ $N_\ell \geq \ell$.
\end{enumerate}
The coefficients of the polynomials \ $\bQ$ \ and \ $\bR$ \ depend on the
 moments \ $\EE(\bzeta_{j_1,i_1} \otimes \cdots \otimes \bzeta_{j_\ell,i_\ell})$,
 \ $i_1 \in \{1, \ldots, N_1\}$, \ldots, $i_\ell \in \{1, \ldots, N_\ell\}$.
\end{Lem}

\begin{Lem}\label{LEM_moments_X}
Let \ $(\bX_k)_{k\in\ZZ_+}$ \ be a 2-type Galton--Watson process with immigration
 such that \ $\alpha, \gamma \in [0, 1)$ \ and \ $\beta, \gamma \in (0, \infty)$
 \ with \ $\alpha + \delta > 0$ \ and \ $\beta \gamma = (1- \alpha) (1 - \gamma)$
 \ (hence it is critical and positively regular).
Suppose \ $\bX_0 = \bzero$, \ and \ $\EE(\|\bxi_{1,1,1}\|^\ell) < \infty$,
 \ $\EE(\|\bxi_{1,1,2}\|^\ell) < \infty$, \ $\EE(\|\bvare_1\|^\ell) < \infty$
 \ with some \ $\ell \in \NN$.
\ Then \ $\EE(\|\bX_k\|^\ell) = \OO(k^\ell)$, \ i.e.,
 \ $\sup_{k \in \NN} k^{-\ell} \EE(\|\bX_k\|^\ell) < \infty$.
\end{Lem}

\noindent
\textbf{Proof.} \
The statement is clearly equivalent with
 \ $\EE\bigl(|P(X_{k,1},X_{k,2})|\bigr) \leq c_{P,\ell} \, k^\ell$, \ $k \in \NN$,
 \ for all polynomials \ $P$ \ of two variables having degree at most \ $\ell$,
 \ where \ $c_{P,\ell}$ \ depends only on \ $P$ \ and \ $\ell$.

If \ $\ell = 1$ \ then \eqref{EXk} and \eqref{Putzer} imply
 \begin{align*}
  \EE(\bX_k)
  = \sum_{j=0}^{k-1} \bm_\bxi^j \, \bm_\bvare 
  = \left(\frac{k}{2 - \alpha - \delta}
          \begin{bmatrix}
           1 - \delta & \beta \\
           \gamma & 1 - \alpha
          \end{bmatrix}
          + \frac{1 - (\alpha+\delta-1)^k}{(2 - \alpha - \delta)^2}
            \begin{bmatrix}
             1 - \alpha  & -\beta \\
             -\gamma & 1 - \delta
            \end{bmatrix}\right)
    \bm_\bvare
 \end{align*}
 for all \ $k \in \NN$, \ which yields the statement.

By \eqref{GWI(2)},
 \begin{align}\label{Xmoment_seged_uj}
  \begin{aligned}
   \bX_k^{\otimes 2}
   &= \left( \sum_{j=1}^{X_{k-1,1}} \bxi_{k,j,1} \right)^{\otimes 2}
      + \left( \sum_{j=1}^{X_{k-1,2}} \bxi_{k,j,2} \right)^{\otimes 2}
      + \bvare_k^{\otimes 2}
      + \left( \sum_{j=1}^{X_{k-1,1}} \bxi_{k,j,1} \right)
        \otimes \left( \sum_{j=1}^{X_{k-1,2}} \bxi_{k,j,2} \right) \\
   &\quad
      + \left( \sum_{j=1}^{X_{k-1,2}} \bxi_{k,j,2} \right)
        \otimes \left( \sum_{j=1}^{X_{k-1,1}} \bxi_{k,j,1} \right)
      + \left( \sum_{j=1}^{X_{k-1,1}} \bxi_{k,j,1} \right) \otimes \bvare_k
      + \bvare_k \otimes \left( \sum_{j=1}^{X_{k-1,1}} \bxi_{k,j,1} \right) \\
   &\quad
      + \left( \sum_{j=1}^{X_{k-1,2}} \bxi_{k,j,1} \right) \otimes \bvare_k
      + \bvare_k \otimes \left( \sum_{j=1}^{X_{k-1,2}} \bxi_{k,j,1} \right) .
  \end{aligned}
 \end{align}
Since for all \ $k \in \NN$, \ the random variables
 \ $\{\bxi_{k,j,1}, \bxi_{k,j,2}, \bvare_k : j \in \NN\}$ \ are independent of
 each other and of the $\sigma$-algebra \ $\cF_{k-1}$, \ we have
 \begin{multline*}
  \EE(\bX_k^{\otimes 2} \mid \cF_{k-1})
  = \EE\left(\left( \sum_{j=1}^M \bxi_{k,j,1} \right)^{\otimes 2}\right)
             \Bigg|_{M=X_{k-1,1}}
    + \EE\left(\left( \sum_{j=1}^N \bxi_{k,j,2} \right)^{\otimes 2}\right)
               \Bigg|_{N=X_{k-1,2}}
    + \EE(\bvare_k^{\otimes 2}) \\
  \begin{aligned}
   &\quad
      + \EE\left(\sum_{j=1}^M \bxi_{k,j,1} \right) \otimes 
        \EE\left(\sum_{j=1}^N \bxi_{k,j,2} \right)
        \Bigg|_{\overset{\SC M=X_{k-1,1}}{N=X_{k-1,2}}}
      + \EE\left(\sum_{j=1}^N \bxi_{k,j,2} \right) \otimes 
        \EE\left(\sum_{j=1}^M \bxi_{k,j,1} \right)
        \Bigg|_{\overset{\SC M=X_{k-1,1}}{N=X_{k-1,2}}} \\
   &\quad
      + \EE\left(\sum_{j=1}^M \bxi_{k,j,1} \right)\Bigg|_{M=X_{k-1,1}}
        \otimes \EE(\bvare_k)
      + \EE(\bvare_n) \otimes
        \EE\left(\sum_{j=1}^M \bxi_{k,j,1} \right)\Bigg|_{M=X_{k-1,1}} \\
   &\quad
      + \EE\left(\sum_{j=1}^N \bxi_{k,j,2} \right)\Bigg|_{N=X_{k-1,2}}
        \otimes \EE(\bvare_k)
      + \EE(\bvare_k) \otimes
        \EE\left(\sum_{j=1}^N \bxi_{k,j,2} \right)\Bigg|_{N=X_{k-1,2}} .
  \end{aligned}
 \end{multline*}
Using part (i) of Lemma \ref{LEM_moments_seged} and separating the terms having
 degree 2 and less than 2, we have
 \begin{align*}
  &\EE(\bX_k^{\otimes 2} \mid \cF_{k-1}) \\
  &= X_{k-1,1}^2 \bm_{\bxi_1}^{\otimes 2} + X_{k-1,2}^2 \bm_{\bxi_2}^{\otimes 2}
     + X_{k-1,1} X_{k-1,2}
       (\bm_{\bxi_1} \otimes \bm_{\bxi_2} + \bm_{\bxi_2} \otimes \bm_{\bxi_1})
     + \bQ_2(X_{k-1,1},X_{k-1,2}) \\
  &= (X_{k-1,1} \bm_{\bxi_1} + X_{k-1,2} \bm_{\bxi_2})^{\otimes 2}
     + \bQ_2(X_{k-1,1},X_{k-1,2})
   = (\bm_\bxi \bX_{k-1})^{\otimes 2} + \bQ_2(X_{k-1,1},X_{k-1,2}) \\
  &= \bm_\bxi^{\otimes 2} \bX_{k-1}^{\otimes 2} + \bQ_2(X_{k-1,1},X_{k-1,2}) ,
 \end{align*}
 where \ $\bQ_2 = (Q_{2,1}, Q_{2,2}, Q_{2,3}, Q_{2,4}) : \RR^2 \to \RR^4$, \ and
 \ $Q_{2,1}$, \ $Q_{2,2}$, \ $Q_{2,3}$ \ and \ $Q_{2,4}$ \ are polynomials of two
 variables having degree at most \ $1$.
\ Hence
 \[
   \EE(\bX_k^{\otimes 2})
   = \bm_\bxi^{\otimes 2} \EE(\bX_{k-1}^{\otimes 2}) + \EE[\bQ_2(X_{k-1,1},X_{k-1,2})] .
 \]
In a similar way,
 \[
   \EE(\bX_k^{\otimes \ell})
   = \bm_\bxi^{\otimes \ell} \EE(\bX_{k-1}^{\otimes \ell})
     + \EE[\bQ_\ell(X_{k-1,1},X_{k-1,2})] ,
 \]
 where \ $\bQ_\ell = (Q_{\ell,1}, \ldots, Q_{\ell,2^\ell}) : \RR^2 \to \RR^{2^\ell}$,
 \ and \ $Q_{\ell,1}$, \ldots, $Q_{\ell,2^\ell}$ \ are polynomials of two variables
 having degree at most \ $\ell-1$, \ implying
 \begin{equation}\label{EEbX_k^ell}
  \begin{aligned}
  \EE(\bX_k^{\otimes \ell})
  &= \sum_{j=1}^k (\bm_\bxi^{\otimes \ell})^{k-j} \EE[\bQ_\ell(X_{j-1,1},X_{j-1,2})] \\
  &= \sum_{j=1}^k
      (\bm_\bxi^{k-j})^{\otimes \ell} \EE[\bQ_\ell(X_{j-1,1},X_{j-1,2})] .
  \end{aligned}
 \end{equation}
Let us suppose now that the statement holds for \ $1, \ldots, \ell-1$.
\ Then
 \[
   \EE[|Q_{\ell,i}(X_{j-1,1},X_{j-1,2})|] \leq c_{Q_{\ell,i},\ell-1} k^{\ell-1} ,
   \qquad k \in \NN, \quad i \in \{1, \ldots, 2^\ell\} .
 \]
Formula \eqref{Putzer} clearly implies
 \ $\|(\bm_\bxi^i)^{\otimes \ell}\| = \OO(1)$, \ i.e.,
 \ $\sup_{i\in\ZZ_+,\,\ell\in\NN} \|(\bm_\bxi^i)^{\otimes \ell}\| < \infty$,
 \ hence we obtain the assertion for \ $\ell$.
\proofend

\begin{Cor}\label{EEX_EEU_EEV}
Let \ $(\bX_k)_{k\in\ZZ_+}$ \ be a 2-type Galton--Watson process with immigration
 such that \ $\alpha, \delta \in [0, 1)$ \ and \ $\beta, \gamma \in (0, \infty)$
 \ with \ $\alpha + \delta > 0$ \ and \ $\beta \gamma = (1 - \alpha) (1 - \gamma)$
 \ (hence it is critical and positively regular).
Suppose \ $\bX_0 = \bzero$, \ and \ $\EE(\|\bxi_{1,1,1}\|^\ell) < \infty$,
 \ $\EE(\|\bxi_{1,1,2}\|^\ell) < \infty$, \ $\EE(\|\bvare_1\|^\ell) < \infty$
 \ with some \ $\ell \in \NN$.
\ Then
 \begin{gather*}
  \EE(\|\bX_k\|^i) = \OO(k^i) , \qquad
  \EE(\bM_k^{\otimes i}) = \OO(k^{\lfloor i/2 \rfloor}) , \qquad
  \EE(U^i_k ) = \OO(k^i) , \qquad
  \EE(V^{2j}_k ) = \OO(k^j)
 \end{gather*}
 for \ $i, j \in \ZZ_+$ \ with \ $i \leq \ell$ \ and \ $2 j \leq \ell$.

If, in addition, \ $\langle \obV_{\!\!\bxi} \bv_\LEFT, \bv_\LEFT \rangle = 0$, \ then
 \[
   \EE(\langle \bv_\LEFT, \bM_k \rangle^i) = \OO(1) , \qquad
   \EE(V_k^{2j}) = \OO(1)
 \]
 for \ $i, j \in \ZZ_+$ \ with \ $i \leq \ell$ \ and \ $2 j \leq \ell$.
\end{Cor}

\noindent
\textbf{Proof.}
The first statement is just Lemma \ref{LEM_moments_X}.
Next we turn to prove \ $\EE(\bM_k^{\otimes i}) = \OO(k^{\lfloor i/2 \rfloor})$.
\ Using \eqref{Mdeco}, part (ii) of Lemma \ref{LEM_moments_seged}, and that the
 random vectors
 \ $\big\{\bxi_{k,j,1} - \EE(\bxi_{k,j,1}), \, \bxi_{k,j.2} - \EE(\bxi_{k,j,2}), \,
          \bvare_k - \EE(\bvare_k)
          : j \in \NN \big\}$
 \ are independent of each other, independent of \ $\cF_{k-1}$, \ and have
 zero mean vector, we obtain
 \[
   \EE(\bM_k^{\otimes i} \mid \cF_{k-1}) = \bR(X_{k-1,1},X_{k-1,2}) ,
 \]
 with \ $\bR = (R_1, \ldots, R_{2^i}) : \RR^2 \to \RR^{2 i}$, \ where
 \ $R_1$, \ldots, $R_{2^i}$ \ are polynomials of two variables having degree
 at most \ $i/2$.
\ Hence
 \[
   \EE(\bM_k^{\otimes i}) = \EE(\bR(X_{k-1,1},X_{k-1,2})) .
 \]
By Lemma \ref{LEM_moments_X}, we conclude
 \ $\EE(\bM_k^{\otimes i}) = \OO(k^{\lfloor i/2 \rfloor})$.

Lemma \ref{LEM_moments_X} implies
 \ $\EE(U^i_k )
    = \EE\left[\left(\frac{(1-\delta)X_{k,1} + \beta X_{k,2}}
                          {(2-\alpha-\delta)\beta}\right)^i\right]
    = \OO(k^i)$.

Next, for \ $j \in \ZZ_+$ \ with \ $2j \leq \ell$, \ we prove
 \ $\EE(V^{2j}_k) = \OO(k^{j})$ \ using induction in \ $k$.
\ By the recursion
 \ $V_k = (\alpha + \delta - 1) V_{k-1}
          + \langle \bv_\LEFT, \bM_k + \bm_\bvare \rangle$,
 \ $k \in \NN$, \ we have
 \ $\EE(V_k) = (\alpha + \delta - 1) \EE(V_{k-1})
               + \langle \bv_\LEFT, \bm_\bvare \rangle$,
 \ $k \in \NN$, \ with initial value \ $\EE(V_0) = 0$, \ hence
 \[
   \EE(V_k)
   = \langle \bv_\LEFT, \bm_\bvare \rangle
     \sum_{i=0}^{k-1} (\alpha + \delta - 1)^i,
   \qquad k \in \NN ,
 \]
 which yields \ $|\EE(V_k)| = \OO(1)$. 
\ Indeed, for all \ $k \in \NN$,
 \[
   \left| \sum_{i=0}^{k-1} (\alpha + \delta - 1)^i \right|
   \leq \frac{1}{1 - |\alpha + \delta - 1|} .
 \]
The rest of the proof of \ $\EE(V^{2j}_k) = \OO(k^{j})$ \ can be carried out as
 in Corollary 9.1 of Barczy et al.~\cite{BarIspPap2}.

By \eqref{Mdeco} and Remark \ref{REMARK0},
 \ $\langle \obV_{\!\!\bxi} \bv_\LEFT, \bv_\LEFT \rangle = 0$ \ implies
 \[
   \langle \bv_\LEFT, \bM_k \rangle
   = \langle \bv_\LEFT, \bvare_k - \EE(\bvare_k) \rangle , \qquad
   k \in \NN ,
 \]
 implying
 \ $\EE(\langle \bv_\LEFT, \bM_k \rangle^i)
    = \EE(\langle \bv_\LEFT, \bvare_1 - \EE(\bvare_1) \rangle^i) = \OO(1)$
 \ for \ $i \in \ZZ_+$ \ with \ $i \leq \ell$.

Finally, by \eqref{regr}, we obtain
 \[
   \langle \bv_\LEFT, \bX_k \rangle
   = \langle \bv_\LEFT, \bm_\bxi \bX_{k-1} \rangle
     + \langle \bv_\LEFT, \bm_\bvare \rangle
     + \langle \bv_\LEFT, \bvare_k - \EE(\bvare_k) \rangle
   = \langle \bm_\bxi^\top \! \bv_\LEFT, \bX_{k-1} \rangle
     + \langle \bv_\LEFT, \bvare_k \rangle
 \]
Using \ $\bm_\bxi^\top \! \bv_\LEFT = (\alpha + \delta - 1) \bv_\LEFT$, \ we
 conclude
 \begin{equation}\label{rec_V_mod}
  V_k = (\alpha + \delta - 1) V_{k-1} + \langle \bv_\LEFT, \bvare_k \rangle ,
  \qquad k \in \NN . 
 \end{equation}
Thus we get a recursion
 \ $\tV_k = (\alpha + \delta - 1) \tV_{k-1} + \langle \bv_\LEFT, \bM_k \rangle$,
 \ $k \in \NN$, \ for the sequence \ $\tV_k := V_k - \EE(V_k)$, \ $k \in \NN$,
 \ and rest of the proof of \ $\EE(V^{2j}_k) = \OO(1)$ \ for \ $j \in \ZZ_+$
 \ with \ $2 j \leq \ell$ \ can be carried out again by the method
 Corollary 9.1 of Barczy et al.~\cite{BarIspPap2}, applying
 \ $\EE(\langle \bv_\LEFT, \bM_k \rangle^i) = \OO(1)$ \ for \ $i \in \ZZ_+$
 \ with \ $i \leq \ell$.
\proofend

The next corollary can be derived exactly as Corollary 9.2 of
 Barczy et al.~\cite{BarIspPap2}.

\begin{Cor}\label{LEM_UV_UNIFORM}
Let \ $(\bX_k)_{k\in\ZZ_+}$ \ be a 2-type Galton--Watson process with immigration
 with offspring means \ $(\alpha, \delta) \in (0, 1)^2$ \ and
 \ $\beta \gamma = (1 - \alpha) (1 - \gamma)$ \ (hence it is critical and
 positively regular).
Suppose \ $\bX_0 = \bzero$, \ and \ $\EE(\|\bxi_{1,1,1}\|^\ell) < \infty$,
 \ $\EE(\|\bxi_{1,1,2}\|^\ell) < \infty$, \ $\EE(\|\bvare_1\|^\ell) < \infty$
 \ with some \ $\ell \in \NN$.
\ Then
 \begin{itemize}
  \item[\textup{(i)}]
   for all \ $i,j\in\ZZ_+$ \ with \ $\max\{i,j\} \leq \lfloor \ell/2 \rfloor$,
    \ and for all \ $\kappa > i + \frac{j}{2} + 1$, \ we have
    \begin{align}\label{seged_UV_UNIFORM1}
     n^{-\kappa}
     \sum_{k=1}^n\vert U_k^i V_k^j\vert
     \stoch 0
     \qquad \text{as \ $n\to\infty$,}
    \end{align}
  \item[\textup{(ii)}]
   for all \ $i,j\in\ZZ_+$ \ with \ $\max\{i,j\} \leq \ell$, \ for all \ $T>0$,
    \ and for all \ $\kappa > i + \frac{j}{2} + \frac{i+j}{\ell}$, \ we have
    \begin{align}\label{seged_UV_UNIFORM2}
     n^{-\kappa} \sup_{t\in[0,T]} \vert U_\nt^i V_\nt^j \vert \stoch 0
     \qquad \text{as \ $n\to\infty$,}
    \end{align}
  \item[\textup{(iii)}]
   for all \ $i,j\in\ZZ_+$ \ with \ $\max\{i,j\} \leq \lfloor \ell/4 \rfloor$,
    \ for all \ $T>0$, \ and \ for all
    \ $\kappa > i + \frac{j}{2} + \frac{1}{2}$, \ we have
    \begin{align}\label{seged_UV_UNIFORM4}
     n^{-\kappa} \sup_{t\in[0,T]}
     \left|\sum_{k=1}^\nt [U_k^i V_k^j - \EE(U_k^i V_k^j \mid \cF_{k-1})] \right|
     \stoch 0
     \qquad \text{as \ $n\to\infty$.}
    \end{align}
 \end{itemize}

If, in addition, \ $\langle \obV_{\!\!\bxi} \bv_\LEFT, \bv_\LEFT \rangle = 0$, \ then
 \begin{itemize}
  \item[\textup{(iv)}]
   for all \ $i,j\in\ZZ_+$ \ with \ $\max\{i,j\} \leq \lfloor \ell/2 \rfloor$,
    \ and for all \ $\kappa > i + 1$, \ we have
    \begin{align}\label{seged_UV_UNIFORM1_mod}
     n^{-\kappa}
     \sum_{k=1}^n\vert U_k^i V_k^j\vert
     \stoch 0
     \qquad \text{as \ $n\to\infty$,}
    \end{align}
  \item[\textup{(v)}]
   for all \ $i,j\in\ZZ_+$ \ with \ $\max\{i,j\} \leq \ell$, \ for all \ $T>0$,
    \ and for all \ $\kappa > i + \frac{i+j}{\ell}$, \ we have
    \begin{align}\label{seged_UV_UNIFORM2_mod}
     n^{-\kappa} \sup_{t\in[0,T]} \vert U_\nt^i V_\nt^j \vert \stoch 0
     \qquad \text{as \ $n\to\infty$,}
    \end{align}
  \item[\textup{(vi)}]
   for all \ $i,j\in\ZZ_+$ \ with \ $\max\{i,j\} \leq \lfloor \ell/4 \rfloor$,
    \ for all \ $T>0$, \ and \ for all
    \ $\kappa > i + \frac{1}{2}$, \ we have
    \begin{align}\label{seged_UV_UNIFORM4_mod}
     n^{-\kappa} \sup_{t\in[0,T]}
     \left|\sum_{k=1}^\nt [U_k^i V_k^j - \EE(U_k^i V_k^j \mid \cF_{k-1})] \right|
     \stoch 0
     \qquad \text{as \ $n\to\infty$.}
    \end{align}
  \item[\textup{(vii)}]
   for all \ $j\in\ZZ_+$ \ with \ $j \leq \lfloor \ell/2 \rfloor$, \ for all
    \ $T>0$, \ and \ for all \ $\kappa > \frac{1}{2}$, \ we have
    \begin{align}\label{seged_UV_UNIFORM4_mod_V}
     n^{-\kappa} \sup_{t\in[0,T]}
     \left|\sum_{k=1}^\nt [V_k^j - \EE(V_k^j \mid \cF_{k-1})] \right|
     \stoch 0
     \qquad \text{as \ $n\to\infty$.}
    \end{align}
 \end{itemize}
\end{Cor}

\begin{Rem}
In the special case \ $\ell=2$, $i=1$, $j=0$, \ one can improve
 \eqref{seged_UV_UNIFORM2}, namely, one can show
 \begin{equation}\label{seged_UV_UNIFORM5}
  n^{-\kappa} \sup_{t\in[0,T]} U_{\nt} \stoch 0
  \qquad \text{as \ $n\to\infty$ \ for \ $\kappa > 1$,}
 \end{equation}
 see Barczy et al.~\cite{BarIspPap2}. 
\end{Rem}

\begin{Lem}\label{LEM_moments_X_sub}
Let \ $(\bX_k)_{k\in\ZZ_+}$ \ be a 2-type Galton--Watson process with immigration such
 that \ $\alpha, \gamma \in [0, 1)$ \ and \ $\beta, \gamma \in (0, \infty)$
 \ with \ $\alpha + \delta > 0$ \ and \ $\beta \gamma < (1- \alpha) (1 - \gamma)$
 \ (hence it is subcritical and positively regular).
Suppose \ $\bX_0 = \bzero$, \ and \ $\EE(\|\bxi_{1,1,1}\|^\ell) < \infty$,
 \ $\EE(\|\bxi_{1,1,2}\|^\ell) < \infty$, \ $\EE(\|\bvare_1\|^\ell) < \infty$
 \ with some \ $\ell \in \NN$.
\ Then \ $\EE(\|\bX_k\|^\ell) = \OO(1)$, \ i.e.,
 \ $\sup_{k\in\NN} \EE(\|\bX_k\|^\ell) < \infty$.
\end{Lem}

\noindent
\textbf{Proof.} \
The statement is clearly equivalent with
 \ $\EE\bigl(|P(X_{k,1},X_{k,2})|\bigr) \leq c_{P,\ell}$, \ $k \in \NN$, \ for
 all polynomials \ $P$ \ of two variables having degree at most \ $\ell$,
 \ where \ $c_{P,\ell}$ \ depends only on \ $P$ \ and \ $\ell$.

By \eqref{EXk} and \eqref{Putzer},
 \begin{align*}
  \EE(\bX_k)
  = \sum_{j=0}^{k-1} \bm_\bxi^j \, \bm_\bvare 
  = \sum_{j=0}^{k-1} \lambda_+^j \bu_\RIGHT \bu_\LEFT^\top \bm_\bvare
    + \sum_{j=0}^{k-1} \lambda_-^j \bv_\RIGHT \bv_\LEFT^\top \bm_\bvare
 \end{align*}
 for all \ $k \in \NN$, \ which, by \ $|\lambda_-| \leq \lambda_+ < 1$, \ yields
 \ $\|\EE(\bX_k)\| = \OO(1)$, \ and hence \ $\EE(\|\bX_k\|) = \OO(1)$.
 
Let us suppose now that the statement holds for \ $1, \ldots, \ell-1$.
\ By \eqref{EEbX_k^ell},
 \[
   \EE(\bX_k^{\otimes \ell})
   = \sum_{j=1}^k (\bm_\bxi^{k-j})^{\otimes \ell} \EE[\bQ_\ell(X_{j-1,1},X_{j-1,2})] ,
 \]
 where \ $\bQ_\ell = (Q_{\ell,1}, \ldots, Q_{\ell,2^\ell}) : \RR^2 \to \RR^{2^\ell}$,
 \ and \ $Q_{\ell,1}$, \ldots, $Q_{\ell,2^\ell}$ \ are polynomials of two variables
 having degree at most \ $\ell-1$.
\ By the induction hypothesis,
 \[
   \EE[|Q_{\ell,i}(X_{j-1,1},X_{j-1,2})|] \leq c_{Q_{\ell,i},\ell-1} , \qquad
   k \in \NN, \quad i \in \{1, \ldots, 2^\ell\} . 
 \]
Formula \eqref{Putzer} clearly implies
 \ $\|(\bm_\bxi^i)^{\otimes \ell}\| = \OO(\lambda_+^{i\ell})$, \ i.e.,
 \ $\sup_{i\in\ZZ_+,\,\ell\in\NN}
     \lambda_+^{-i\ell} \|(\bm_\bxi^i)^{\otimes \ell}\| < \infty$,
 \ hence, by \ $0 < \lambda_+ < 1$, \  we obtain the assertion for \ $\ell$.
\proofend

\section{CLS estimators}
\label{section_estimators}

For each \ $n \in \NN$, \ a CLS estimator \ $\hbmbxi^{(n)}$ \ of \ $\bm_\bxi$
 \ based on a sample \ $\bX_1, \ldots, \bX_n$ \ can be obtained by minimizing
 the sum of squares
 \[
   \sum_{k=1}^n \big\| \bX_k - \EE(\bX_k \mid \cF_{k-1}) \big\|^2
   = \sum_{k=1}^n \left\| \bX_k - \bm_\bxi \bX_{k-1} - \bm_\bvare \right\|^2
 \]
 with respect to \ $\bm_\bxi$ \ over \ $\RR^{2\times2}$.
\ In what follows, we use the notation \ $\bx_0 := \bzero$.
\ For all \ $n \in \NN$, \ we define the
 function \ $Q_n : (\RR^2)^n \times \RR^{2\times2} \to \RR$ \ by
 \[
   Q_n(\bx_1, \ldots, \bx_n ; \bm_\bxi')
   := \sum_{k=1}^n \left\| \bx_k - \bm_\bxi' \bx_{k-1} - \bm_\bvare \right\|^2
 \]
 for all \ $\bm_\bxi' \in \RR^{2\times2}$ \ and
 \ $\bx_1, \ldots, \bx_n \in \RR^2$.
\ By definition, for all \ $n \in \NN$, \ a CLS estimator of \ $\bm_\bxi$ \ is
 a measurable function \ $F_n : (\RR^2)^n \to \RR^{2\times2}$ \ such that
 \begin{align*}
  Q_n(\bx_1, \ldots, \bx_n; F_n(\bx_1, \ldots, \bx_n))
  = \inf_{\bm_\bxi' \in \RR^{2\times2}} Q_n(\bx_1, \ldots, \bx_n ; \bm_\bxi')
 \end{align*}
 for all \ $\bx_1, \ldots, \bx_n \in \RR^2$.
\ Next we give the solutions of this extremum problem.

\begin{Lem}\label{CLSE1}
For each \ $n \in \NN$, \ any CLS estimator of \ $\bm_\bxi$ \ is a measurable
 function \ $F_n : (\RR^2)^n \to \RR^{2\times2}$ \ for which
 \begin{equation}\label{CLSE_alpha_beta_gamma}
  F_n(\bx_1, \ldots, \bx_n)
  = H_n(\bx_1, \ldots, \bx_n) G_n(\bx_1, \ldots, \bx_n)^{-1}
 \end{equation}
 on the set
 \[
   \bigl\{ (\bx_1, \ldots, \bx_n) \in (\RR^2)^n
           : \det(G_n(\bx_1, \ldots, \bx_n)) > 0 \bigr\} ,
 \]
 where
 \[
   G_n(\bx_1, \ldots, \bx_n) := \sum_{k=1}^n \bx_{k-1} \bx_{k-1}^\top , \qquad
   H_n(\bx_1, \ldots, \bx_n) := \sum_{k=1}^n (\bx_k - \bm_\bvare) \bx_{k-1}^\top .
 \]
\end{Lem}

\noindent
\textbf{Proof of Lemma \ref{CLSE1}.}
In the proof we write
 \[
   \bm_\bxi'
   = \begin{bmatrix} \alpha' & \beta' \\ \gamma' & \delta' \end{bmatrix} ,
   \qquad
   F_n
   = \begin{bmatrix}
      f_{1,1}^{(n)} & f_{1,2}^{(n)} \\
      f_{2,1}^{(n)} & f_{2,2}^{(n)}
     \end{bmatrix} , \qquad
   h_{n,i}
   = \sum_{k=1}^n
      \begin{bmatrix}
       (x_{k,i} - m_{\bvare,i}) x_{k-1,1} \\
       (x_{k,i} - m_{\bvare,i}) x_{k-1,2}
      \end{bmatrix} , \quad i \in \{1, 2\} ,
 \]
 $G_n(\bx_1, \ldots, \bx_n) = G_n$, \ $H_n(\bx_1, \ldots, \bx_n) = H_n$, \ and
 \ $Q_n(\bx_1, \ldots, \bx_n ; \bm_\bxi') = Q_n$. 
\ The quadratic function \ $Q_n$ \ can be written in the form
 \begin{align*}
  Q_n &= \sum_{k=1}^n
          (x_{k,1} - \alpha' x_{k-1,1} - \beta' x_{k-1,2} - m_{\bvare,1})^2
         + \sum_{k=1}^n
            (x_{k,2} - \gamma' x_{k-1,1} - \delta' x_{k-1,2} - m_{\bvare,2})^2 \\
      &= \sum_{k=1}^n
          \begin{bmatrix} \alpha' \\ \beta' \end{bmatrix}^\top
          \begin{bmatrix}
           x_{k-1,1}^2 & x_{k-1,1} x_{k-1,2} \\
           x_{k-1,1} x_{k-1,2} & x_{k-1,2}^2
          \end{bmatrix}
          \begin{bmatrix} \alpha' \\ \beta' \end{bmatrix}
         - \sum_{k=1}^n
            \begin{bmatrix} \alpha' \\ \beta' \end{bmatrix}^\top
            \begin{bmatrix}
             (x_{k,1} - m_{\bvare,1}) x_{k-1,1} \\
             (x_{k,1} - m_{\bvare,1}) x_{k-1,2}
           \end{bmatrix} \\
      &\quad 
         - \sum_{k=1}^n
            \begin{bmatrix}
             (x_{k,1} - m_{\bvare,1}) x_{k-1,1} \\
             (x_{k,1} - m_{\bvare,1}) x_{k-1,2}
            \end{bmatrix}^\top
            \begin{bmatrix} \alpha' \\ \beta' \end{bmatrix}
         + \sum_{k=1}^n (x_{k,1} - m_{\bvare,1})^2 \\
      &\quad
         + \sum_{k=1}^n
            \begin{bmatrix} \gamma' \\ \delta' \end{bmatrix}^\top
            \begin{bmatrix}
             x_{k-1,1}^2 & x_{k-1,1} x_{k-1,2} \\
             x_{k-1,1} x_{k-1,2} & x_{k-1,2}^2
            \end{bmatrix}
            \begin{bmatrix} \gamma' \\ \delta' \end{bmatrix}
         - \sum_{k=1}^n
            \begin{bmatrix} \gamma' \\ \delta' \end{bmatrix}^\top
            \begin{bmatrix}
             (x_{k,2} - m_{\bvare,2}) x_{k-1,1} \\
             (x_{k,2} - m_{\bvare,2}) x_{k-1,2}
           \end{bmatrix} \\
      &\quad 
         - \sum_{k=1}^n
            \begin{bmatrix}
             (x_{k,2} - m_{\bvare,2}) x_{k-1,1} \\
             (x_{k,2} - m_{\bvare,2}) x_{k-1,2}
            \end{bmatrix}^\top
            \begin{bmatrix} \gamma' \\ \delta' \end{bmatrix}
         + \sum_{k=1}^n (x_{k,2} - m_{\bvare,2})^2 \\
      &= \left( \begin{bmatrix} \alpha' \\ \beta' \end{bmatrix}
                - G_n^{-1} h_{n,1} \right)^\top
         G_n
         \left( \begin{bmatrix} \alpha' \\ \beta' \end{bmatrix}
                - G_n^{-1} h_{n,1} \right)
         - h_{n,1}^\top G_n^{-1} h_{n,1} \\
      &\quad
         + \left( \begin{bmatrix} \gamma' \\ \delta' \end{bmatrix}
                  - G_n^{-1} h_{n,2} \right)^\top
           G_n
           \left( \begin{bmatrix} \gamma' \\ \delta' \end{bmatrix}
                  - G_n^{-1} h_{n,2} \right)
         - h_{n,2}^\top G_n^{-1} h_{n,2}
         + \sum_{k=1}^n \| \bx_k - \bm_\bvare \|^2  .
 \end{align*}
We can check that the matrix \ $G_n$ \ is strictly positive definite.
Indeed, 
 \ $\sum_{k=1}^n x_{k-1,1}^2 \sum_{k=1}^n x_{k-1,2}^2
    - \left(\sum_{k=1}^n x_{k-1,1} x_{k-1,2}\right)^2 > 0$
 \ implies \ $\sum_{k=1}^n x_{k-1,1}^2 \sum_{k=1}^n x_{k-1,2}^2 > 0$, \ and hence,
 \ $\sum_{k=1}^n x_{k-1,1}^2 > 0$ \ and \ $\sum_{k=1}^n x_{k-1,2}^2 > 0$. 
\ Consequently, 
 \[
   \begin{bmatrix} f_{1,1}^{(n)} \\ f_{1,2}^{(n)} \end{bmatrix} = G_n^{-1} h_{n,1} ,
   \qquad
   \begin{bmatrix} f_{2,1}^{(n)} \\ f_{2,2}^{(n)} \end{bmatrix} = G_n^{-1} h_{n,2} ,
   \qquad
   F_n
   = \begin{bmatrix} h_{n,1}^\top \\ h_{n,2}^\top \end{bmatrix} (G_n^{-1})^\top
   = H_n G_n^{-1} ,
 \]
 since
 \[
   \begin{bmatrix} h_{n,1}^\top \\ h_{n,2}^\top \end{bmatrix}
   = \sum_{k=1}^n
      \begin{bmatrix}
       (x_{k,1} - m_{\bvare,1}) x_{k-1,1} & (x_{k,1} - m_{\bvare,1}) x_{k-1,2} \\
       (x_{k,2} - m_{\bvare,2}) x_{k-1,1} & (x_{k,2} - m_{\bvare,2}) x_{k-1,2}
      \end{bmatrix}
   = \sum_{k=1}^n (\bx_k - \bm_\bvare) \bx_{k-1}^\top
   = H_n ,
 \]
 hence we obtain \eqref{CLSE_alpha_beta_gamma}.
\proofend

For the existence of these CLS estimators in case of a critical symmetric
 2-type Galton--Watson process, i.e., when \ $\varrho = 1$, \ we need the
 following approximations.

\begin{Lem}\label{main_VV}
Suppose that the assumptions of Theorem \ref{main} hold.
For each \ $T > 0$, \ we have
 \[
   n^{-2}
   \sup_{t\in[0,T]}
    \left| \sum_{k=1}^\nt V_k^2
           - \frac{\langle \obV_{\!\!\bxi} \, \bv_\LEFT, \bv_\LEFT \rangle}
                  {1-\lambda^2}
             \sum_{k=1}^\nt U_{k-1} \right|
   \stoch 0 \qquad \text{as \ $n \to \infty$.}
 \]
\end{Lem}

\noindent
\textbf{Proof.}
In order to prove the satement, we derive a decomposition of
 \ $\sum_{k=1}^\nt V_k^2$ \ as a sum of a martingale and some other terms.
Using recursion \eqref{rec_V}, Lemma \ref{Moments} and \eqref{XUV}, we obtain
 \begin{align*}
  \EE(V_k^2 \mid \cF_{k-1})
  &=\EE\left[ \left. \left( \lambda V_{k-1}
                        + \langle \bv_\LEFT, \bM_k + \bm_\bvare \rangle \right)^2
              \,\right|\, \cF_{k-1} \right] \\
  &= \lambda^2 V_{k-1}^2
     + 2 \lambda \langle \bv_\LEFT, \bm_\bvare \rangle V_{k-1}
     + \langle \bv_\LEFT, \bm_\bvare \rangle^2
     + \bv_\LEFT^\top \EE(\bM_k \bM_k^\top \mid \cF_{k-1}) \bv_\LEFT \\
  &= \lambda^2 V_{k-1}^2
     + \bv_\LEFT^\top (X_{k-1,1} \bV_{\bxi_1} + X_{k-1,2} \bV_{\bxi_2}) \bv_\LEFT
     + \text{constant}
     + \text{constant $\times$ $V_{k-1}$} \\
  &= \lambda^2 V_{k-1}^2
     + \bv_\LEFT^\top \obV_{\!\!\bxi} \bv_\LEFT U_{k-1}
     + \text{constant}
     + \text{constant $\times$ $V_{k-1}$.}
 \end{align*}
Thus
 \begin{align*}
  \sum_{k=1}^\nt V_k^2
  &= \sum_{k=1}^\nt \big[ V_k^2 - \EE(V_k^2 \mid \cF_{k-1}) \big]
     + \sum_{k=1}^\nt \EE(V_k^2 \mid \cF_{k-1}) \\
  &= \sum_{k=1}^\nt \big[ V_k^2 - \EE(V_k^2 \mid \cF_{k-1}) \big]
     + \lambda^2 \sum_{k=1}^\nt V_{k-1}^2
     + \bv_\LEFT^\top \obV_{\!\!\bxi} \bv_\LEFT \sum_{k=1}^\nt U_{k-1} \\
  &\quad
     + \OO(n)
     + \text{constant $\times$ $\sum_{k=1}^\nt V_{k-1}$.}
 \end{align*}
Consequently,
 \begin{align*}
   \sum_{k=1}^\nt V_k^2
   &= \frac{1}{1-\lambda^2}
      \sum_{k=1}^\nt \big[ V_k^2 - \EE(V_k^2 \mid \cF_{k-1}) \big]
      + \frac{1}{1-\lambda^2}
        \langle \obV_{\!\!\bxi} \, \bv_\LEFT, \bv_\LEFT \rangle
        \sum_{k=1}^\nt U_{k-1} \\
   &\quad
      - \frac{\lambda^2}{1-\lambda^2} V_\nt^2
      + \OO(n)
      + \text{constant $\times$ $\sum_{k=1}^\nt V_{k-1}$.}
 \end{align*}
Using \eqref{seged_UV_UNIFORM4} with \ $(\ell, i, j) = (8, 0, 2)$ \ we obtain
 \begin{align*}
  n^{-2}
  \sup_{t\in[0,T]}
   \left|\sum_{k=1}^\nt \big[ V_k^2 - \EE(V_k^2 \mid \cF_{k-1}) \big]\right|
  \stoch 0 \qquad \text{as \ $n \to \infty$.}
 \end{align*}
Using \eqref{seged_UV_UNIFORM2} with \ $(\ell, i, j) = (3, 0, 2)$ \ we obtain
 \ $n^{-2} \sup_{t\in[0,T]} V_\nt^2 \stoch 0$.
\ Moreover, \ $n^{-2} \sum_{k=1}^\nt V_{k-1} \stoch 0$ \ as \ $n \to \infty$
 \ follows by \eqref{seged_UV_UNIFORM1} with the choice
 \ $(\ell, i, j) = (8, 0, 1)$.
\ Consequently, we obtain the statement.
\proofend

\begin{Lem}\label{limsupUV}
Suppose that the assumptions of Theorem \ref{main} hold. 
Then for each \ $T > 0$,
 \[
   n^{-5/2} \sup_{t\in[0,T]} \biggl| \sum_{k=1}^{\nt} U_{k-1} V_{k-1} \biggr| \stoch 0
   \qquad \text{as \ $n \to \infty$.}  
 \]
\end{Lem}

\noindent
\textbf{Proof.}
The aim of the following discussion is to decompose
 \ $\sum_{k=1}^{\nt} U_{k-1} V_{k-1}$ \ as a sum of a martingale and some other
 terms.
Using the recursions \eqref{rec_V}, \eqref{rec_U} and Lemma \ref{Moments}, we
 obtain
 \begin{align*}
  \EE(U_{k-1}V_{k-1} \mid \cF_{k-2})
   &= \EE\Bigl((U_{k-2} + \langle \bu_\LEFT, \bM_{k-1} + \bm_\bvare \rangle)
               \bigl(\lambda V_{k-2}
                     + \langle \bv_\LEFT, \bM_{k-1} + \bm_\bvare \rangle\bigr)
               \,\Big|\, \cF_{k-2}\Bigr) \\
   &= \lambda U_{k-2} V_{k-2}
      + \langle \bv_\LEFT, \bm_\bvare \rangle U_{k-2}
      + \lambda \langle \bu_\LEFT, \bm_\bvare \rangle V_{k-2}
      + \bu_\LEFT^\top \bm_\bvare \bm_\bvare^\top \bv_\LEFT \\
   &\quad
      + \bu^\top \EE(\bM_{k-1} \bM_{k-1}^\top \mid \cF_{k-2}) \bv \\
   &= \lambda U_{k-2} V_{k-2}
      + \text{constant}
      + \text{linear combination of \ $U_{k-2}$ \ and \ $V_{k-2}$.}
 \end{align*}
Thus
 \begin{align*}
  \sum_{k=1}^{\nt} U_{k-1} V_{k-1}
  &= \sum_{k=2}^{\nt} \big[ U_{k-1} V_{k-1} - \EE(U_{k-1}V_{k-1} \mid \cF_{k-2}) \big]
     + \sum_{k=2}^{\nt} \EE(U_{k-1}V_{k-1} \mid \cF_{k-2}) \\
  &= \sum_{k=2}^{\nt} \big[ U_{k-1} V_{k-1} - \EE(U_{k-1}V_{k-1} \mid \cF_{k-2}) \big]
     + \lambda \sum_{k=2}^{\nt} U_{k-2} V_{k-2} \\
  &\quad
     + \OO(n)
     + \text{linear combination of \ $\sum_{k=2}^{\nt} U_{k-2}$ \ and
             \ $\sum_{k=2}^{\nt} V_{k-2}$.}
 \end{align*}
Consequently
 \begin{multline*}
  \sum_{k=2}^{\nt} U_{k-1} V_{k-1}
  = \frac{1}{1-\lambda}
    \sum_{k=2}^{\nt}
     \big[ U_{k-1} V_{k-1} - \EE(U_{k-1}V_{k-1} \mid \cF_{k-2}) \big] \\
  - \frac{\lambda}{1-\lambda} U_{\nt-1} V_{\nt-1}
  + \OO(n)
  + \text{linear combination of \ $\sum_{k=2}^{\nt} U_{k-2}$ \ and
          \ $\sum_{k=2}^{\nt} V_{k-2}$.}
 \end{multline*}
Using \eqref{seged_UV_UNIFORM4} with \ $(\ell, i , j) = (4, 1, 1)$ \ we have
 \begin{align*}
    n^{-5/2}\sup_{t \in [0,T]}\,
           \Biggl\vert \sum_{k=2}^\nt
                    \big[U_{k-1} V_{k-1}
                         - \EE(U_{k-1} V_{k-1} \mid \cF_{k-2}) \big]
           \Biggr\vert
\stoch 0 \qquad \text{as \ $n\to\infty$.}
 \end{align*}
Thus, in order to show the statement, it suffices to prove
 \begin{gather}\label{5supsumU_5supsumV_5supUV}
  n^{-5/2} \sum_{k=1}^{\nT} U_k \stoch 0 , \qquad
  n^{-5/2} \sum_{k=1}^{\nT} |V_k| \stoch 0 , \qquad
  n^{-5/2} \sup_{t \in [0,T]} | U_\nt V_\nt | \stoch 0  
 \end{gather}
 as \ $n \to \infty$.
\ Using \eqref{seged_UV_UNIFORM1} with \ $(\ell, i , j) = (2, 1, 0)$ \ and
 \ $(\ell, i , j) = (2, 0, 1)$, and \eqref{seged_UV_UNIFORM2} with
 \ $(\ell, i , j) = (3, 1, 1)$ \ we have
 \eqref{5supsumU_5supsumV_5supUV}, thus we conclude the statement.
\proofend

\begin{Lem}\label{main_VVt}
Suppose that the assumptions of Theorem \ref{main} hold.
If \ $\langle \obV_{\!\!\bxi} \bv_\LEFT, \bv_\LEFT \rangle = 0$, \ then for each
 \ $T > 0$,
 \[
   n^{-1} \sup_{t\in[0,T]} \left| \sum_{k=1}^\nt V_k^2 - Mt \right| \stoch 0 \qquad
   \text{as \ $n \to \infty$,}
 \]
 where \ $M$ \ is defined in Theorem \ref{maint_Ad}.

Moreover, \ $M = 0$ \ if and only if
 \ $\langle \bV_{\!\!\bvare} \bv_\LEFT, \bv_\LEFT \rangle
    + \langle \bv_\LEFT, \bm_\bvare \rangle^2 = 0$,
 \ which is equivalent to \ $(1-\alpha) X_{k,1} \ase \beta X_{k,2}$ \ for all
 \ $k \in \NN$.
\end{Lem}

\noindent
\textbf{Proof.}
First we show
 \begin{equation}\label{sumV}
  n^{-1}
  \sup_{t\in[0,T]}
   \left|\sum_{k=1}^\nt V_k
         - \frac{\langle \bv_\LEFT, \bm_\bvare \rangle}{1-\lambda} t\right|
  \stoch 0 \qquad \text{as \ $n \to \infty$}
 \end{equation}
 for each \ $T > 0$.
\ Using recursion \eqref{rec_V_mod}, we obtain
 \[
   \EE(V_k \mid \cF_{k-1})
   = \lambda V_{k-1} + \langle \bv_\LEFT, \bm_\bvare \rangle ,
   \qquad k \in \NN . 
 \]
Thus
 \[
   \sum_{k=1}^\nt V_k
   = \sum_{k=1}^\nt [V_k - \EE(V_k \mid \cF_{k-1})]
     + \lambda \sum_{k=1}^n V_{k-1}
     + \nt \langle \bv_\LEFT, \bm_\bvare \rangle .
 \]
Consequently,
 \[
   \sum_{k=1}^\nt V_k
   = \frac{1}{1-\lambda} \sum_{k=1}^\nt [V_k - \EE(V_k \mid \cF_{k-1})]
     - \frac{\lambda}{1-\lambda} V_\nt
     + \nt \frac{\langle \bv_\LEFT, \bm_\bvare \rangle}{1-\lambda} .
 \]
Using \eqref{seged_UV_UNIFORM4_mod_V} with \ $(\ell, j) = (2, 1)$ \ we obtain
 \begin{align*}
  n^{-1}
  \sup_{t\in[0,T]}
   \left|\sum_{k=1}^\nt \big[ V_k - \EE(V_k \mid \cF_{k-1}) \big]\right|
  \stoch 0 \qquad \text{as \ $n \to \infty$.}
 \end{align*}
Using \eqref{seged_UV_UNIFORM2_mod} with \ $(\ell, i, j) = (2, 0, 1)$ \ we
 obtain \ $n^{-1} \sup_{t\in[0,T]} |V_\nt| \stoch 0$ \ as \ $n \to \infty$, \ and
 hence we conclude \eqref{sumV}.

In order to prove the convergence in the statement, we derive a decomposition
 of \ $\sum_{k=1}^\nt V_k^2$ \ as a sum of a martingale and some other terms.
Using recursion \eqref{rec_V_mod}, we obtain
 \begin{align*}
  \EE(V_k^2 \mid \cF_{k-1})
  &=\EE\left[ \left. \left( \lambda V_{k-1}
                        + \langle \bv_\LEFT, \bvare_k \rangle \right)^2
              \,\right|\, \cF_{k-1} \right] \\
  &= \lambda^2 V_{k-1}^2
     + 2 \lambda \langle \bv_\LEFT, \bm_\bvare \rangle V_{k-1}
     + \EE(\langle \bv_\LEFT, \bvare_k \rangle^2) .
 \end{align*}
Thus
 \begin{align*}
  \sum_{k=1}^\nt V_k^2
  &= \sum_{k=1}^\nt \big[ V_k^2 - \EE(V_k^2 \mid \cF_{k-1}) \big]
     + \sum_{k=1}^\nt \EE(V_k^2 \mid \cF_{k-1}) \\
  &= \sum_{k=1}^\nt \big[ V_k^2 - \EE(V_k^2 \mid \cF_{k-1}) \big]
     + \lambda^2 \sum_{k=1}^\nt V_{k-1}^2 \\
  &\quad
     + 2 \lambda \langle \bv_\LEFT, \bm_\bvare \rangle
       \sum_{k=1}^\nt V_{k-1}
     + \nt \EE(\langle \bv_\LEFT, \bvare_k \rangle^2) .
 \end{align*}
Consequently,
 \begin{align*}
   \sum_{k=1}^\nt V_k^2
   &= \frac{1}{1-\lambda^2}
      \sum_{k=1}^\nt \big[ V_k^2 - \EE(V_k^2 \mid \cF_{k-1}) \big]
      - \frac{\lambda^2}{1-\lambda^2} V_\nt^2 \\
   &\quad
      + \frac{2 \lambda \langle \bv_\LEFT, \bm_\bvare \rangle}
             {1-\lambda^2}
        \sum_{k=1}^\nt V_{k-1}
      + \nt \frac{\EE(\langle \bv_\LEFT, \bvare_k \rangle^2)}
                 {1-\lambda^2} .
 \end{align*}
Using \eqref{seged_UV_UNIFORM4_mod_V} with \ $(\ell, j) = (4, 2)$ \ we obtain
 \begin{align*}
  n^{-1}
  \sup_{t\in[0,T]}
   \left|\sum_{k=1}^\nt \big[ V_k^2 - \EE(V_k^2 \mid \cF_{k-1}) \big]\right|
  \stoch 0 \qquad \text{as \ $n \to \infty$.}
 \end{align*}
Using \eqref{seged_UV_UNIFORM2_mod} with \ $(\ell, i, j) = (3, 0, 2)$ \ we
 obtain \ $n^{-1} \sup_{t\in[0,T]} V_\nt^2 \stoch 0$ \ as \ $n \to \infty$.
\ Consequently, by \eqref{sumV}, we obtain
 \begin{align*}
  &n^{-1}
   \sup_{t\in[0,T]}
    \left|\sum_{k=1}^\nt V_k^2
          - \frac{2 \lambda \langle \bv_\LEFT, \bm_\bvare \rangle^2}
                 {(1-\lambda)(1-\lambda^2)} t
          - \frac{\EE(\langle \bv_\LEFT, \bvare_k \rangle^2)}
                 {1-\lambda^2} t\right|\\
  &=n^{-1}
    \sup_{t\in[0,T]}
     \left|\sum_{k=1}^\nt V_k^2 - M t\right|
    \stoch 0
    \qquad \text{as \ $n \to \infty$.}
 \end{align*}

Clearly \ $M = 0$ \ if and only if
 \ $\langle \bV_{\!\!\bvare} \bv_\LEFT, \bv_\LEFT \rangle = 0$ \ and
 \ $\langle \bv_\LEFT, \bm_\bvare \rangle = 0$, \ which is equivalent to
 \ $\EE(\langle \bv_\LEFT, \bvare_i \rangle^2)
    = \langle \bV_{\!\!\bvare} \bv_\LEFT, \bv_\LEFT \rangle
      + \langle \bv_\LEFT, \bm_\bvare \rangle^2 = 0$
 \ for all \ $i \in \NN$, \ which is equivalent to
 \ $\langle \bv_\LEFT, \bvare_k \rangle
    = (1 - \alpha) \vare_{k,1} - \beta \vare_{k,2} \ase 0$ \ for all
 \ $k \in \NN$, \ which is equivalent to
 \ $(1 - \alpha) X_{k,1} - \beta X_{k,2} \ase 0$ \ for all \ $k \in \NN$.
\proofend

\begin{Lem}\label{limsupUVt}
Suppose that the assumptions of Theorem \ref{main} hold. 
If \ $\langle \obV_{\!\!\bxi} \bv_\LEFT, \bv_\LEFT \rangle = 0$, \ then for each
 \ $T > 0$,
 \[
   n^{-2}
   \sup_{t\in[0,T]}
    \biggl| \sum_{k=1}^{\nt} U_{k-1} V_{k-1}
            - \frac{\langle\bv_\LEFT,\bm_\bvare\rangle}{1-\lambda}
              \sum_{k=1}^{\nt} U_{k-2} \biggr|
   \stoch 0 \qquad \text{as \ $n \to \infty$.}  
 \]
\end{Lem}

\noindent
\textbf{Proof.}
The aim of the following discussion is to decompose
 \ $\sum_{k=1}^{\nt} U_{k-1} V_{k-1}$ \ as a sum of a martingale and some other
 terms.
Using the recursions \eqref{rec_U}, \eqref{rec_V_mod}, and Lemma
 \ref{Moments}, we obtain
 \begin{align*}
  \EE(U_{k-1} V_{k-1} \mid \cF_{k-2})
   &= \EE\bigl((U_{k-2} + \langle \bu_\LEFT, \bM_{k-1} + \bm_\bvare \rangle)
               (\lambda V_{k-2}
                     + \langle \bv_\LEFT, \bvare_{k-1} \rangle)
               \,\big|\, \cF_{k-2}\bigr) \\
   &= \lambda U_{k-2} V_{k-2}
      + \langle \bv_\LEFT, \bm_\bvare \rangle U_{k-2}
      + \lambda \langle \bu_\LEFT, \bm_\bvare \rangle V_{k-2} \\
   &\quad
      + \langle \bu_\LEFT, \bm_\bvare \rangle 
        \langle \bv_\LEFT, \bm_\bvare \rangle
      + \bu_\LEFT^\top \EE(\bM_{k-1} \bvare_{k-1}^\top \mid \cF_{k-2}) \bv_\LEFT \\
   &= \lambda U_{k-2} V_{k-2}
      + \langle \bv_\LEFT, \bm_\bvare \rangle U_{k-2}
      + \text{constant}
      + \text{constant$\times$$V_{k-2}$.}
 \end{align*}
Thus
 \begin{align*}
  \sum_{k=1}^{\nt} U_{k-1} V_{k-1}
  &= \sum_{k=2}^{\nt} \big[ U_{k-1} V_{k-1} - \EE(U_{k-1}V_{k-1} \mid \cF_{k-2}) \big]
     + \sum_{k=2}^{\nt} \EE(U_{k-1}V_{k-1} \mid \cF_{k-2}) \\
  &= \sum_{k=2}^{\nt} \big[ U_{k-1} V_{k-1} - \EE(U_{k-1}V_{k-1} \mid \cF_{k-2}) \big]
     + \lambda \sum_{k=2}^{\nt} U_{k-2} V_{k-2} \\
  &\quad
     + \langle \bv_\LEFT, \bm_\bvare \rangle \sum_{k=2}^{\nt} U_{k-2}
     + \OO(n)
     + \text{constant$\times$$\sum_{k=2}^{\nt} V_{k-2}$.}
 \end{align*}
Consequently,
 \begin{align*}
  \sum_{k=2}^{\nt} U_{k-1} V_{k-1}
  &= \frac{1}{1-\lambda}
     \sum_{k=2}^{\nt}
      \big[ U_{k-1} V_{k-1} - \EE(U_{k-1}V_{k-1} \mid \cF_{k-2}) \big]
     - \frac{\lambda}{1-\lambda} U_{\nt-1} V_{\nt-1} \\
  &\quad
     + \frac{\langle\bv_\LEFT,\bm_\bvare\rangle}{1-\lambda} \sum_{k=2}^{\nt} U_{k-2}
     + \OO(n)
     + \text{constant$\times$$\sum_{k=2}^{\nt} V_{k-2}$.}
 \end{align*}
Using \eqref{seged_UV_UNIFORM4_mod} with \ $(\ell, i , j) = (4, 1, 1)$ \ we
 have
 \begin{align*}
  n^{-2} \sup_{t \in [0,T]} \,
        \Biggl| \sum_{k=2}^\nt
                 \big[U_{k-1} V_{k-1}
                      - \EE(U_{k-1} V_{k-1} \mid \cF_{k-2}) \big] \Biggr|
  \stoch 0 \qquad \text{as \ $n\to\infty$.}
 \end{align*}
Thus, in order to show the statement, it suffices to prove
 \begin{gather}\label{5supsumV_5supUV}
  n^{-2} \sum_{k=1}^{\nT} |V_k| \stoch 0 , \qquad
  n^{-2} \sup_{t \in [0,T]} | U_\nt V_\nt | \stoch 0  
 \end{gather}
 as \ $n \to \infty$.
\ Using \eqref{seged_UV_UNIFORM1_mod} with \ $(\ell, i , j) = (2, 0, 1)$,
 \ and \eqref{seged_UV_UNIFORM2_mod} with \ $(\ell, i , j) = (3, 1, 1)$ \ we
 have \eqref{5supsumV_5supUV}, thus we conclude the statement.
\proofend

Now we can prove asymptotic existence and uniqueness of CLS estimators of the
 offspring mean matrix and of the criticality parameter in the critical case.

\begin{Pro}\label{ExUn}
Suppose that the assumptions of Theorem \ref{main} hold, and
 \ $\langle \obV_{\!\!\bxi} \bv_\LEFT, \bv_\LEFT \rangle
    + \langle \bV_{\!\!\bvare} \bv_\LEFT, \bv_\LEFT \rangle
    + \langle \bv_\LEFT, \bm_\bvare \rangle^2 > 0$.
\ Then \ $\lim_{n \to \infty} \PP(\Omega_n) = 1$, \ where \ $\Omega_n$ \ is
 defined in \eqref{H_n}, and hence the probability of the existence of a
 unique CLS estimator \ $\hbmbxi^{(n)}$ \ converges to 1 as \ $n \to \infty$.
\ If \ $\langle \obV_{\!\!\bxi} \bv_\LEFT, \bv_\LEFT \rangle > 0$ \ then
\ $\lim_{n \to \infty} \PP(\tOmega_n) = 1$, \ where \ $\tOmega_n$ \ is
 defined in \eqref{tH_n}, and hence the probability of the existence of the
 estimator \ $\hvarrho_n$ \ converges to 1 as \ $n \to \infty$.
\end{Pro}

\noindent
\textbf{Proof.}
Recall convergence \ $(n^{-1} U_{\nt})_{t\in\RR_+} \distr (\cY_t)_{t\in\RR_+}$ \ from
 \eqref{convU}.
Using Lemmas \ref{main_VV}, \ref{limsupUV}, \ref{Conv2Funct} and \ref{Marci},
 one can show
 \[
   \sum_{k=1}^n
    \begin{bmatrix}
     n^{-3} U_{k-1}^2 \\
     n^{-5/2} U_{k-1} V_{k-1} \\
     n^{-2} V_{k-1}^2
    \end{bmatrix}
   \distr
   \begin{bmatrix}
    \int_0^1 \cY_t^2 \, \dd t \\[1mm]
    0 \\[1mm]
    \frac{\langle\bv_\LEFT,\bm_\bvare\rangle}{1-\lambda} \,
    \int_0^1 \cY_t \, \dd t  \end{bmatrix}
   \qquad \text{as \ $n \to \infty$.}
 \]
By \eqref{A_n_deco} and continuous mapping theorem,
 \begin{equation}\label{seged2}
  n^{-5} \det(\bA_n)
  \distr
  \frac{\langle\obV_{\!\!\bxi}\bv_\LEFT,\bv_\LEFT\rangle}{1-\lambda^2}
  \int_0^1 \cY_t^2 \, \dd t
  \int_0^1 \cY_t \, \dd t
  \qquad \text{as \ $n \to \infty$.}
 \end{equation}
Since \ $\bm_\bvare \ne \bzero$, \ by the SDE \eqref{SDE_Y}, we have
 \ $\PP(\text{$\cY_t = 0$ \ for all \ $t \in [0,1]$}) = 0$, \ which implies
 that
 \ $\PP\bigl( \int_0^1 \cY_t^2 \, \dd t \int_0^1 \cY_t \, \dd t > 0 \bigr) = 1$.
\ Consequently, the distribution function of
 \ $\int_0^1 \cY_t^2 \, \dd t \int_0^1 \cY_t \, \dd t$ \ is continuous at 0.

If \ $\langle \obV_{\!\!\bxi} \bv_\LEFT, \bv_\LEFT \rangle > 0$ \ then, by
 \eqref{seged2},
 \begin{align*}
  \PP(\Omega_n)
  &= \PP\left( \det(\bA_n) > 0 \right)
   = \PP\left( n^{-5} \det(\bA_n) > 0 \right) \\
  &\to
   \PP\left(\frac{\langle\obV_{\!\!\bxi}\bv_\LEFT,\bv_\LEFT\rangle}
                 {1-\lambda^2}
            \int_0^1 \cY_t^2 \, \dd t \int_0^1 \cY_t \, \dd t > 0 \right)
   = \PP\left( \int_0^1 \cY_t^2 \, \dd t \int_0^1 \cY_t \, \dd t > 0 \right)
   = 1
 \end{align*}
 as \ $n \to \infty$.

If \ $\langle \obV_{\!\!\bxi} \bv_\LEFT, \bv_\LEFT \rangle = 0$, \ then, using
 convergence \eqref{convU} and Lemmas \ref{main_VVt}, \ref{limsupUVt},
 \ref{Conv2Funct} and \ref{Marci}, one can show
 \[
   \sum_{k=1}^n
    \begin{bmatrix}
     n^{-3} U_{k-1}^2 \\
     n^{-2} U_{k-1} V_{k-1} \\
     n^{-1} V_{k-1}^2
    \end{bmatrix}
   \distr
   \begin{bmatrix}
    \int_0^1 \cY_t^2 \, \dd t \\[1mm]
    \frac{\langle\bv_\LEFT,\bm_\bvare\rangle}{1-\lambda}
    \int_0^1 \cY_t \, \dd t \\[1mm]
    \frac{\langle\bv_\LEFT,\bm_\bvare\rangle^2}{(1-\lambda)^2}
    + \frac{\langle\bV_{\!\!\bvare}\bv_\LEFT,\bm_\bvare\rangle}{1-\lambda^2}
   \end{bmatrix}
   \qquad \text{as \ $n \to \infty$.}
 \]
By \eqref{A_n_deco} and continuous mapping theorem,
 \begin{align*}
  n^{-4} \det(\bA_n)
  &\distr
   \left(\frac{\langle\bv_\LEFT,\bm_\bvare\rangle^2}{(1-\lambda)^2}
         + \frac{\langle\bV_{\!\!\bvare}\bv_\LEFT,\bm_\bvare\rangle}
                {1-\lambda^2}\right)
           \int_0^1 \cY_t^2 \, \dd t
         - \left(\frac{\langle\bv_\LEFT,\bm_\bvare\rangle}{1-\lambda}
                 \int_0^1 \cY_t \, \dd t\right)^2 \\
  &= \frac{\langle\bV_{\!\!\bvare}\bv_\LEFT,\bm_\bvare\rangle}{1-\lambda^2}
     \int_0^1 \cY_t^2 \, \dd t
     + \frac{\langle\bv_\LEFT,\bm_\bvare\rangle^2}{(1-\lambda)^2}
       \left( \int_0^1 \cY_t^2 \, \dd t
              - \left(\int_0^1 \cY_t \, \dd t\right)^2 \right)
 \end{align*}
 as \ $n \to \infty$.
\ As above, \ $\PP\bigl( \int_0^1 \cY_t^2 \, \dd t > 0 \bigr) = 1$.
\ It is also known that
 \ $\PP\bigl( \int_0^1 \cY_t^2 \, \dd t
              - \bigl(\int_0^1 \cY_t \, \dd t\bigr)^2 > 0 \bigr) = 1$,
 \ see, e.g., the proof of Theorem 3.2 in Barczy et al.\ \cite{BarDorLiPap}.
Consequently, the distribution function of the above limit distribution is
 continuous at 0.
Since
 \ $\langle \obV_{\!\!\bxi} \bv_\LEFT, \bv_\LEFT \rangle
    + \langle \bV_{\!\!\bvare} \bv_\LEFT, \bv_\LEFT \rangle
    + \langle \bv_\LEFT, \bm_\bvare \rangle^2 > 0$
 \ and \ $\langle \obV_{\!\!\bxi} \bv_\LEFT, \bv_\LEFT \rangle = 0$, \ we have
 \ $\langle \bV_{\!\!\bvare} \bv_\LEFT, \bv_\LEFT \rangle > 0$ \ or
 \ $\langle \bv_\LEFT, \bm_\bvare \rangle^2 > 0$,
 \ hence,
 \begin{align*}
  &\PP(\Omega_n)
   = \PP\left( \det(\bA_n) > 0 \right)
   = \PP\left( n^{-4} \det(\bA_n) > 0 \right) \\
  &\to
   \PP\left( \frac{\langle\bV_{\!\!\bvare}\bv_\LEFT,\bm_\bvare\rangle}
                  {1-\lambda^2}
             \int_0^1 \cY_t^2 \, \dd t
             + \frac{\langle\bv_\LEFT,\bm_\bvare\rangle^2}{(1-\lambda)^2}
               \left( \int_0^1 \cY_t^2 \, \dd t
                      - \left(\int_0^1 \cY_t \, \dd t\right)^2 \right)
             > 0 \right)
   = 1
 \end{align*}
 as \ $n \to \infty$.

If \ $\langle \obV_{\!\!\bxi} \bv_\LEFT, \bv_\LEFT \rangle > 0$, \ then
 \eqref{m_xi} yields \ $\hbmbxi^{(n)} \distr \bm_\bxi$ \ as \ $n \to \infty$,
 \ and hence \ $\hbmbxi^{(n)} \stoch \bm_\bxi$ \ as \ $n \to \infty$, \ thus
 \ $(\halpha_n - \hdelta_n)^2 + 4 \hbeta_n \hgamma_n
    \stoch (\alpha - \delta)^2 + 4 \beta \gamma = (1 - \lambda)^2 > 0$,
 \ implying
 \[
   \PP(\tOmega_n)
   = \PP\bigl((\halpha_n - \hdelta_n)^2 + 4 \hbeta_n \hgamma_n \geq 0\bigr)
   \to 1 \qquad \text{as \ $n \to \infty$,}
 \] 
 hence we obtain the satement.
\proofend

Next we prove asymptotic existence and uniqueness of CLS estimators of the
 offspring mean matrix and of the criticality parameter in the subcritical case.

\begin{Pro}\label{ExUn_sub}
Suppose that the assumptions of Theorem \ref{main_sub} hold.
Then \ $\lim_{n \to \infty} \PP(\Omega_n) = 1$ \ and
 \ $\lim_{n \to \infty} \PP(\tOmega_n) = 1$, \ where \ $\Omega_n$ \ and \ $\tOmega_n$
 \ are defined in \eqref{H_n} and \eqref{tH_n}, respectively, and hence the probability
 of the existence of unique CLS estimators \ $\hbmbxi^{(n)}$ \ and \ $\hvarrho_n$
 \ converges to 1 as \ $n \to \infty$.
\end{Pro}

\noindent
\textbf{Proof.}
Recall convergence
 \[
   n^{-1} \bA_n \as \EE\Bigl(\tbX \tbX^\top\Bigr) , \qquad \text{as \ $n \to \infty$,}
 \]
 from \eqref{SLLN1}.
Consequently, 
 \[
   n^{-1} \det(\bA_n) \as \det\Bigl(\EE\Bigl(\tbX \tbX^\top\Bigr)\Bigr) , \qquad 
   \text{as \ $n \to \infty$.}
 \]
The matrix \ $\EE\bigl(\tbX \tbX^\top\bigr)$ \ is invertible, see the proof of Theorem
 \ref{main_sub}.
Thus \ $\det\bigl(\EE\bigl(\tbX \tbX^\top\bigr)\bigr) > 0$, \ and hence we obtain 
 \[
   \PP(\Omega_n)
   = \PP(\det(\bA_n) > 0) = \PP(n^{-1} \det(\bA_n) > 0) \to 1 , \qquad 
   \text{as \ $n \to \infty$.} 
 \]
The estimator \ $\hbmbxi^{(n)}$ \ is strongly consistent, hence
 \ $\hbmbxi^{(n)} \as \bm_\bxi$ \ as \ $n \to \infty$.
\ Consequently, 
 \[
   (\halpha_n - \hdelta_n)^2 + 4 \hbeta_n \hgamma_n
   \as (\alpha - \delta)^2 + 4 \beta \gamma > 0 , \qquad \text{as \ $n \to \infty$,}
 \]
 and hence we obtain \ $\lim_{n \to \infty} \PP(\tOmega_n) = 1$.
\proofend

\section{A version of the continuous mapping theorem}
\label{app_C}

The following version of continuous mapping theorem can be found for example
 in Kallenberg \cite[Theorem 3.27]{K}.

\begin{Lem}\label{Lem_Kallenberg}
Let \ $(S, d_S)$ \ and \ $(T, d_T)$ \ be metric spaces and
 \ $(\xi_n)_{n \in \NN}$, \ $\xi$ \ be random elements with values in \ $S$
 \ such that \ $\xi_n \distr \xi$ \ as \ $n \to \infty$.
\ Let \ $f : S \to T$ \ and \ $f_n : S \to T$, \ $n \in \NN$, \ be measurable
 mappings and \ $C \in \cB(S)$ \ such that \ $\PP(\xi \in C) = 1$ \ and
 \ $\lim_{n \to \infty} d_T(f_n(s_n), f(s)) = 0$ \ if
 \ $\lim_{n \to \infty} d_S(s_n,s) = 0$ \ and \ $s \in C$.
\ Then \ $f_n(\xi_n) \distr f(\xi)$ \ as \ $n \to \infty$.
\end{Lem}

For the case \ $S = \DD(\RR_+, \RR^d)$ \ and \ $T = \RR^q$
 \ (or \ $T = \DD(\RR_+,\RR^q)$), \ where \ $d$, \ $q \in \NN$, \ we formulate
 a consequence of Lemma \ref{Lem_Kallenberg}.

For functions \ $f$ \ and \ $f_n$, \ $n \in \NN$, \ in \ $\DD(\RR_+, \RR^d)$,
 \ we write \ $f_n \lu f$ \ if \ $(f_n)_{n \in \NN}$ \ converges to \ $f$
 \ locally uniformly, i.e., if \ $\sup_{t \in [0,T]} \|f_n(t) - f(t)\| \to 0$ \ as
 \ $n \to \infty$ \ for all \ $T > 0$.
\ For measurable mappings
 \ $\Phi : \DD(\RR_+, \RR^d) \to \RR^q$
 \ (or \ $\Phi : \DD(\RR_+, \RR^d) \to \DD(\RR_+,\RR^q)$) \ and
 \ $\Phi_n : \DD(\RR_+, \RR^d) \to \RR^q$
 \ (or \ $\Phi_n : \DD(\RR_+, \RR^d) \to \DD(\RR_+,\RR^q)$), \ $n \in \NN$,
 \ we will denote by \ $C_{\Phi, (\Phi_n)_{n \in \NN}}$ \ the set of all functions
 \ $f \in \CC(\RR_+, \RR^d)$ \ such that
 \ $\Phi_n(f_n) \to \Phi(f)$ \ (or \ $\Phi_n(f_n) \lu \Phi(f)$) \ whenever
 \ $f_n \lu f$ \ with \ $f_n \in \DD(\RR_+, \RR^d)$, \ $n \in \NN$.

We will use the following version of the continuous mapping theorem several
 times, see, e.g., Barczy et al. \cite[Lemma 4.2]{BarIspPap1} and Isp\'any and
 Pap \cite[Lemma 3.1]{IspPap}.

\begin{Lem}\label{Conv2Funct}
Let \ $d, q \in \NN$, \ and \ $(\bcU_t)_{t\in\RR_+}$ \ and
 \ $(\bcU^{(n)}_t)_{t\in\RR_+}$, \ $n \in \NN$, \ be \ $\RR^d$-valued stochastic
 processes with c\`adl\`ag paths such that \ $\bcU^{(n)} \distr \bcU$.
\ Let \ $\Phi : \DD(\RR_+, \RR^d) \to \RR^q$
 \ (or \ $\Phi : \DD(\RR_+, \RR^d) \to \DD(\RR_+,\RR^q)$) \ and
 \ $\Phi_n : \DD(\RR_+, \RR^d) \to \RR^q$
 \ (or \ $\Phi_n : \DD(\RR_+, \RR^d) \to \DD(\RR_+,\RR^q)$), \ $n \in \NN$,
 \ be measurable mappings such that there exists
 \ $C \subset C_{\Phi,(\Phi_n)_{n\in\NN}}$ \ with \ $C \in \cB(\DD(\RR_+, \RR^d))$
 \ and \ $\PP(\bcU \in C) = 1$.
\ Then \ $\Phi_n(\bcU^{(n)}) \distr \Phi(\bcU)$.
\end{Lem}

In order to apply Lemma \ref{Conv2Funct}, we will use the following statement
 several times, see Barczy et al.~\cite[Lemma B.3]{BarIspPap2}.

\begin{Lem}\label{Marci}
Let \ $d, p, q \in \NN$, \ $h : \RR^d \to \RR^q$ \ be a continuous function
 and \ $K : [0,1] \times \RR^{2d} \to \RR^p$ \ be a function such that for all
 \ $R > 0$ \ there exists \ $C_R > 0$ such that
 \begin{equation}\label{Lipschitz}
  \| K(s, x) - K(t, y) \| \leq C_R \left( | t - s | + \| x - y \| \right)
 \end{equation}
 for all \ $s, t \in [0, 1]$ \ and \ $x, y \in \RR^{2d}$ \ with
 \ $\| x \| \leq R$ \ and \ $\| y \| \leq R$.
\ Moreover, let us define the mappings
 \ $\Phi, \Phi_n : \DD(\RR_+, \RR^d) \to \RR^{q+p}$, \ $n \in \NN$, \ by
 \begin{align*}
  \Phi_n(f)
  &:= \left( h(f(1)),
             \frac{1}{n}
             \sum_{k=1}^n
              K\left( \frac{k}{n}, f\left( \frac{k}{n} \right),
                      f\left( \frac{k-1}{n} \right) \right) \right) , \\
  \Phi(f)
  & := \left(  h(f(1)), \int_0^1 K( u, f(u), f(u) ) \, \dd u \right)
 \end{align*}
 for all \ $f \in \DD(\RR_+, \RR^d)$.
\ Then the mappings \ $\Phi$ \ and \ $\Phi_n$, \ $n \in \NN$, \ are measurable,
 and \ $C_{\Phi,(\Phi_n)_{n \in \NN}} = \CC(\RR_+, \RR^d) \in \cB(\DD(\RR_+, \RR^d))$.
\end{Lem}

\section{Convergence of random step processes}
\label{section_conv_step_processes}

We recall a result about convergence of random step processes towards a
 diffusion process, see Isp\'any and Pap \cite{IspPap}.
This result is used for the proof of convergence \eqref{conv_Z}.

\begin{Thm}\label{Conv2DiffThm}
Let \ $\bgamma : \RR_+ \times \RR^d \to \RR^{d \times r}$ \ be a continuous
 function.
Assume that uniqueness in the sense of probability law holds for the SDE
 \begin{equation}\label{SDE}
  \dd \, \bcU_t
  = \gamma (t, \bcU_t) \, \dd \bcW_t ,
  \qquad t \in \RR_+,
 \end{equation}
 with initial value \ $\bcU_0 = \bu_0$ \ for all \ $\bu_0 \in \RR^d$, \ where
 \ $(\bcW_t)_{t \in \RR_+}$ \ is an $r$-dimensional standard Wiener process.
Let \ $(\bcU_t)_{t \in \RR_+}$ \ be a solution of \eqref{SDE} with initial value
 \ $\bcU_0 = \bzero \in \RR^d$.

For each \ $n \in \NN$, \ let \ $(\bU^{(n)}_k)_{k \in \NN}$ \ be a sequence of
 $d$-dimensional martingale differences with respect to a filtration
 \ $(\cF^{(n)}_k)_{k \in \ZZ_+}$, i.e., \ $\EE(\bU^{(n)}_k \mid \cF^{(n)}_{k-1}) = 0$,
 \ $n \in \NN$, \ $k \in \NN$.
\ Let
 \[
   \bcU^{(n)}_t := \sum_{k=1}^{\nt} \bU^{(n)}_k \, ,
   \qquad t \in \RR_+, \quad n \in \NN .
 \]
Suppose that \ $\EE \big( \|\bU^{(n)}_k\|^2 \big) < \infty$ \ for all
 \ $n, k \in \NN$.
\ Suppose that for each \ $T > 0$,
 \begin{enumerate}
  \item [\textup{(i)}]
        $\sup\limits_{t\in[0,T]}
         \left\| \sum\limits_{k=1}^{\nt}
                  {\var\bigl(\bU^{(n)}_k \mid \cF^{(n)}_{k-1}\bigr)}
                 - \int_0^t
                    \bgamma(s,\bcU^{(n)}_s) \bgamma(s,\bcU^{(n)}_s)^\top
                    \dd s \right\|
         \stoch 0$,\\
  \item [\textup{(ii)}]
        $\sum\limits_{k=1}^{\lfloor nT \rfloor}
          \EE \big( \|\bU^{(n)}_k\|^2 \bbone_{\{\|\bU^{(n)}_k\| > \theta\}}
                    \bmid \cF^{(n)}_{k-1} \big)
         \stoch 0$
        \ for all \ $\theta>0$,
 \end{enumerate}
 where \ $\stoch$ \ denotes convergence in probability.
Then \ $\bcU^{(n)} \distr \bcU$ \ as \ $n\to\infty$.
\end{Thm}

Note that in (i) of Theorem \ref{Conv2DiffThm}, \ $\|\cdot\|$ \ denotes
 a matrix norm, while in (ii) it denotes a vector norm.

\end{document}